\documentclass[10pt]{article} 

\usepackage{amsmath}
\usepackage{amssymb} 
\usepackage{latexsym} 
\usepackage{epic,eepic,epsfig,psfrag} 
\usepackage{amscd}  
\usepackage{index}

\newcommand{\QED}{\hspace*{\fill}$\Box$\medskip}

\def\rd{{\rm d}}
\def\rT{{\rm T}}

\def\rG{{\rm G}}
\def\p{\phi}
\def\ph{\varphi}
\def\a{\alpha}
\def\b{\beta} 
\def\d{\delta} 

\def\ep{\varepsilon} 
\def\e{\eta} 
\def\g{\gamma}
\def\th{\theta} 

\def\s{\sigma} 
\def\t{\tau} 
\def\i{\iota} 
\def\x{\xi} 
 
\def\n{\nu} 
 
\def\z{\zeta} 
\def\o{\omega} 
\def\l{\lambda} 
\def\D{\Delta} 
\def\L{\Lambda} 
\def\G{\Gamma} 
\def\pg{\mathhexbox278}
\def\S{\Sigma}
\def\Si{\Sigma} 
\def\Th{\Theta} 
 
\def\Om{\Omega} 
\def\P{\Phi} 
\def\cg{\mathfrak{g}}
\def\cA{{\mathcal A}}
 
\def\cC{{\mathcal C}}

\def\cG{{\mathcal G}}

\def\cK{{\mathcal K}}
\def\cL{{\mathcal L}} 
\def\cM{{\mathcal M}}

\def\cS{{\mathcal S}}
\def\cT{{\mathcal T}}
\def\cU{{\mathcal U}} 
\def\cV{{\mathcal V}}
\def\cW{{\mathcal W}}

\def\R{{\mathbb R}}

\def\N{{\mathbb N}} 
 
\def\H{{\mathbb H}} 
\def\Z{{\mathbb Z}}
 
\def\half{{\textstyle{\frac 12}}} 
\def\im{{\rm im}\,}
\def\laplace{\Delta} 
 
\def\st{\: \big| \:} 
\DeclareMathOperator{\supp}{supp}

\def\dt{{\rm d}t}
\def\ds{{\rm d}s}
\def\pd{\partial}

\def\comp{\circ}

\def\ta{{\tilde\alpha}}
\def\tg{{\tilde g}}

\def\tA{{\tilde A}}

\def\tM{{\tilde M}}
\def\tU{{\tilde U}}
\def\la{\langle\,}
\def\ra{\,\rangle}
\def\dvol{\,\rd{\rm vol}}

\def\HF{{\rm HF}}

\def\SU{{\rm SU}}

\def\tsum{\textstyle\sum}

\newtheorem{dfn}{Definition}[section] 
\newtheorem{lem}[dfn]{Lemma} 
\newtheorem{prp}[dfn]{Proposition} 
\newtheorem{thm}[dfn]{Theorem} 
\newtheorem{rmk}[dfn]{Remark} 
 
\newtheorem{ex}[dfn]{Example}

\begin{document}

\bibliographystyle{plain}

\author{Katrin Wehrheim 
\thanks{wehrheim@math.ethz.ch; supported by Swiss National Science Foundation grant 21-64937.01;
} }

\title{Anti-self-dual instantons with Lagrangian boundary conditions}
%
%

\maketitle

\begin{abstract}
We study a nonlocal boundary value problem for anti-self-dual instantons on
$4$-manifolds with a space-time splitting of the boundary. The model case is $\R\times Y$, 
where $Y$ is a compact oriented $3$-manifold with boundary $\Sigma$.
The restriction of the instanton to each time slice $\{t\}\times\Sigma$ is
required to lie in a fixed (singular) Lagrangian submanifold of the moduli
space of flat connections over $\Sigma$. We establish the basic regularity and
compactness properties (assuming $L^p$-bounds on the curvature for $p>2$) as well as the
Fredholm theory in a compact model case. The motivation for studying this
boundary value problem lies in the construction of an instanton Floer homology for
$3$-manifolds with boundary. The present paper is part of a program proposed by
Salamon for the proof of the Atiyah-Floer conjecture for homology-$3$-spheres.

\end{abstract}

\section{Introduction}

Let $X$ be a manifold with boundary, let $\rG$ be a compact Lie group, and consider
a principal $\rG$-bundle $P\to M$.
The natural boundary condition for the Yang-Mills equation $\rd_A^*F_A=0$ on $P$
is $*F_A|_{\pd X}=0$. 
For this boundary value problem there are regularity and compactness results,
see for example \cite{U1,U2,W}.
Every solution is gauge equivalent to a smooth solution, and Uhlenbeck 
compactness holds: Every sequence of solutions with $L^p$-bounded curvature 
(where $2p>{\rm dim}\,X$) is gauge equivalent to a sequence that contains a
$\cC^\infty$-convergent subsequence.
On an oriented $4$-manifold, the anti-self-dual instantons, i.e.\ connections satisfying
$F_A+*F_A=0$, are special first order solutions of the Yang-Mills equation.
An important application of Uhlenbeck's theorem is the compactification
of the moduli space of anti-self-dual instantons over a manifold without boundary
leading to the Donaldson invariants of smooth $4$-manifolds~\cite{D2} and to 
the instanton Floer homology groups of $3$-manifolds~\cite{F1}. 

On a $4$-manifold with boundary the boundary condition $*F_A|_{\pd X}=0$ for anti-self-dual 
instantons implies that the curvature vanishes altogether at the boundary. 
This is an overdetermined boundary value problem comparable to Dirichlet boundary 
conditions for holomorphic maps.
As in the latter case it is natural to consider weaker Lagrangian boundary conditions.
The Cauchy-Riemann equation becomes elliptic when augmented with Lagrangian or more 
generally totally real boundary conditions.
We consider a version of such Lagrangian boundary conditions for anti-self-dual instantons 
on a $4$-manifold with a space-time splitting of the boundary,
and prove that they suffice to obtain the analogue of the above mentioned 
regularity and compactness results for Yang-Mills connections.

More precisely, we consider oriented $4$-manifolds $X$ such that each connected component of 
the boundary $\pd X$ is diffeomorphic to $\cS\times\S$, where $\cS$ is a $1$-manifold and 
$\S$ is a closed Riemann surface.
We shall study a boundary value problem associated to a gauge invariant
Lagrangian submanifold $\cL$ of the space of flat connections on $\S$:
The restriction of the anti-self-dual instanton to each time-slice of the boundary
is required to belong to $\cL$.
This boundary condition arises naturally from examining the Chern-Simons 
functional on a $3$-manifold $Y$ with boundary $\Si$. 
Namely, the Langrangian boundary condition renders the Chern-Simons $1$-form on the space of 
connections closed, see \cite{Sa1}.
The resulting gradient flow equation leads to the boundary value problem studied
in this paper (for the case $X=\R\times Y$).
Our main results establish the basic regularity and compactness properties as well
as the Fredholm theory, the latter for the compact model case $X=S^1\times Y$.

Boundary value problems for Yang-Mills connections were already considered by
Donaldson in \cite{D3}. He studies the Hermitian Yang-Mills equation 
for connections induced by Hermitian metrics on holomorphic bundles over a compact 
K\"ahler manifold $Z$ with boundary. Here the unique solubility of the Dirichlet 
problem (prescribing the metric over the boundary) leads to an identification
between framed holomorphic bundles over $Z$ (meaning a holomorphic bundle with a fixed 
trivialization over $\pd Z$) and Hermitian Yang-Mills connections over $Z$.
In particular, when $Z$ has complex dimension $1$ and boundary $\pd Z=S^1$, this links loop 
groups to moduli spaces of flat connections over $Z$. 
This observation suggests an alternative approach to Atiyah's \cite{A} 
correspondence between holomorphic curves in the loop group of a compact Lie group $\rG$
and anti-self-dual instantons on $\rG$-bundles over the $4$-sphere. 
The correspondence might be established via an adiabatic limit relating holomorphic spheres in  
the moduli space of flat connections over the disc to anti-self-dual instantons over the 
product (of sphere and disc).
Our motivation for studying the present boundary value problem lies more in the
direction of another such correspondence between holomorphic curves in moduli spaces 
of flat connections and anti-self-dual instantons -- the Atiyah-Floer conjecture for
Heegard splittings of a homology-$3$-sphere.

A Heegard splitting $Y=Y_0\cup_\S Y_1$ of a homology $3$-sphere $Y$ into two handlebodies 
$Y_0$ and $Y_1$ with common boundary $\S$ gives rise to two Floer homologies (i.e.\ 
generalized Morse homologies) as follows:
Firstly, the moduli space $M_\S$ of gauge equivalence classes of flat connections on the 
trivial $\SU(2)$-bundle over $\S$ is a symplectic manifold (with singularities) and the 
moduli spaces $L_{Y_i}$ of flat connections over $\S$ that extend to $Y_i$ are Lagrangian 
submanifolds of $M_\S$ as explained in \cite{W Cauchy}.
The symplectic Floer homology $\HF^{\rm symp}_*(M_\S,L_{Y_0},L_{Y_1})$ is now generated by
the intersection points of the Lagrangian submanifolds, and the generalized connecting orbits 
(that define the boundary operator) are pseudoholomorphic strips with boundary values in the 
two Lagrangian submanifolds.
It was conjectured by Atiyah \cite{A1} and Floer that this should be isomorphic to the instanton
Floer homology $\HF^{\rm inst}_*(Y)$. For the latter, the critical points are 
the flat $\SU(2)$-connections over $Y$. These are the actual critical points of the 
Chern-Simons functional, and the connecting orbits are given by its generalized flow lines, 
i.e.\ anti-self-dual instantons on $\R\times Y$.

The program by Salamon \cite{Sa1} for the proof of this conjecture is to define 
the instanton Floer homology $\HF^{\rm inst}_*(Y,L)$ for $3$-manifolds with boundary $\pd Y=\Si$ 
using boundary conditions associated with a Lagrangian submanifold $L\subset M_\Si$.
Then the conjectured isomorphism might be established in two steps via the 
intermediate $\HF^{\rm inst}_*([0,1]\times\Si,L_{Y_0}\times L_{Y_1})$.

Fukaya was the first to suggest the use of Lagrangian boundary conditions in order 
to define a Floer homology for $3$-manifolds $Y$ with boundary, \cite{Fu}.
He studies a slightly different equation, involving a degeneration of the metric
in the anti-self-duality equation, and uses ${\rm SO}(3)$-bundles that are
nontrivial over the boundary $\pd Y$.
Now there are interesting examples, where one has to work with the trivial bundle.
For example, on a handlebody $Y$ there exists no nontrivial $\rG$-bundle for connected 
$\rG$. So if one considers any $3$-manifold $Y$ with the Lagrangian submanifold 
$\cL_{Y'}$, the space of flat connections on $\pd Y=\pd Y'$ that extend over a handlebody $Y'$,
then one also deals with the trivial bundle.
Consequently, if one wants to use Floer homology on $3$-manifolds with boundary to prove the 
Atiyah-Floer conjecture, then it is crucial to extend this construction to the
case of trivial ${\rm SU}(2)$-bundles.
There are two approaches that suggest themselves for such a generalization.
One would be the attempt to extend Fukaya's construction to the case of trivial
bundles, and another would be to follow the alternative construction outlined
in \cite{Sa1}.
The present paper follows the second route and sets up the basic analysis for 
this theory. We will only consider trivial $\rG$-bundles. However, our main 
theorems~A, B, and C below generalize directly to nontrivial bundles 
-- just the notation becomes more cumbersome.

The main theorems are described below; they are proven in sections \ref{reg,comp} 
and~\ref{Fredholm}. The appendix reviews the regularity theory for the Neumann and
Dirichlet problem in the weak formulation that will be needed throughout this paper.
Here we moreover introduce a technical tool for extracting regularity results for
single components of a $1$-form from weak equations that are related to a combination 
of Neumann and Dirichlet problems.

\section*{The main results}
Throughout this paper, $P$ is the trivial $\rG$-bundle over a $4$-manifold $X$.
So a connection on $P$ is a $1$-form $A\in\Om^1(X,\cg)$ with values in the Lie algebra $\cg$. 
Its curvature is given by $F_A=\rd A + \half [A\wedge A]$.
For more details on gauge theory and the notation used in this paper see 
\cite{W Cauchy} or \cite{W}.

We consider the following class of Riemannian $4$-manifolds.
Here and throughout all Riemann surfaces are closed $2$-dimensional manifolds.
Moreover, unless otherwise mentioned, all manifolds are allowed to have a smooth 
boundary.

\begin{dfn} \label{def mfd}
A {\bf $\mathbf{4}$-manifold with a boundary space-time splitting} is a pair $(X,\t)$
with the following properties:
\begin{enumerate}
\item
$X$ is an oriented $4$-manifold which can be exhausted by a nested sequence of 
compact deformation retracts.
\item
$\t=(\t_1,\ldots,\t_n)$ is an $n$-tuple of embeddings $\t_i:\cS_i\times\Si_i\to X$
with disjoint images, where $\Si_i$ is a Riemann surface and $\cS_i$ is either an 
open interval in $\R$ or is equal to $S^1=\R/\Z$.
\item
The boundary $\pd X$ is the union 
$$
\pd X = \bigcup_{i=1}^n \t_i(\cS_i\times\Si_i).
$$
\end{enumerate}
\end{dfn}

\begin{dfn} \label{mfd}
Let $(X,\t)$ be a $4$-manifold with a boundary space-time splitting. 
A Riemannian metric $g$ on $X$ is called {\bf compatible} with $\t$ if for each 
$i=1,\ldots n$ there exists a neighbourhood $\cU_i\subset\cS_i\times[0,\infty)$ of
$\cS_i\times\{0\}$ and an extension of $\t_i$ to an embedding  
$\bar\t_i:\cU_i\times\Si_i\to X$ such that
$$
\bar\t_i^*g = \ds^2 + \dt^2 + g_{s,t} .
$$
Here $g_{s,t}$ is a smooth family of metrics on $\Si_i$ and we denote by $s$ the 
coordinate on $\cS_i$ and by $t$ the coordinate on $[0,\infty)$.

We call a triple $(X,\t,g)$ with these properties a
{\bf Riemannian $\mathbf{4}$-manifold with a boundary space-time splitting}.
\end{dfn}

\begin{rmk}\rm
In definition~\ref{mfd} the extended embeddings $\bar\t_i$ are 
uniquely determined by the metric as follows.
The restriction $\bar\t_i|_{t=0}=\t_i$ to the boundary is prescribed, and the paths
$t\mapsto\bar\t_i(s,t,z)$ are normal geodesics.
\end{rmk}

\begin{ex}\rm   \label{ex:interpolate metric}
Let $X:=\R\times Y$, where $Y$ is a compact oriented $3$-manifold with boundary 
$\pd Y=\Si$, and let $\t:\R\times\Si\to X$ be the obvious inclusion.
Given any two metrics $g_-$ and $g_+$ on $Y$ there exists a metric $g$ on $X$ such
that $g=\ds^2+g_-$ for $s\leq -1$,\, $g=\ds^2+g_+$ for $s\geq 1$, and
$(X,\t,g)$ satisfies the conditions of definition~\ref{mfd}.
The metric $g$ cannot necessarily be chosen in the form $\ds^2+g_s$
(one has to homotop the embeddings and the metrics).
\end{ex}

Now let $(X,\t,g)$ be a Riemannian $4$-manifold with a boundary space-time 
splitting and consider a trivial $\rG$-bundle over $X$ for a compact Lie group $\rG$.
The Sobolev spaces of connections and gauge transformations are denoted by
\begin{align*}
\cA^{k,p}(X) &= W^{k,p}(X,\rT^*X\otimes\cg), \\
\cG^{k,p}(X) &= W^{k,p}(X,\rG).
\end{align*}
Let $p>2$, then for each $i=1,\ldots,n$ the Banach space of connections $\cA^{0,p}(\S_i)$ 
carries the symplectic form $\o(\a,\b)=\int_{\S_i} \la \a\wedge\b\ra$.
Note that the Hodge $*$ operator for any metric on $\S_i$ is an $\o$-compatible complex 
structure on $\cA^{0,p}(\S_i)$ since $**=-{\rm id}$ and $\o(\cdot,*\cdot)$ defines a positive 
definite inner product -- the $L^2$-metric.
We call a submanifold $\cL\subset\cA^{0,p}(\S_i)$ Lagrangian if it is isotropic, i.e.\ 
$\o|_{\cL}\equiv 0$, and if $\rT_A\cL$ is maximal for all $A\in\cL$ in the following sense: 
If $\a\in\cA^{0,p}(\S_i)$ satisfies $\o(\a,\rT_A\cL)=\{0\}$, then $\a\in\rT_A\cL$.

We fix an $n$-tuple $\cL=(\cL_1,\ldots,\cL_n)$ of Lagrangian submanifolds
$\cL_i\subset\cA^{0,p}(\Si_i)$ that are contained in the space of flat connections and that
are gauge invariant,
$$
\cL_i\subset\cA_{\rm flat}^{0,p}(\Si_i) \qquad\text{and}\qquad
u^* \cL_i = \cL_i \quad \forall u\in\cG^{1,p}(\Si_i) .
$$
Here $\cA_{\rm flat}^{0,p}(\Si_i)$ is the space of weakly flat $L^p$-connections
on $\Si_i$ as introduced in \cite{W Cauchy}. It is shown in \cite[Lemma 4.2]{W Cauchy} 
that the assumptions on the $\cL_i$ imply that they are totally real with respect to the Hodge 
$*$ operator for any metric on $\S_i$, i.e.\ for all $A\in\cL_i$ one has the topological sum
$$
\cA^{0,p}(\S_i) = \rT_A \cL_i \oplus * \rT_A \cL_i  .
$$
Now we consider the following boundary value problem for connections
$A\in\cA^{1,p}_{\rm loc}(X)$
\footnote{
The subscript $\scriptstyle{\rm loc}$ indicates that the regularity
only holds on all compact subsets of the noncompact domain of definition.
}
\begin{equation}\label{ASD bvp intro}
\left\{\begin{array}{l}
*F_A + F_A = 0,\\
\t_i^*A|_{\{s\}\times\Si_i} \in\cL_i \quad\forall s\in\cS_i ,\, i=1,\ldots,n .
\end{array}\right.
\end{equation}
Observe that the boundary condition is meaningful since for every neighbourhood $\cU\times\S$
of a boundary component one has the continuous embedding
$W^{1,p}(\cU\times\S)\subset W^{1,p}(\cU,L^p(\S))\hookrightarrow \cC^0(\cU,L^p(\S))$.
The first nontrivial observation is that every connection in $\cL_i$ is gauge
equivalent to a smooth connection on $\Si_i$ and hence $\cL_i\cap\cA(\Si)$ is 
dense in $\cL_i$, as shown in \cite[Theorem~3.1]{W Cauchy}.
Moreover, the $\cL_i$ are modelled on $L^p$-spaces, and every $W^{1,p}_{\rm loc}$-connection on $X$ 
satisfying the boundary condition in (\ref{ASD bvp intro}) can be locally approximated by smooth 
connections satisfying the same boundary condition, see \cite[Corollary~4.4,~4.5]{W Cauchy}.

Note that the present boundary value problem is a first order equation with first order 
boundary conditions (flatness in each time-slice). Moreover, the boundary conditions contain 
some crucial nonlocal (i.e.\ Lagrangian) information.
We moreover emphasize that while $\cL_i$ is a smooth Banach submanifold of $\cA^{0,p}(\Si_i)$, 
the quotient $\cL_i/\cG^{1,p}(\Si_i)$ is not required to be a smooth submanifold of the 
moduli space $M_{\S_i}:=\cA^{0,p}_{\rm flat}(\Si_i)/\cG^{1,p}(\Si_i)$, which itself might be
singular.

For example, $\cL_i$ could be the set of flat connections on $\S_i$ that extend to
flat connections over a handlebody with boundary $\S_i$, as introduced in \cite[Lemma~4.6]{W Cauchy}.
To overcome the difficulties arising from the singularities in the quotient, we work with the 
(smooth) quotient by the based gauge group.

The following two theorems are the main regularity and compactness results for the
solutions of (\ref{ASD bvp intro}) generalizing the regularity theorem and the Uhlenbeck
compactness for Yang-Mills connections on $4$-manifolds without boundary.
They will be proven in section~\ref{reg,comp}.\\

\noindent
{\bf Theorem A \; (Regularity)}\\  {\it
Let $p>2$. Then every solution $A\in\cA^{1,p}_{\rm loc}(X)$ of 
(\ref{ASD bvp intro}) is gauge equivalent to a smooth solution, i.e.\
there exists a gauge transformation $u\in\cG^{2,p}_{\rm loc}(X)$ such that $u^*A\in\cA(X)$
is smooth.
}\\

\noindent
{\bf Theorem B \; (Compactness)}\\  {\it
Let $p>2$ and let $g^\n$ be a sequence of metrics compatible with $\t$ that 
uniformly converges with all derivatives on every compact set to a smooth metric.
Suppose that $A^\n\in\cA^{1,p}_{\rm loc}(X)$ is a sequence of solutions 
of (\ref{ASD bvp intro}) with respect to the metrics $g^\n$ such that for 
every compact subset $K\subset X$ there is a uniform bound on the curvature
$\|F_{A^\n}\|_{L^p(K)}$.
Then there exists a subsequence (again denoted $A^\n$) and a sequence of gauge
transformations $u^\n\in\cG^{2,p}_{\rm loc}(X)$ such that $u^{\n\;*}A^\n$ 
converges uniformly with all derivatives on every compact set to a smooth 
connection $A\in\cA(X)$.
}\\

The difficulty of these results lies in the global nature of the boundary condition.
This makes it impossible to directly generalize the proof of the regularity and compactness
theorems for Yang-Mills connections, where one chooses suitable local gauges, obtains the 
higher regularity and estimates from an elliptic boundary value problem, and then patches the
gauges together. With our global Lagrangian boundary condition one cannot obtain local
regularity results.

However, an approach by Salamon can be generalized to manifolds with boundary. Firstly,
Uhlenbeck's weak compactness theorem yields a weakly $W^{1,p}_{\rm loc}$-convergent 
subsequence. Its limit serves as reference connection with respect to which a further
subsequence can be put into relative Coulomb gauge globally (on large compact sets).
Then one has to establish elliptic estimates and regularity results for the given boundary
value problem together with the relative Coulomb gauge equations.
The crucial point in this last step is to establish the higher regularity for the 
$\Si$-component of the connections in a neighbourhood $\cU\times\Si$ of a boundary component.
The global nature of the boundary condition forces us to deal with a 
Cauchy-Riemann equation on $\cU$ with values in the Banach space $\cA^{0,p}(\Si)$ 
and with Lagrangian boundary conditions.
At this point we will make use of the regularity results in \cite{W Cauchy} that are
established in the general framework of a Cauchy-Riemann equation for functions with
values in a complex Banach space and with totally real boundary conditions.
Some more general analytic tools for this approach will be taken from \cite{W}.

The case $2<p\leq 4$, when $W^{1,p}$-functions are not automatically continuous, poses some 
special difficulties in this last step. Firstly, in order to obtain regularity results from
the Cauchy-Riemann equation, one has to straighten out the Lagrangian submanifold by going to
suitable coordinates. This requires a $\cC^0$-convergence of the connections, which in case
$p>4$ is given by a standard Sobolev embedding. In case $p>2$ one still obtains a special
compact embedding $W^{1,p}(\cU\times\S)\hookrightarrow\cC^0(\cU,L^p(\S))$ that suits our
purposes. Secondly, the straightening of the Lagrangian introduces a nonlinearity in the 
Cauchy-Riemann equation that already poses some problems in case $p>4$. In case $p\leq 4$ this
forces us to deal with the Cauchy-Riemann equation with values in an $L^2$-Hilbert space and
then use some interpolation inequalities for Sobolev norms.

For the definition of the standard instanton Floer homology it suffices to prove a compactness result 
like theorem~B for $p=\infty$.
In our case however the bubbling analysis \cite{W bubbling} requires the compactness result for some
$p<3$. This is why we have taken some care to deal with this case.\\

In order to define a Floer homology for 3-manifolds with boundary as outlined in \cite{Sa1}
one has to consider the moduli space of finite energy solutions of (\ref{ASD bvp intro}),
writing $\cL$ for the $n$-tuple of Lagrangian submanifolds $\cL_i$, 
$$
\cM(\cL) := 
\bigl\{ A\in\cA^{1,p}_{\rm loc}(X) \st A \;\text{satisfies}\; (\ref{ASD bvp intro}),
                                       \|F_A\|_{L^2} <\infty \bigr\}
/ \cG^{2,p}_{\rm loc}(X).
$$
Theorem~A implies that for every equivalence class $[A]\in\cM(\cL)$ one can 
find a smooth representative $A\in\cA(X)$.
Theorem~B is one step towards a compactness result for $\cM(\cL)$: 
Every closed subset of $\cM(\cL)$ with a uniform $L^p$-bound for the curvature 
is compact. 
In addition, theorem~B allows the metric to vary, which is relevant for the
metric-independence of the Floer homology.

Our third main result is the Fredholm theory in section \ref{Fredholm}. It is a step towards 
proving that the moduli space $\cM(\cL)$ of solutions of (\ref{ASD bvp intro}) is a manifold 
whose components have finite (but possibly different) dimensions. 
This also exemplifies our hope that the further analytical details of Floer theory will
work out along the usual lines once the right analytic setup has been found in the proof 
of theorems~A and B.

In the context of Floer homology and in Floer-Donaldson theory it is important
to consider $4$-manifolds with cylindrical ends. 
This requires an analysis of the asymptotic behaviour which will be carried out
elsewhere.
Here we shall restrict the discussion of the Fredholm theory to the compact case. 
The crucial point is the behaviour of the linearized operator near the
boundary; in the interior we are dealing with the usual anti-self-duality equation.
Hence it suffices to consider the following model case.
Let $Y$ be a compact oriented $3$-manifold with boundary $\pd Y=\Si$
and suppose that $(g_s)_{s\in S^1}$ is a smooth family of metrics on $Y$
such that
$$
X=S^1\times Y, \qquad  \t:S^1\times\Si\to X, \qquad  g=\ds^2 + g_s
$$ 
satisfy the assumptions of definition \ref{mfd}.
Here the space-time splitting $\t$ of the boundary is the obvious inclusion
$\t:S^1\times\S \hookrightarrow \pd X = S^1\times\S$, where $\S=\bigcup_{i=1}^n \S_i$ might 
be a disjoint union of Riemann surfaces $\S_i$. An $n$-tuple of Lagrangian submanifolds
$\cL_i\subset\cA^{0,p}(\S_i)$ as above then defines a gauge invariant Lagrangian submanifold 
$\cL:=\cL_1\times\ldots\times\cL_n$ of the symplectic Banach space
$\cA^{0,p}(\S)=\cA^{0,p}(\S_1)\times\ldots\times\cA^{0,p}(\S_n)$ 
such that $\cL\subset\cA^{0,p}_{\rm flat}(\S)$.

In order to linearize the boundary value problem (\ref{ASD bvp intro}) together with the
local slice condition, fix a smooth connection $A+\P\ds\in\cA(S^1\times Y)$ such that
$A_s:=A(s)|_{\pd Y}\in\cL$ for all $s\in S^1$.
Here $\P\in\cC^\infty(S^1\times Y,\cg)$, and 
$A\in\cC^\infty(S^1\times Y,\rT^*Y\otimes\cg)$ is an $S^1$-family
of $1$-forms on $Y$ (not a $1$-form on $X$ as previously).
Now let $E_A^{1,p}$ be the space of $S^1$-families of $1$-forms
$\a\in W^{1,p}(S^1\times Y,\rT^*Y\otimes\cg)$ that satisfy the boundary 
conditions
$$
*\a(s)|_{\pd Y} = 0 \qquad\text{and}\qquad
 \a(s)|_{\pd Y}\in\rT_{A_s}\cL \qquad\text{for all}\; s\in S^1 .
$$
Then the linearized operator
$$
D_{(A,\P)} :
E_A^{1,p} \times W^{1,p}(S^1\times Y,\cg)  \longrightarrow 
L^p(S^1\times Y,\rT^*Y\otimes\cg) \times L^p(S^1\times Y,\cg) 
$$
is given with $\nabla_s = \pd_s + [\P,\cdot]$ by
$$
D_{(A,\P)}(\a,\ph) \;=\; \bigl( \nabla_s\a  - \rd_A\ph + *\rd_A\a \,,\, 
                                \nabla_s\ph - \rd_A^*\a \bigr) .
$$
The second component of this operator is $-\rd_{A+\P\ds}^*(\a+\ph\ds)$, and the
first boundary condition is $*(\a+\ph\ds)|_{\pd X}=0$, corresponding to the 
choice of a local slice at $A+\P\ds$.
In the first component of $D_{(A,\P)}$ we have used the global space-time splitting 
of the metric on $S^1\times Y$ to identify the self-dual $2$-forms
$*\g_s - \g_s\wedge\ds$ with families $\g_s$ of $1$-forms on $Y$. The vanishing
of this component is equivalent to the linearization $\rd_{A+\P\ds}^+(\a+\ph\ds)=0$
of the anti-self-duality equation.
Furthermore, the boundary condition $\a(s)|_{\pd Y}\in\rT_{A_s}\cL$ is the linearization
of the Lagrangian boundary condition in the boundary value problem (\ref{ASD bvp intro}).\\

\noindent
{\bf Theorem C \; (Fredholm properties)}\\  {\it
Let $Y$ be a compact oriented $3$-manifold with boundary $\pd Y=\Si$ 
and let $S^1\times Y$ be equipped with a 
product metric $\ds^2+g_s$ that is compatible with $\t:S^1\times\Si \to S^1\times Y$.
Let $A+\P\ds\in\cA(S^1\times Y)$ such that $A(s)|_{\pd Y}\in\cL$ for all $s\in S^1$. 
Then the following holds for all $p>2$.
\begin{enumerate}
\item $D_{(A,\P)}$ is Fredholm.
\item
There is a constant $C$ such that for all
$\a\in E^{1,p}_A $ and $\ph\in W^{1,p}(S^1\times Y,\cg)$
$$
\| (\a,\ph) \|_{W^{1,p}}
\leq C \bigl( \| D_{(A,\P)}(\a,\ph) \|_{L^p}  +\| (\a,\ph) \|_{L^p}  \bigr).
$$
\item
Let $q\geq p^*$ such that $q\neq 2$.
Suppose that $\b\in L^q(S^1\times Y,\rT^*Y\otimes\cg)$, 
$\z\in L^q(S^1\times Y,\cg)$, and assume that there exists a constant $C$ such that 
for all $\a\in E_A^{1,p}$ and $\ph\in W^{1,p}(S^1\times Y,\cg)$
$$
\left| \int_{S^1\times Y} \la D_{(A,\P)}(\a,\ph) \,,\, (\b,\z) \ra  \right|
\;\leq\; C \, \|(\a,\ph)\|_{L^{q^*}} .
$$
Then in fact $\b\in W^{1,q}(S^1\times Y,\rT^*Y\otimes\cg)$
and $\z\in W^{1,q}(S^1\times Y,\cg)$ .
\end{enumerate}
}

Here and throughout we use the notation $\frac 1p + \frac 1{p^*} = 1$
for the conjugate exponent $p^*$ of $p$.
The above inner product $\la\cdot,\cdot\ra$ is the pointwise inner product in 
$\rT^*Y\otimes\cg \times \cg$.
The reason for our assumption $q\neq 2$ in theorem~C~(iii) is a technical
problem in dealing with the singularities of $\cL/\cG^{1,p}(\Si)$. 
We resolve these singularities by dividing only by the based gauge group.
This leads to coordinates of $L^p(\Si,\rT^*\Si\otimes\cg)$ in a Banach space that
comprises based Sobolev spaces $W^{1,p}_z(\Si,\cg)$ of functions vanishing at a
fixed-point $z\in\Si$. So these coordinates that straighten out $\rT\cL$ along 
$A|_{S^1\times\pd Y}$ are welldefined only for $p>2$.
Now in order to prove the regularity claimed in theorem~C~(iii) we have to use such 
coordinates either for $\b$ or for the test $1$-forms $\a$, i.e.\ we have to
assume that either $q>2$ or $q^*>2$.
This is completely sufficient for our purposes -- concluding a higher regularity of elements
of the cokernel. This will be done via an iteration of theorem~C~(iii) that can always be
chosen such as to jump across $q=2$.
However, we believe that the use of different coordinates should permit to extend this result.\\

\noindent
{\bf Conjecture}  {\it
Theorem C (iii) continues to hold for $q=2$.
}\\

One indication for this conjecture is that the $L^2$-estimate in theorem~C~(ii) is true
(for $W^{1,p}$-regular $\a$ and $\p$ with $p>2$), as will be shown in section~\ref{Fredholm}.
This $L^2$-estimate can be proven by a much more elementary method than the general
$L^p$-regularity and -estimates. In fact, it was already stated in \cite{Sa1} as an indication
for the wellposedness of the boundary value problem (\ref{ASD bvp intro}).

\section*{Outlook}

We give a brief sketch of Salamon's program for the proof of the Atiyah-Floer conjecture
(for more details see \cite{Sa1}) in order to point out the significance of this paper
for the whole program.

The first step of the program is to define the instanton Floer homology 
$\HF^{\rm inst}_*(Y,L)$ for a 3-manifold with boundary $\pd Y=\S$ and a 
Lagrangian submanifold $L=\cL/\cG^{1,p}(\S)\subset M_\Si$ in the moduli space of flat connections.
The Floer complex will be generated by the gauge equivalence classes of irreducible flat 
connections $A\in\cA(Y)$ with Lagrangian boundary conditions $A|_\S\in\cL$. 
\footnote{
A connection $A\in\cA_{\rm flat}(Y)$ is called irreducible if its isotropy subgroup
of $\cG(Y)$ (the group of gauge transformations that leave $A$ fixed) is discrete,
i.e.\ $\rd_A|_{\Om^0}$ is injective. 
There should be no reducible flat connections with Lagrangian boundary conditions other
than the gauge orbit of the trivial connection. 
This will be guaranteed by certain conditions on $Y$ and $L$, 
for example this is the case when $L=L_{Y'}$ for a handlebody $Y'$ with 
$\pd Y'=\bar\Si$ such that $Y\cap_\Si Y'$ is a homology-3-sphere.
}
For any two such connections $A^+,A^-$ one then has to study the moduli
space of Floer connecting orbits,
$$
\cM(A^-,A^+) = \bigl\{ \tA\in\cA(\R\times Y) \st \tA \;\text{satisfies}\; 
                          (\ref{ASD bvp intro}),
                          \lim_{s\to\pm\infty}\tA = A^\pm \bigr\} / \cG(\R\times Y) .
$$
Theorem A shows that the boundary value problem (\ref{ASD bvp intro}) is wellposed.
In particular, the spaces of smooth connections and gauge transformations in the definition
of the above moduli space can be replaced by suitable Sobolev completions.
The next step in the construction of the Floer homology groups is the analysis of the asymptotic 
behaviour of the finite energy solutions of (\ref{ASD bvp intro}) on $\R\times Y$, which will be
carried out elsewhere.
Combining this with theorem~C one obtains an appropriate Fredholm theory and proves that 
for a suitably generic perturbation the spaces $\cM(A^-,A^+)$ are smooth manifolds.
In the monotone case the connections in the $k$-dimensional part $\cM^k(A^-,A^+)$ 
have a fixed energy.

Theorem B is a major step towards a compactification of these moduli spaces $\cM^k(A^-,A^+)$. 
It proves their compactness under the assumption of an $L^p$-bound on the curvature for 
$p>2$, whereas the energy is only the $L^2$-norm.
So the key remaining analytic task is an analysis of the possible bubbling phenomena.
This will be carried out in \cite{W bubbling} and draws upon the techniques developed in this paper.
When this is understood, the construction of the Floer homology groups should be routine.
In particular, for the metric independence note that one can interpolate between different metrics 
on $Y$ as in example \ref{ex:interpolate metric}, and theorem~B allows for the variation of metrics
on $X$.
So this paper sets up the basic analytic framework for the Floer theory of $3$-manifolds with
boundary.

The further steps in the program for the proof of the Atiyah-Floer conjecture are to 
consider a Heegard splitting $Y=Y_0\cup_\Sigma Y_1$ of a homology 3-sphere, and 
identify $\HF^{\rm inst}_*([0,1]\times\Si,L_{Y_0}\times L_{Y_1})$ with
$\HF^{\rm inst}_*(Y)$ and $\HF^{\rm symp}_*(M_\Si,L_{Y_0},L_{Y_1})$ respectively.
In both cases, the Floer complexes can be identified by elementary arguments, so the main 
task is to identify the connecting orbits.

In the case of the two instanton Floer homologies, the idea is to choose an embedding 
$(0,1)\times\Si\hookrightarrow Y$ starting from a tubular neighbourhood of 
$\Si\subset Y$ at $t=\half$ and shrinking $\{t\}\times\Si$ to the $1$-skeleton of $Y_t$ 
for $t=0,1$. Then the anti-self-dual instantons on $\R\times Y$ pull back to anti-self-dual 
instantons on $\R\times[0,1]\times\Si$ with a degenerate metric for $t=0$ and $t=1$.
On the other hand, one can consider anti-self-dual instantons on $\R\times[\ep,1-\ep]\times\Si$ 
with boundary values in $\cL_{Y_0}$ and $\cL_{Y_1}$. As $\ep\to 0$, one should be able 
to pass from this genuine boundary value problem to solutions on the closed manifold $Y$. 
This is a limit process for the boundary value problem studied in this paper.

The identification of the instanton and symplectic Floer homologies requires an adaptation 
of the adiabatic limit argument in~\cite{DS} to boundary value problems for anti-self-dual 
instantons and pseudoholomorphic curves respectively 
Here one again deals with the boundary value problem (\ref{ASD bvp intro}) studied in this 
paper. As the metric on $\Si$ converges to zero, the solutions, i.e.\ anti-self-dual instantons
on $\R\times [0,1]\times\Si$ with Lagrangian boundary conditions in 
$\cL_{Y_0},\cL_{Y_1}$ should be in one-to-one correspondence with connections 
on $\R\times [0,1]\times\Si$ that descend to pseudoholomorphic strips in $M_\Si$ with 
boundary values in $L_{Y_0}$ and~$L_{Y_1}$.
The basic elliptic properties of the boundary value problem (\ref{ASD bvp intro}) that
are established in this paper will also play an important role in this adiabatic limit
analysis.\\

I would like to thank Dietmar Salamon for his constant help and encouragement in pursueing 
this project.

\section{Regularity and compactness}
\label{reg,comp}

Let $(X,\t)$ be a $4$-manifold with boundary space-time splitting.
This means that $X$ is oriented and 
$$
X=\bigcup_{k\in\N} X_k,
$$ 
where all $X_k$ are compact submanifolds and deformation retracts of $X$ such that
$X_k\subset{\rm int}\,X_{k+1}$ for all $k\in\N$.
Here the interior of a submanifold $X'\subset X$ is to be understood with respect to
the relative topology, i.e.\ we define ${\rm int}\,X'=X\setminus {\rm cl}(X\setminus X')$.
Moreover, 
$$
\pd X = \bigcup_{i=1}^n \t_i(\cS_i\times\Si_i),
$$ 
where each $\Si_i$ is a Riemann surface, each $\cS_i$ is either an open interval in 
$\R$ or is equal to $S^1=\R/\Z$, and the embeddings $\t_i:\cS_i\times\Si_i\to X$
have disjoint images.
We then consider the trivial $\rG$-bundle over $X$, where $\rG$ is a compact Lie 
group with Lie algebra $\cg$.
For $i=1,\ldots,n$ let $\cL_i\subset\cA^{0,p}(\Si_i)$ be a Lagrangian submanifold
and suppose that
$$
\cL_i\subset\cA^{0,p}_{\rm flat}(\Si_i) \qquad\text{and} \qquad
\cG^{1,p}(\S_i)^*\cL_i \subset\cL_i .
$$
Furthermore, let $X$ be equipped with a metric $g$ that is compatible with the
space-time splitting $\t$.
This means that for each $i=1,\ldots,n$ the map 
$\cS_i\times[0,\infty)\times\Si_i \to X$, $(s,t,z)\mapsto \g_{(s,z)}(t)$ given by 
the normal geodesics $\g_{(s,z)}$ starting at $\g_{(s,z)}(0)=\t_i(s,z)$ restricts 
to an embedding $\bar\t_i:\cU_i\times\Si_i \hookrightarrow X$ for some neighbourhood
$\cU_i\subset\cS_i\times[0,\infty)$ of $\cS_i\times\{0\}$.
Now consider the boundary value problem (\ref{ASD bvp intro}) for connections
$A\in\cA^{1,p}_{\rm loc}(X)$, restated below.
\begin{equation}\label{ASD bvp}
\left\{\begin{array}{l}
*F_A + F_A =0 ,\\
\t_i^*A|_{\{s\}\times\Si_i} \in\cL_i \quad\forall s\in\cS_i,\,i=1,\ldots,n .
\end{array}\right.
\end{equation}
The anti-self-duality equation is welldefined for 
$A\in\cA^{1,p}_{\rm loc}(X)$ with any $p\geq 1$, but in order to be able to state the 
boundary condition correctly we have to assume $p>2$. 
Then the trace theorem for Sobolev spaces (e.g. \cite[Theorem~6.2]{Adams}) 
ensures that $\t_i^*A |_{\{s\}\times\Si_i}\in \cA^{0,p}(\Si_i)$ for all $s\in\cS_i$.

%
%
%
%
%

The aim of this section is to prove the regularity theorem~A and the compactness 
theorem~B for this boundary value problem.
Both theorems are dealing with the noncompact base manifold $X$.
However, we shall use an extension argument by Donaldson and Kronheimer
\cite[Lemma~4.4.5]{DK} to reduce the problem to compact base manifolds. 
For the following special version of this argument a detailed proof can be found in
\cite[Propositions 8.6,10.8]{W}.
At this point, the assumption that the exhausting compact submanifolds $X_k$ are 
deformation retracts of $X$ comes in crucially. It ensures that every gauge transformation
on $X_k$ can be extended to $X$, which is a central point in the argument of Donaldson and
Kronheimer that proves the following proposition.

\begin{prp} \label{DKargument}
Let the $4$-manifold $\tM=\bigcup_{k\in\N}M_k$ be exhausted by compact 
submanifolds $M_k\subset{\rm int}\,M_{k+1}$ that are deformation retracts of 
$\tM$, and let $p>2$.
\begin{enumerate}
\item
Let $A\in\cA^{1,p}_{\rm loc}(\tM)$ and suppose that for each $k\in\N$ there exists
a gauge transformation $u_k\in\cG^{2,p}(M_k)$ such that $u_k^*A|_{M_k}$ is 
smooth.
Then there exists a gauge transformation $u\in\cG^{2,p}_{\rm loc}(\tM)$ such that 
$u^*A$ is smooth. 
\item
Let a sequence of connections $(A^\n)_{\n\in\N}\subset\cA^{1,p}_{\rm loc}(\tM)$ be 
given and suppose that the following holds:

For every $k\in\N$ and every subsequence of $(A^\n)_{\n\in\N}$ there exist a 
further subsequence $(\n_{k,i})_{i\in\N}$ and gauge transformations 
$u^{k,i}\in\cG^{2,p}(M_k)$ such that
$$
\sup_{i\in\N}\; \bigl\| u^{k,i\;*}A^{\n_{k,i}} \bigr\|_{W^{\ell,p}(M_k)}
< \infty \qquad\forall \ell\in\N .
$$
Then there exists a subsequence $(\n_i)_{i\in\N}$ and a sequence of gauge 
transformations $u^i\in\cG^{2,p}_{\rm loc}(\tM)$ such that
$$
\sup_{i\in\N}\; \bigl\| u^{i\;*}A^{\n_i} \bigr\|_{W^{\ell,p}(M_k)} < \infty
\qquad\forall k\in\N, \ell\in\N .
$$
\end{enumerate}
\end{prp}

So in order to prove theorem~A it suffices to find smoothing gauge 
transformations on the compact submanifolds $X_k$ in view of proposition~\ref{DKargument}~(i).
For that purpose we shall use the so-called local slice theorem.
The following version is proven e.g.\ in \cite[Theorem 9.1]{W}.
Note that we are dealing with trivial bundles, so we will be using the 
product connection as reference connection in the definition of the Sobolev 
norms of connections.

\begin{prp} \label{local slice thm} {\bf (Local Slice Theorem)} \\
Let $M$ be a compact $4$-manifold, let $p>2$, and let $q>4$ be such that
$\frac 1q > \frac 1p - \frac 14$ (or $q=\infty$ in case $p>4$).
Fix $\hat A\in\cA^{1,p}(M)$ and let a constant $c_0>0$ be given. 
Then there exist constants $\ep>0$ and $C_{\scriptscriptstyle CG}$ such that 
the following holds.
For every $A\in\cA^{1,p}(M)$ with
\begin{equation*}
\|A-\hat A\|_q \leq \ep   \qquad\text{and}\qquad
\|A-\hat A\|_{W^{1,p}} \leq c_0
\end{equation*}
there exists a gauge transformation $u\in\cG^{2,p}(M)$ such that 
\[
\left\{\begin{aligned} 
\rd_{\hat A}^*(u^*A-\hat A)&=0, \\
*(u^*A-\hat A)|_{\pd M}&=0, 
\end{aligned}\right.
\qquad\text{and}\qquad
\begin{aligned} 
\|u^*A-\hat A\|_q \quad\;
&\leq C_{\scriptscriptstyle CG} \|A-\hat A\|_q ,\\
\|u^*A-\hat A\|_{W^{1,p}} 
&\leq C_{\scriptscriptstyle CG} \|A-\hat A\|_{W^{1,p}}.
\end{aligned}
\]
\end{prp}

\begin{rmk} \hspace{1mm} \label{local slice rmk} \\ \rm
\vspace{-5mm}
\begin{enumerate}
\item 
If the boundary value problem in proposition~\ref{local slice thm} is satisfied
one says that $u^*A$ is in Coulomb gauge relative to $\hat A$.
This is equivalent to $v^*\hat A$ being in Coulomb gauge relative to $A$ for
$v=u^{-1}$, i.e.\ the boundary value problem can be replaced by
\[
\left\{\begin{aligned} 
\rd_A^*(v^*\hat A-A)&=0, \\
*(v^*\hat A- A)|_{\pd M}&=0 .
\end{aligned}\right.
\]
\item
The assumptions in proposition~\ref{local slice thm} on $p$ and $q$ guarantee that one
has a compact Sobolev embedding 
$$
W^{1,p}(M)\hookrightarrow L^q(M).
$$
\item
One can find uniform constants for varying metrics in the following sense.
Fix a metric $g$ on $M$. 
Then there exist constants $\ep,\d>0$, and $C_{CG}$ such that the assertion of 
proposition~\ref{local slice thm}
holds for all metrics $g'$ with $\|g-g'\|_{\cC^1}\leq\d$.
\end{enumerate}
\end{rmk}

In the following we shall briefly outline the proof of theorem~A.
Given a solution $A\in\cA^{1,p}_{\rm loc}(X)$ of (\ref{ASD bvp}) one fixes $k\in\N$ 
and proves the assumption of proposition~\ref{DKargument}~(i) as follows.
One finds some sufficiently large compact submanifold $M\subset X$ with $X_k \subset M$.
Then one chooses a smooth connection $A_0\in\cA(M)$ sufficiently 
$W^{1,p}$-close to $A$ and applies the local slice theorem with the reference connection 
$\hat A=A$ to find a gauge tranformation that puts $A_0$ into relative Coulomb gauge with 
respect to $A$. This is equivalent to finding a gauge transformation that puts $A$ into 
relative Coulomb gauge with respect to $A_0$.  
We denote this gauge transformed connection again by $A\in\cA^{1,p}(M)$. 
It satisfies the following boundary value problem:
\begin{equation}\label{bvp}
\left\{\begin{aligned} 
\rd_{A_0}^*(A-A_0)&=0, \\
*F_A + F_A &=0, \\
*(A-A_0)|_{\pd M}&=0 , \\
\t_i^*A|_{\{s\}\times\Si_i} &\in\cL_i \quad\forall s\in\cS_i,\, i=1,\ldots,n .
\end{aligned}\right.
\end{equation}
More precisely, the Lagrangian boundary condition only holds for those
\hbox{$s\in\cS_i$} and $i\in\{1,\ldots n\}$ for which $\t_i(\{s\}\times\Si_i)$ is entirely 
contained in $\pd M$.
If $M$ was chosen large enough, then the regularity theorem~\ref{reg} below will assert the 
smoothness of $\tA$ on $X_k$.

The proof outline of the proof of theorem~A goes along similar lines. 
We will use proposition~\ref{DKargument}~(ii) to reduce the problem to compact base 
manifolds. On these, we shall use the following weak Uhlenbeck compactness theorem 
(see \cite{U1}, \cite[Theorem 8.1]{W}) to find a subsequence of gauge 
equivalent connections that converges $W^{1,p}$-weakly.

\begin{prp} \label{weak comp} {\bf (Weak Uhlenbeck Compactness)} \\
Let $M$ be a compact $4$-manifold and let $p>2$.
Suppose that the sequence of connections $A^\n \in \cA^{1,p}(M)$ is such that $\|F_{A^\n}\|_p$ 
is uniformly bounded. 
Then there exists a subsequence (again denoted $(A^\n)_{\n\in\N}$) and a 
sequence $u^\n\in\cG^{2,p}(M)$ of gauge transformations such that 
$u^{\n\;*}A^\n$ weakly converges in $\cA^{1,p}(M)$.
\end{prp}

The limit $A_0$ of the convergent subsequence then serves as reference connection $\hat A$
in the local slice theorem, proposition~\ref{local slice thm}, and this way one obtains a 
$W^{1,p}$-bounded sequence of connections $\tA^\n$ that solve the boundary 
value problem (\ref{bvp}).
This makes crucial use of the compact Sobolev embedding $W^{1,p}\hookrightarrow L^q$
on compact $4$-manifolds (with $q$ from the local slice theorem).
The estimates in the subsequent theorem~\ref{reg} then provide the higher $W^{k,p}$-bounds
on the connections that will imply the compactness.
One difficulty in the proof of this regularity theorem is that due to the global nature of the
boundary conditions one has to consider the $\Si$-components of the connections near the 
boundary as maps into the Banach space $\cA^{0,p}(\Si)$ that solve a Cauchy-Riemann equation 
with Lagrangian boundary conditions.
In order to prove a regularity result for such maps one has to straighten out the 
Lagrangian submanifold by using coordinates in $\cA^{0,p}(\Si)$.
(This is done in \cite{W Cauchy}.)
Thus on domains $\cU\times\Si$ at the boundary a crucial assumption is that the 
$\Si$-components of the connections all lie in one such coordinate chart, that
is one needs the connections to converge strongly in the $L^\infty(\cU,L^p(\Si))$-norm.
In the case $p>4$ this is ensured by the compact embedding 
$W^{1,p}\hookrightarrow L^\infty$ on $\cU\times\Si$. 
To treat the case $2<p\leq 4$ we shall make use of the following special Sobolev embedding. 
The proof mainly uses techniques from \cite{Adams}.

\begin{lem} \label{Sob emb}
Let $M,N$ be compact manifolds and let $p>m=\dim M$ and $p>n=\dim N$. 
Then the following embedding is compact,
$$
W^{1,p}(M\times N) \hookrightarrow L^\infty(M,L^p(N)).
$$
\end{lem}

\noindent
{\bf Proof of lemma \ref{Sob emb}: } \\
Since $M$ is compact it suffices to prove the embedding in (finitely many) coordinate charts.
These can be chosen as either balls $B_2\subset\R^m$ in the interior or half balls 
$D_2=B_2\cap\H^m$ in the half space \hbox{$\H^m=\{x\in\R^m\st x_1\geq 0\}$} at the boundary
of $M$. We can choose both of radius $2$ but cover $M$ by balls and half balls of radius $1$.
So it suffices to consider a bounded set \hbox{$\cK\subset W^{1,p}(B_2\times N)$} and prove
that it restricts to a precompact set in $L^\infty(B_1,L^p(N))$, and similarly with the
half balls.
Here we use the Euclidean metric on $\R^m$, which is equivalent to the metric induced from $M$.

For a bounded subset $\cK\subset W^{1,p}(D_2\times N)$ over the half ball define the subset
$\cK'\subset W^{1,p}(B_2\times N)$ by extending all $u\in\cK$ to $B_2\setminus\H^m$ by 
$u(x_1,x_2,\ldots,x_m):=u(-x_1,x_2,\ldots,x_m)$ for $x_1\leq 0$.
The thus extended function is still $W^{1,p}$-regular with twice the norm of $u$.
So $\cK'$ also is a bounded subset, and if this restricts to a precompact set in 
$L^\infty(B_1,L^p(N))$, then also $\cK\subset L^\infty(D_1,L^p(N))$ is compact.
Hence it suffices to consider the interior case of the full ball.

The claimed embedding is continuos by the standard Sobolev estimates -- check for example
in \cite{Adams} that the estimates generalize directly to functions with values in 
a Banach space. In fact, one obtains an embedding
$$
W^{1,p}(B_2\times N) \subset W^{1,p}(B_2,L^p(N)) \hookrightarrow \cC^{0,\l}(B_2,L^p(N))
$$
into some H\"older space with $\l=1-\frac mp>0$. 
One can also use this Sobolev estimate for $W^{1,p}(N)$ with $\l'=1-\frac np>0$ combined 
with the inclusion $L^p\hookrightarrow L^1$ on $B_2$ to obtain a continuous embedding
$$
W^{1,p}(B_2\times N) \subset L^p(B_2,W^{1,p}(N)) \hookrightarrow L^p(B_2,\cC^{0,\l'}(N)) 
\subset L^1(B_2,\cC^{0,\l'}(N)) .
$$
Now consider a bounded subset $\cK\subset W^{1,p}(B_2\times N)$.
The first embedding ensures that the functions $u\in\cK$, $u:B_2 \to L^p(N)$ are
equicontinuous. For some constant $C$
\begin{equation}  \label{equicont}
\| u(x) - u(y) \|_{L^p(N)} \leq C |x-y|^\l \qquad\forall x,y\in B_2, u\in\cK .
\end{equation}
The second embedding asserts that for some constant $C'$
\begin{equation} \label{bound}
\int_{B_2} \| u \|_{\cC^{0,\l'}(N)}  \leq  C'  \qquad\forall u\in\cK .
\end{equation}
In order to prove that $\cK\subset L^\infty(B_1,L^p(N))$ is precompact we now fix any
$\ep>0$ and show that $\cK$ can be covered by finitely many $\ep$-balls.

Let $J\in\cC^\infty(\R^m,[0,\infty))$ be such that $\supp J\subset B_1$ and $\int J = 1$. 
Then $J_\d(x):=\d^{-m}J(x/\d)$ are mollifiers for $\d>0$ with 
$\supp J_\d \subset B_\d$ and $\int J_\d = 1$.
Let $\d\leq 1$, then $J_\d * u |_{B_1}\in\cC^\infty(B_1,L^p(N))$ is welldefined.
Moreover, choose $\d>0$ sufficiently small such that for all $u\in\cK$
\begin{align*}
\bigl\| J_\d * u - u \bigr\|_{L^\infty(B_1,L^p(N))}
&= \sup_{x\in B_1} \Bigl\| \int_{B_\d} J_\d(y) \,
                                     ( u(x-y) - u(x) ) \,\rd^m y \Bigr\|_{L^p(N)} \\
&\leq \sup_{x\in B_1} \int_{B_\d} J_\d(y)\, C |y|^\l \,\rd^m y  \\
&\leq C \d^\l \;\leq\;{\textstyle\frac 14}\ep .
\end{align*}
Now it suffices to prove the precompactness of $\cK_\d:=\{J_\d*u \st u\in\cK\}$, then this
set can be covered by $\half\ep$-balls around $J_\d * u_i$ with $u_i\in\cK$ for 
$i=1,\ldots,I$
\footnote{
If a subset $K\subset (X,\rd)$ of a metric space is precompact, then for fixed $\ep>0$
one firstly finds $v_1,\ldots,v_I\in X$ such that for each $x\in K$ one has $\rd(x,v_i)\leq\ep$
for some $v_i$. For each $v_i$ choose one such $x_i\in K$, or simply drop $v_i$ if this does not
exist. Then $K$ is covered by $2\ep$-balls around the $x_i$: For each $x\in K$ one has
$\rd(x,x_i)\leq\rd(x,v_i)+\rd(v_i,x_i)$ for some $i=1,\ldots,I$.
}
and above estimate shows that $\cK$ is covered by the $\ep$-balls around the $u_i$.
Indeed, for each $u\in\cK$ one has $\|J_\d*u-J_\d*u_i\|_{L^\infty(B_1,L^p(N))}\leq\frac \ep 2$
for some $i=1,\ldots,I$ and thus
$$
\|u-u_i\|\;\leq\; \|u-J_\d*u\| + \|J_\d*u-J_\d*u_i\| + \|J_\d*u_i-u_i\|
\;\leq\; \ep .
$$
The precompactness of $\cK_\d\subset L^\infty(B_1,L^p(N))$ will follow from the
Arz\'ela-Ascoli theorem (see e.g.\ \cite[IX \pg 4]{Lang}).
Firstly, the smoothened functions $J_\d * u$ are still equicontinuous on $B_1$.
For all $u\in\cK$ and $x,y\in B_1$ use (\ref{equicont}) to obtain
\begin{align*}
\left\| (J_\d * u)(x) - (J_\d * u)(y) \right\|_{L^p(N)} 
&\leq \int_{B_\d} J_\d(z) \, \| u(x-z) - u(y-z) \|_{L^p(N)} \,\rd^m z \\
&\leq \int_{B_\d} J_\d(z) \, C | x - y |^\l \,\rd^m z 
\;=\;  C | x - y |^\l .
\end{align*}
Secondly, the $L^\infty$-norm of the smoothened functions is bounded by the $L^1$-norm
of the original ones, so for fixed $\d>0$ one obtains a uniform bound from (\ref{bound}) :
For all $u\in\cK$ and $x\in B_1$
\begin{align*}
\left\| (J_\d * u)(x) \right\|_{\cC^{0,\l'}(N)} 
&\leq \int_{B_2} J_\d(x-y) \, \| u(y) \|_{\cC^{0,\l'}(N)}  \,\rd^m y \\
&\leq C' \| J_\d \|_\infty .
\end{align*}
Now the embedding $\cC^{0,\l'}(N)\hookrightarrow L^p(N)$ is a standard compact Sobolev
embedding, so this shows that the subset $\{ (J_\d * u)(x) \st u\in\cK \}\subset L^p(N)$ 
is precompact for all $x\in B_1$.
Thus the Arz\'ela-Ascoli theorem asserts that \hbox{$\cK_\d \subset L^\infty(B_1,L^p(N))$} 
is compact, and this finishes the proof of the lemma.
\QED

In the proof of theorem~B, the weak Uhlenbeck compactness together with the local
slice theorem and this lemma will put us in the position to apply the following crucial
regularity theorem that also is the crucial point in the proof of theorem~A.
Here $(X,\t)$ is a $4$-manifold with a boundary space-time splitting as described in 
definition~\ref{def mfd} and in the beginning of this section.

\begin{thm} \label{reg}
For every compact subset $K\subset X$ there exists a compact submanifold 
$M\subset X$ such that $K \subset M$ and the following holds for all $p>2$.
\begin{enumerate}
\item
Suppose that $A\in\cA^{1,p}(M)$ solves the boundary value problem (\ref{bvp}).
Then $A|_K\in\cA(K)$ is smooth.
\item
Fix a metric $g_0$ that is compatible with $\t$ and a smooth connection $A_0\in\cA(M)$
such that $\t_i^*A_0|_{\{s\}\times\Si_i}\in\cL_i$ for all $s\in\cS_i$ and 
$i=1,\ldots,n$. 
Moreover, fix a compact neighbourhood $\cV=\bigcup_{i=1}^n \bar\t_{0,i}(\cU_i\times\Si_i)$
of $K\cap\pd X$. (Here $\bar\t_{0,i}$ denotes the extension of $\t_i$ given by the
geodesics of $g_0$.)
Then for every given constant $C_1$ there exist constants $\d>0$, $\d_k>0$, and 
$C_k$ for all $k\geq 2$ such that the following holds:

Fix $k\geq 2$ and let $g$ be a metric that is compatible with $\t$ and satisfies
$\|g-g_0\|_{\cC^{k+2}(M)}\leq\d_k$.
Suppose that $A\in\cA^{1,p}(M)$ solves the boundary value problem (\ref{bvp})
with respect to the metric $g$ and satisfies
\begin{align*}
\|A-A_0\|_{W^{1,p}(M)}&\leq C_1 , \\
\|\bar\t_{0,i}^*(A-A_0)|_{\Si_i}\|_{L^\infty(\cU_i,\cA^{0,p}(\Si_i))}&\leq\d 
\quad\forall i=1,\ldots,n.
\end{align*}
Then $A|_K\in\cA(K)$ is smooth by (i) and
$$
\|A-A_0\|_{W^{k,p}(K)}\leq C_k .
$$
\end{enumerate}
\end{thm}

We first give some preliminary results for the proof of theorem~\ref{reg}.
The interior regularity as well as the regularity of the $\cU_i$-components on a 
neighbourhood $\cU_i\times\Si_i$ of a boundary component $\cS_i\times\S_i$ will be a 
consequence of the following regularity result for Yang-Mills connections. 
The proof is similar to that of lemma~\ref{Hodge prop} and can be 
found in full detail in \cite[Proposition 10.5]{W}.
Here $M$ is a compact Riemannian manifold with boundary $\pd M$ 
and outer unit normal $\n$.
One then deals with two different spaces of test functions, 
\begin{align*}
\cC^\infty_\d(M,\cg) &:= \bigl\{ \p\in\cC^\infty(M,\cg) \st \p|_{\pd M} =0 \bigr\},  \\
\cC^\infty_\n(M,\cg) &:= \bigl\{ \p\in\cC^\infty(M,\cg) \st 
                                \tfrac{\pd\p}{\pd\n}\bigr|_{\pd M} =0 \bigr\}.
\end{align*}

\begin{prp}  \label{component reg}
Let $(M,g)$ be a compact Riemannian $4$-manifold. Fix a smooth reference 
connection $A_0\in\cA(M)$.
Let $X\in\G(\rT M)$ be a smooth vector field that is either perpendicular to 
the boundary, i.e.\ $X|_{\pd M}=h\cdot\n$ for some 
$h\in\cC^\infty(\pd M)$, or is tangential, i.e.\ 
$X|_{\pd M}\in\G(\rT\pd M)$.
In the first case let $\cT=\cC^\infty_\d(M,\cg)$, in the latter case let
$\cT=\cC^\infty_\n(M,\cg)$.
Moreover, let $N\subset\pd M$ be an open subset such that $X$ vanishes in a 
neighbourhood of $\pd M\setminus N\subset M$.
Let $1<p<\infty$ and $k\in\N$ be such that either $kp>4$
or $k=1$ and $2<p<4$.
In the first case let $q:=p$, in the latter case let $q:=\frac{4p}{8-p}$.
Then there exists a constant $C$ such that the following holds.

Let $A=A_0 + \a\in \cA^{k,p}(M)$ be a connection. 
Suppose that it satisfies 
\begin{equation}\label{eqn1}
\left\{\begin{array}{rl}
\rd_{A_0}^*\a &\hspace{-2.5mm}= 0, \\
*\a|_{\pd M} &\hspace{-2.5mm}= 0 \quad\text{on}\; N\subset\pd M ,
\end{array}\right.
\end{equation}
and that for all $1$-forms $\b=\p\cdot \i_X g$ with $\p\in\cT$
\begin{equation}\label{eqn2}
\int_M \la F_A \,,\, \rd_A \b \ra = 0 .
\end{equation}
Then $\a(X)\in W^{k+1,q}(M,\cg)$ and
$$
\|\a(X)\|_{W^{k+1,q}}\leq 
C\left( 1 + \|\a\|_{W^{k,p}} + \|\a\|_{W^{k,p}} ^3 \right) .
$$
Moreover, the constant $C$ can be chosen such that it depends continuously on 
the metric $g$ and the vector field $X$ with respect to the 
$\cC^{k+1}$-topology.
\end{prp}

\begin{rmk} \label{component reg rmk} \rm
In the case $k=1$ and $2<p<4$ the iteration of proposition~\ref{component reg}
also allows to obtain $W^{2,p}$-regularity and -estimates from initial 
$W^{1,p}$-regularity and -estimates.

Indeed, the Sobolev embedding $W^{2,q}\hookrightarrow W^{1,p'}$ holds with 
$p'=\frac{4q}{4-q}$ since $q<4$. Now as long as $p'<4$ one can iterate the 
proposition and Sobolev embedding to obtain regularity and estimates 
in $W^{1,p_i}$ with $p_0=p$ and
$$
p_{i+1}\;=\;\frac{4q_i}{4-q_i}\;=\;\frac{2p_i}{4-p_i} \;\geq\; \th p_i \;>\; p_i.
$$
Since $\th:=\frac 2{4-p}>1$ this sequence terminates after finitely many steps 
at some $p_N\geq 4$. Now in case $p_N>4$ the proposition even yields 
$W^{2,p_N}$-regularity and -estimates.
In case $p_N=4$ one only uses $W^{1,p_N}$ for some smaller $p'_N>\frac83$ in order
to conclude $W^{2,p'_{N+1}}$-regularity and -estimates for $p'_{N+1}>4$.

Similarly, in case $k=1$ and $p=4$ one only needs two steps to reach $W^{2,p'}$ for 
some $p'>4$.
\end{rmk}

The above proposition and remark can be used on all components of the connections 
in theorem~\ref{reg} except for the $\S$-components in small neighbourhoods 
$\cU\times\S$ of boundary components $\cS\times\S$.
For the regularity of their higher derivatives in $\S$-direction we shall use the following 
lemma. 
The crucial regularity of the derivatives in direction of $\cU$ of the $\S$-components will
then follow from the general regularity theory for Cauchy-Riemann equations in \cite{W Cauchy}.

\begin{lem} \label{Sigma reg}
Let $k\in\N_0$ and $1<p<\infty$.
Let $\Om$ be a compact manifold, let $\Si$ be a Riemann surface, and 
equip $\Om\times\Si$ with a product metric $g_\Om \oplus g$, where
$g=(g_x)_{x\in \Om}$ is a smooth family of metrics on $\Si$.
Then there exists a constant $C$ such that the following holds:

Suppose that $\a\in W^{k,p}(\Om\times\Si,\rT^*\Si)$ such that both $\rd_\Si\a$ 
and $\rd_\Si^*\a$ are of class $W^{k,p}$ on $\Om\times\S$.
Then $\nabla_\Si\a$ also is of class $W^{k,p}$ and one has the following estimate 
on $\Om\times\S$
$$
\|\nabla_\Si\a\|_{W^{k,p}} \leq C \bigl( \|\rd_\Si\a\|_{W^{k,p}} 
+ \|\rd_\Si^*\a\|_{W^{k,p}} +  \|\a\|_{W^{k,p}} \bigr).
$$
Here $\nabla_\S$ denotes the family of Levi-Civita connections on $\S$ that is given by
the family of metrics $g$.
Moreover, for every fixed family of metrics $g$ one finds a $\cC^k$-neighbourhood 
of metrics for which this estimate holds with a uniform constant $C$.
\end{lem}

\noindent
{\bf Proof of lemma \ref{Sigma reg}: } \\
We first prove this for $k=0$, i.e.\ suppose that 
$\a\in L^p(\Om\times\Si,\rT^*\Si)$ and that $\rd_\Si\a,\rd_\Si^*\a$
(defined as weak derivatives) are also of class $L^p$.
We introduce the following functions 
$$
f := \rd_\Si^*\a \;\in\; L^p(\Om\times\Si), \qquad
g := - *_{\scriptscriptstyle\Si}\rd_\Si\a \;\in\; L^p(\Om\times\Si),
$$
and choose sequences $f^\n,g^\n\in\cC^\infty(\Om\times\Si)$, and 
$\a^\n\in\cC^\infty(\Om\times\Si,\rT^*\Si)$ that converge to $f,g$, and $\a$ 
respectively in the $L^p$-norm.
Note that \hbox{$\int_\Om f = \int_\Om g = 0$} in $L^p(\S)$, so the $f^\n$ and $g^\n$ can be 
chosen such that their mean value over $\Om$ also vanishes for all $z\in\S$.
Then fix $z\in\Si$ and find
$\x^\n,\z^\n\in\cC^\infty(\Om\times\Si)$ such that
\[
\left\{\begin{aligned} 
\laplace_\Si \x^\n &= f^\n, \\
\x^\n(x,z) &= 0 \quad\forall x\in \Om ,
\end{aligned} \right.
\qquad\qquad
\left\{\begin{aligned} 
\laplace_\Si \z^\n &= g^\n , \\
\z^\n(x,z) &= 0 \quad\forall x\in \Om .
\end{aligned} \right.
\]
These solutions are uniquely determined since
$\laplace_\Si: W^{j+2,p}_z(\Si) \to W^{j,p}_m(\Si)$ is a bounded isomorphism
for every $j\in\N_0$ depending smoothly on the metric, i.e.\ on $x\in \Om$.
Here $W^{j,p}_m(\Si)$ denotes the space of $W^{j,p}$-functions with mean value zero
and $W^{j+2,p}_z(\Si)$ consists of those functions that vanish at $z\in\Si$.

Furthermore, let $\pi_x:\Om^1(\Si)\to h^1(\Si,g_x)$ be the projection of the smooth $1$-forms
to the harmonic part $h^1(\S)=\ker\laplace_\S=\ker\rd_\S\cap\ker\rd_\S^*$ with respect to the 
metric $g_x$ on $\Si$. 
Then $\pi$ is a family of bounded operators from $L^p(\Si,\rT^*\Si)$ to 
$W^{j,p}(\Si,\rT^*\Si)$ for any $j\in\N_0$, and it depends smoothly on $x\in\Om$.
So the harmonic part of $\tilde\a^\n$ is also smooth, 
$\pi\comp\tilde\a^\n \in\cC^\infty(\Om\times\Si,\rT^*\Si)$.
Now consider
$$
\a^\n := \rd_\Si\x^\n + *_{\scriptscriptstyle\Si}\rd_\Si\z^\n + \pi\comp\tilde\a^\n
\;\in\; \cC^\infty(\Om\times\Si,\rT^*\Si) .
$$
We will show that the sequence $\a^\n$ of $1$-forms converges to $\a$ in the 
$L^p$-norm and that moreover $\nabla_\Si\a^\n$ is an $L^p$-Cauchy sequence.
For that purpose we will use the following estimate. For all $1$-forms
$\b\in W^{1,p}(\S,\rT^*\S)$ abbreviating $\rd_\S=\rd$
\begin{align} 
\| \b \|_{W^{1,p}(\Si)} 
&\leq C \bigl( \| \rd^* \b \|_{L^p(\Si)} 
             + \| \rd \b \|_{L^p(\Si)} 
             + \| \pi(\b) \|_{W^{1,p}(\Si)} \bigr)  \nonumber\\
&\leq C \bigl( \| \rd^* \b \|_{L^p(\Si)} 
              + \| \rd \b \|_{L^p(\Si)} 
              + \| \b \|_{L^p(\Si)} \bigr) .  \label{Hodge}
\end{align}
Here and in the following $C$ denotes any finite constant that is uniform for all 
metrics $g_x$ on $\Si$ in a family of metrics that lies in a sufficiently small 
$\cC^k$-neighbourhood of a fixed family of metrics.
To prove (\ref{Hodge}) we use the Hodge decomposition $\b=\rd\x+*\rd\z+\pi(\b)$.
(See e.g.\ \cite[Theorem 6.8]{Wa} and recall that one can identify 2-forms on $\S$ with
functions via the Hodge $*$ operator.)
Here one chooses $\x,\z\in W^{2,p}_z(\Si)$ such that they solve $\laplace\x=\rd^*\b$ 
and $\laplace\z=*\rd\b$ respectively and concludes from proposition~\ref{Laplace reg}
for some uniform constant $C$
\begin{align*}
\|\rd\x\|_{W^{1,p}(\S)} &\leq \|\x\|_{W^{2,p}(\S)} \leq C \|\rd^*\b\|_{L^p(\S)} , \\
\|\text{$*\rd\z$}\|_{W^{1,p}(\S)} &\leq \|\z\|_{W^{2,p}(\S)} \leq C \|\rd\b\|_{L^p(\S)} .
\end{align*}
The second step in (\ref{Hodge}) moreover uses the fact that the projection to the harmonic
part is bounded as map $\pi:L^p(\Si,\rT^*\Si)\to W^{1,p}(\Si,\rT^*\Si)$.

Now consider $\a-\a^\n\in L^p(\Om\times\Si,\rT^*\Si)$. For almost all 
$x\in \Om$ we have $\a(x,\cdot)-\a^\n(x,\cdot)\in L^p(\Si,\rT^*\Si)$ as well as
$*\rd_\Si(\a(x,\cdot)-\a^\n(x,\cdot)) \in L^p(\Si)$ and 
$\rd_\Si^*(\a(x,\cdot)-\a^\n(x,\cdot)) \in L^p(\Si)$.
Then for these $x\in\Om$ one concludes from the Hodge decomposition that in fact
$\a(x,\cdot)-\a^\n(x,\cdot)\in W^{1,p}(\Si,\rT^*\Si)$. So we can apply 
(\ref{Hodge}) and integrate over $x\in \Om$ to obtain for all $\n\in\N$
\begin{align*}
&\| \a - \a^\n \|_{L^p(\Om\times\Si)}^p \\
&\leq \int_\Om \| \a(x,\cdot) - \a^\n(x,\cdot) \|_{L^p(\Si,g_x)}^p  \\
&\leq C \int_\Om \bigl( \| \rd_\Si^*(\a-\a^\n) \|_{L^p(\Si)}^p 
                  + \| \rd_\Si(\a-\a^\n) \|_{L^p(\Si)}^p 
                  + \| \pi(\a-\tilde\a^\n) \|_{W^{1,p}(\Si)}^p \bigr) \\
&\leq C \bigl( \| f-f^\n \|_{L^p(\Om\times\Si)}^p
              + \| g-g^\n \|_{L^p(\Om\times\Si)}^p 
              + \| \a-\tilde\a^\n \|_{L^p(\Om\times\Si)}^p \bigr) .
\end{align*}
In the last step we again used the continuity of $\pi$. This proves
the convergence $\a^\n\to\a$ in the $L^p$-norm, and hence 
$\nabla_\Si\a^\n\to\nabla_\Si\a$ in the distributional sense.
Next, we use (\ref{Hodge}) to estimate for all $\n\in\N$
\begin{align*}
\| \nabla_\Si\a^\n \|_{L^p(\Om\times\Si)}^p
&= \int_\Om \| \nabla_\Si\a^\n(x,\cdot) \|_{L^p(\Si,g_x)}^p   \\
&\leq C \int_\Om \bigl( \| \rd_\Si^*\a^\n \|_{L^p(\Si)} 
                    + \| \rd_\Si\a^\n \|_{L^p(\Si)} 
                    + \| \a^\n \|_{L^p(\Si)} \bigr)^p  \\
&\leq C \bigl( \| \rd_\Si^*\a^\n \|_{L^p(\Om\times\Si)}^p 
             + \| \rd_\Si\a^\n \|_{L^p(\Om\times\Si)}^p 
             + \| \a^\n \|_{L^p(\Om\times\Si)}^p \bigr) .
\end{align*}
Here one deals with $L^p$-convergent sequences
$\rd_\Si^*\a^\n=\laplace_\Si\x^\n=f^\n \to f=\rd_\Si^*\a$,
$-*\rd_\Si\a^\n=\laplace_\Si\z^\n=g^\n \to g= - *\rd_\Si\a$, and $\a^\n\to\a$. 
So $(\nabla_\Si\a^\n)_{\n\in\N}$ is uniformly bounded in $L^p(\Om\times\Si)$
and hence contains a weakly $L^p$-convergent subsequence. The limit is 
$\nabla_\Si\a$ since this already is the limit in the distributional sense. 
Thus we have proven the $L^p$-regularity of $\nabla_\Si\a$ on $\Om\times\S$, and moreover 
above estimate is preserved under the limit, which proves the lemma in the case $k=0$,
\begin{align*}
\| \nabla_\Si\a \|_{L^p(\Om\times\Si)}
&\leq \liminf_{\n\to\infty} \| \nabla_\Si\a^\n \|_{L^p(\Om\times\Si)}  \\
&\leq \liminf_{\n\to\infty} C \bigl( \| \rd_\Si^*\a^\n \|_{L^p(\Om\times\Si)}
                                   + \| \rd_\Si\a^\n \|_{L^p(\Om\times\Si)} 
                                   + \| \a^\n \|_{L^p(\Om\times\Si)} \bigr) \\
&= C \bigl( \| \rd_\Si^*\a \|_{L^p(\Om\times\Si)}
          + \| \rd_\Si  \a \|_{L^p(\Om\times\Si)} 
          + \|         \a \|_{L^p(\Om\times\Si)} \bigr) .
\end{align*}
In the case $k\geq 1$ one can now use the previous result to prove the lemma.
Let $\a\in W^{k,p}(\Om\times\Si,\rT^*\Si)$ and suppose that 
$\rd_\Si\a,\rd_\Si^*\a$ are of class $W^{k,p}$.
We denote by $\nabla$ the covariant derivative on $\Om\times\Si$. Then we have to
show that $\nabla^k\nabla_\Si\a$ is of class $L^p$.
So let $X_1,\ldots,X_k$ be smooth vector fields on $\Om\times\Si$ and introduce
$$
\tilde\a:=\nabla_{X_1}\ldots\nabla_{X_k}\a
\;\in\; L^p(\Om\times\Si,\rT^*\Si) .
$$
Both $\rd_\Si\tilde\a$ and $\rd_\Si^*\tilde\a$ are of class $L^p$ since
\begin{align*}
\rd_\Si\tilde\a 
& = [\rd_\Si,\nabla_{X_1}\ldots\nabla_{X_k}] \a +
    \nabla_{X_1}\ldots\nabla_{X_k}\rd_\Si\a ,\\
\rd_\Si^*\tilde\a 
& = [\rd_\Si^*,\nabla_{X_1}\ldots\nabla_{X_k}] \a +
    \nabla_{X_1}\ldots\nabla_{X_k}\rd_\Si^*\a .
\end{align*}
So the result for $k=0$ implies that $\nabla_\Si\tilde\a$ is of class $L^p$,
hence $\nabla^k\nabla_\Si\a$ also is of class $L^p$ since for all smooth
vector fields
$$
\nabla_{X_1}\ldots\nabla_{X_k}\nabla_\Si\a 
= [\nabla_\Si,\nabla_{X_1}\ldots\nabla_{X_k}] \a + \nabla_\Si\tilde\a .
$$
With the same argument -- using coordinate vector fields $X_i$ and cutting them
off -- one obtains the estimate
$$
\|\nabla^k\nabla_\Si\a\|_{L^p(\Om\times\Si)}
\leq C \bigl( \|\nabla^k\rd_\Si^*\a\|_{L^p(\Om\times\Si)}
            + \|\nabla^k\rd_\Si\a\|_{L^p(\Om\times\Si)}
            + \|\a\|_{W^{k,p}(\Om\times\Si)} \bigr) .
$$
Now this proves the lemma,
\begin{align*}
\|\nabla_\Si\a\|_{W^{k,p}(\Om\times\Si)}
&\leq \|\nabla_\Si\a\|_{W^{k-1,p}(\Om\times\Si)}
    + \|\nabla^k\nabla_\Si\a\|_{L^p(\Om\times\Si)}  \\
&\leq C \bigl( \|\rd_\Si^*\a\|_{W^{k,p}(\Om\times\Si)}
             + \|\rd_\Si\a\|_{W^{k,p}(\Om\times\Si)}
             + \|\a\|_{W^{k,p}(\Om\times\Si)} \bigr) .
\end{align*}
\QED \\

\noindent
{\bf Proof of theorem \ref{reg} : } \\
Recall that a neighbourhood of the boundary $\pd X$ is covered by embeddings
\hbox{$\bar\t_{0,i}:\cU_i\times\Si_i\hookrightarrow X$} such that 
$\bar\t_{0,i}^*g_0=\ds^2+\dt^2+g_{0;s,t}$. (In the case (i) we put $g_0:=g$.)
Since $K\subset X$ is compact one can cover it by a compact subset 
$K_{\text{int}} \subset {\rm int}\,X$ and 
$K_{\text{bdy}}:=\bigcup_{i=1}^n \bar\t_{0,i}(I_{0,i}\times[0,\d_0]\times\Si_i)$
for some $\d_0>0$ and $I_{0,i}\subset\cS_i$ that are either compact intervals in 
$\R$ or equal to $S^1$.
Moreover, one can ensure that $K_{\text{bdy}}\subset{\rm int}\,\cV$ lies in the interior of
the fixed neighbourhood of $K\cap\pd X$.
Since $X$ is exhausted by the compact submanifolds $X_k$ one then finds
$M:=X_k\subset X$ such that both $K_{\text{bdy}}$ and $K_{\text{int}}$ are 
contained in the interior of $M$ (and thus also $K\subset M$).
Now let $A\in\cA^{1,p}(M)$ be a solution of the boundary value problem 
(\ref{bvp}) with respect to a metric $g$ that is compatible with $\t$. 
Then we will prove its regularity and the corresponding estimates
in the interior case on $K_{\text{int}}$ and in the boundary case on $K_{\text{bdy}}$
separately.\\

\noindent
\underline{\bf Interior case :}\\
Firstly, since $K_{\text{int}} \subset {\rm int}\,M$ and 
$K_{\text{int}}\subset {\rm int}\,X = X\setminus\pd X$ we find a sequence of 
compact submanifolds $M_k\subset {\rm int}\,X$ such that
$K_{\text{int}}\subset M_{k+1} \subset{\rm int}\,M_k\subset M$ for all $k\in\N$. 
We will prove inductively $A|_{M_k}\in\cA^{k,p}(M_k)$ for all $k\in\N$
and thus $A|_{K_{\text{int}}}\in\cA(K_{\text{int}})$ is smooth.
Moreover, we inductively find constants $C_k,\d_k>0$ such that the 
additional assumptions of (ii) in the theorem imply
\begin{equation}\label{Ck est}
\|A-A_0\|_{W^{k,p}(M_k)}\leq C_k .
\end{equation}
Here we use the fixed smooth metric $g_0$ to define the Sobolev norms -- 
for a sufficiently small $\cC^k$-neighbourhood of metrics, the Sobolev norms are equivalent 
with a uniform constant independent of the metric.
Moreover, recall that the reference connection $A_0$ is smooth.

To start the induction we observe that this regularity and estimate are satisfied for 
$k=1$ by assumption. 
For the induction step assume this regularity and estimate to hold for some $k\in\N$.
Then we will use proposition~\ref{component reg} on $A|_{M_k}\in\cA^{k,p}(M_k)$
to deduce the regularity and estimate on $M_{k+1}$.

Every coordinate vector field on $M_{k+1}$ can be extended to a vector field 
$X$ on $M_k$ that vanishes near the boundary $\pd M_k$. 
So it suffices to consider such vector fields, i.e.\ use $N=\emptyset$ in the proposition.
Then $\a:=A-A_0$ satisfies the assumption (\ref{eqn1}).
For the weak equation (\ref{eqn2}) we calculate for all 
$\b=\p\cdot \i_X g$ with $\p\in\cT=\cC^\infty_\d(M_k,\cg)$
$$
- \int_{M_k} \la F_A \,,\, \rd_A\b \ra 
\;=\;  \int_{M_k} \la \rd_A(\p\cdot \i_X g) \wedge F_A \ra 
\;=\;  \int_{\pd M_k} \la \p\cdot \i_X g \wedge F_A \ra 
\;=\; 0 .
$$
We have used Stokes' theorem while approximating $A$ by smooth 
connections $\tA$, for which the Bianchi identity $\rd_\tA F_\tA =0$ holds.
Now proposition~\ref{component reg} and remark~\ref{component reg rmk} imply that
$A|_{M_{k+1}}\in\cA^{k+1,p}(M_{k+1})$.
In the case (ii) of the theorem the proposition moreover provides $\d_{k+1}>0$ and a uniform 
constant $C$ for all metrics $g$ with $\|g-g_0\|_{\cC^{k+1}(M_k)}\leq\d_{k+1}$ 
such that the following holds: If (\ref{Ck est}) holds for some constant $C_k$, 
then 
\begin{align*}
\|A-A_0\|_{W^{k+1,p}(M_{k+1})} 
&\leq C\left( 1 + \|A-A_0\|_{W^{k,p}(M_k)} 
                + \|A-A_0\|_{W^{k,p}(M_k)} ^3 \right) \\
&\leq C\left( 1 + C_k + C_k^3 \right) \;=:\, C_{k+1} . 
\end{align*}
Here we have used the fact that the Sobolev norm of a $1$-form is equivalent
to an expression in terms of the Sobolev norms of its components in the 
coordinate charts.
In case $k=1$ and $p\leq 4$, this uniform bound is not found directly but after finitely
many iterations of proposition \ref{component reg} that give estimates on manifolds 
$N_1=M_1$ and $M_2\subset N_{i+1} \subset {\rm int}\,N_i$. In each step one chooses a smaller
$\d_2>0$ and a bigger $C_2$. This iteration uses the same Sobolev embeddings as 
remark~\ref{component reg rmk}.
This proves the induction step on the interior part $K_{\text{int}}$.\\

\noindent
\underline{\bf Boundary case :}\\
It remains to prove the regularity and estimates on $K_{\text{bdy}}$ near the 
boundary.
So consider a single boundary component $K':=\bar\t_0(I_0\times[0,\d_0]\times\Si)$.
We identify $I_0=S^1\cong\R/\Z$ or shift the compact interval such that 
$I_0=[-r_0,r_0]$ and hence $K'=\bar\t_0([-r_0,r_0]\times[0,\d_0]\times\Si)$ for some
$r_0>0$.
Since $K_{\text{bdy}}$ (and thus also $K'$) lies in the interior of $M$ as well as $\cV$, 
one then finds $R_0>r_0$ and $\D_0>\d_0$ such that 
$\bar\t_0([-R_0,R_0]\times[0,\D_0]\times\Si)\subset M\cap\cV$.
Here $\bar\t_0$ is the embedding that brings the metric $g_0$ into the standard 
form $\ds^2+\dt^2+g_{0;s,t}$. 
A different metric $g$ compatible with $\t$ defines a different embedding 
$\bar\t$ such that $\bar\t^*g=\ds^2+\dt^2+g_{s,t}$.
However, if $g$ is sufficiently $\cC^1$-close to $g_0$, then the geodesics are
$\cC^0$-close and hence $\bar\t$ is $\cC^0$-close to $\bar\t_0$. 
(These embeddings are fixed for $t=0$, and for $t>0$ given by the normal geodesics.)
Thus for a sufficiently small choice of $\d_2>0$ one finds $R>r>0$ and $\D>\d>0$ 
such that for all $\t$-compatible metrics $g$ in the $\d_2$-ball around $g_0$
$$
K'\subset\bar\t([-r,r]\times[0,\d]\times\Si)
\qquad\text{and}\qquad 
\bar\t([-R,R]\times[0,\D]\times\Si) \subset M\cap\cV .
$$
(In the case (i) this holds with $r_0,\d_0,R_0$, and $\D_0$ for the fixed metric $g=g_0$.)
We will prove the regularity and estimates for $\bar\t^*A$ on $[-r,r]\times[0,\d]\times\Si$.
This suffices because for $\cC^{k+2}$-close metrics the embedding $\bar\t$ 
will be $\cC^{k+1}$-close to the fixed $\bar\t_0$, so that one obtains uniform constants in the 
estimates between the $W^{k,p}$-norms of $A$ and $\bar\t^*A$.
Furthermore, the families $g_{s,t}$ of metrics on $\Si$ will be $\cC^k$-close to $g_{0;s,t}$ for 
$(s,t)\in[-R,R]\times[0,\D]$ if $\d_k$ is chosen sufficiently small.
Now choose compact submanifolds $\Om_k\subset\H:=\{(s,t)\in\R^2\st t\geq 0\}$ 
such that for all $k\in\N$
$$
[-r,r]\times[0,\d]
\;\subset\; \Om_{k+1}
\;\subset\; {\rm int}\,\Om_k
\;\subset\; [-R,R]\times[0,\D] .
$$
We will prove the theorem by establishing the regularity and estimates for 
$\bar\t^*A$ on the $\Om_k\times\Si$ in Sobolev spaces of increasing differentiability.
We distinguish the cases $p>4$ and $4\geq p>2$. In case $p>4$ one uses the following induction.\\

\noindent
{\bf I)}\; {\it
Let $p>2$ and suppose that $A\in\cA^{1,2p}(M)$ solves (\ref{bvp}).
Then we will prove inductively that
$\bar\t^*A|_{\Om_k\times\Si}\in\cA^{k,q}(\Om_k\times\Si)$ for all $k\in\N$ and with 
$q=p$ or $q=2p$ according to whether $k\geq 2$ or $k=1$.
Moreover, we will find a constant $\d>0$ and constants $C_k,\d_k>0$ for all $k\geq 2$ 
such that the following holds:

If in addition $\|g-g_0\|_{\cC^{k+2}(M)}\leq\d_k$ and 
\begin{align*}
\|A-A_0\|_{W^{1,2p}(M)}&\leq C_1 , \\
\|\bar\t_0^*(A-A_0)|_\S\|_{L^\infty(\cU,\cA^{0,p}(\S))}&\leq\d ,
\end{align*}
then for all $k\in\N$
$$
\|\bar\t^*(A-A_0)\|_{W^{k,q}(\Om_k\times\S)}\leq C_k .
$$
}

This is sufficient to conclude the theorem in case $p>4$ as follows. 
One uses I) with $p$ replaced by $\half p$ to obtain regularity and estimates of $A-A_0$ in 
$\cA^{1,p}(\Om_1\times\Si)$, $\cA^{2,\frac p2}(\Om_2\times\Si)$, and 
$\cA^{k,\frac p2}(\Om_k\times\Si)$ for all $k\geq 3$.
Recall that the component $K'$ of $K_{\rm bdy}$ is contained in each $\bar\t(\Om_k\times\Si)$.
In addition, one has the Sobolev embeddings 
$W^{k+1,\frac p2}\hookrightarrow W^{k,p}\hookrightarrow\cC^{k-1}$ on the compact $4$-manifolds
$\Om_{k+1}\times\S)$, c.f.\ \cite[Theorem 5.4]{Adams}.
So this proves the regularity and estimates on $K_{\rm bdy}$.

In the case $4\geq p>2$ a preliminary iteration is required in order to achieve the
regularity and estimates that are assumed in I). In contrast to I) the iteration is in $p$ 
instead of $k$.\\

\noindent
{\bf II)}\; {\it
Let $4\geq p>2$ and suppose that $A\in\cA^{1,p}(M)$ solves (\ref{bvp}).
Then we will prove inductively that
$\bar\t^*A|_{\Om_j\times\Si}\in\cA^{1,p_j}(\Om_j\times\Si)$ for a sequence
$(p_j)$ with $p_1=p$ and $p_{j+1}=\th(p_j)\cdot p_j$, where 
$\th:(2,4]\to (1,\frac{17}{16}]$ is monotonely increasing
and thus the sequence terminates with $p_N >4$ for some $N\in\N$.

Moreover, we will find constants $\d>0$ and constants $C_{1,j},\d_{1,j}>0$ for $j=2,\ldots,N$
such that the following holds:

If for some $j=1,\ldots,N$ in addition $\|g-g_0\|_{\cC^3(M)}\leq\d_{1,j}$ and 
\begin{align*}
\|A-A_0\|_{W^{1,p}(M)}&\leq C_1 , \\
\|\bar\t_0^*(A-A_0)|_\S\|_{L^\infty(\cU,\cA^{0,p}(\S))}&\leq\d , 
\end{align*}
then
$$
\|\bar\t^*(A-A_0)\|_{W^{1,p_j}(\Om_j\times\Si)}\leq C_{1,j} .
$$
}

Assuming I) and II) we first prove the theorem for the case $4\geq p >2$.
After finitely many steps the iteration of II) gives regularity and estimates in 
$\cA^{1,p_N}(\Om_N\times\Si)$ with $p_N>4$ and under the assumption 
$\|g-g_0\|_{\cC^3(M)}\leq\d_{1,N}$ on the metric.
Now if necessary decrease $p_N$ slightly such that $2p\geq p_N>4$, then one
still has $\cA^{1,p_N}$-regularity and estimates on all components of $K_{\rm bdy}$
as well as on $K_{\rm int}$ (from the previous argument on the interior).
Thus the assumptions of I) are satisfied with $p$ replaced by $\half p_N$ and 
$C_1$ replaced by a combination of $C_{1,N}$ and a constant from the interior iteration
(both of which only depend on $C_1$). One just has to choose $\d_2\leq\d_{1,N}$ and
choose the $\d>0$ in I) smaller than the $\d>0$ from II).
Then the iteration in I) gives regularity and estimates of $A-A_0$ in 
$\cA^{k,\frac 12 p_N}(\Om_k\times\Si)$ for all $k\geq 2$.
This proves the theorem in case $2<p\leq 4$ due to the Sobolev embeddings
$W^{k+1,\frac 12 p_N}\hookrightarrow W^{k,p}\hookrightarrow\cC^{k-2}$.
So it remains to establish I) and II).\\

\noindent
\underline{\bf Proof of I):}\\
The start of the induction $k=1$ is true by assumption
(after replacing $C_1$ by a larger constant to make up for the effect of $\bar\t^*$).
For the induction step assume that the claimed regularity and estimates hold for some 
$k\in\N$ and consider the following decomposition of the connection $A$ and its curvature:
\begin{align}
\bar\t^* A &= \P\,\ds + \Psi\,\dt + B , \nonumber\\
\bar\t^* F_A &= F_B + (\rd_B\P - \pd_s B)\wedge\ds + (\rd_B\Psi - \pd_t B)\wedge\dt 
\label{coos}\\
&\quad   + (\pd_s\Psi  - \pd_t\P + [\P,\Psi] )\,\ds\wedge\dt .  \nonumber
\end{align}
Here $\P,\Psi \in W^{k,q}(\Om_k\times\Si,\cg)$,
and $B\in W^{k,q}(\Om_k\times\Si,\rT^*\Si\otimes\cg)$ is a $2$-parameter 
family of $1$-forms on $\Si$. 
Choose a further compact submanifold $\Om\subset{\rm int}\,\Om_k$ such that
$\Om_{k+1}\subset{\rm int}\,\Om$.
Now we shall use proposition~\ref{component reg} to deduce the higher regularity 
of $\P$ and $\Psi$ on $\Om\times\Si$.
For this purpose one has to extend the vector fields $\pd_s$ and $\pd_t$ on $\Om\times\Si$ to
different vector fields on $\Om_k\times\Si$, both denoted by $X$, and verify the assumptions
(\ref{eqn1}) and (\ref{eqn2}) of proposition~\ref{component reg}.
These extensions will be chosen such that they vanish in a 
neighbourhood of $(\pd\Om_k\setminus\pd\H)\times\Si$. Then $\a:=\bar\t^*(A-A_0)$ 
satisfies (\ref{eqn1}) on $M=\bar\t(\Om_k\times\Si)$ with 
$N=\bar\t((\pd\Om_k\cap\pd\H)\times\S)$. 

Choose a cutoff function $h\in\cC^\infty(\Om_k,[0,1])$ that equals $1$ on 
$\Om$ and vanishes in a neighbourhood of $\pd\Om_k\setminus\pd\H$.
Then firstly, $X:=h\pd_t$ is a vector field as required that is perpendicular
to the boundary $\pd\Om_k\times\Si$. 
For this type of vector field we have to check the assumption (\ref{eqn2}) for all 
$\b=\p h \cdot \dt$ with $\p\in\cC^\infty_\d(\Om_k\times\Si,\cg)$. Note that
$\bar\t_*\b = (\p\cdot h)\comp\bar\t^{-1}\cdot \i_{(\bar\t_*\pd_t)}g$ can be 
trivially extended to $M$ and then vanishes when restricted to $\pd M$.
So we can use partial integration as in the interior case to obtain
\begin{align*}
\int_{\Om_k\times\Si} \la F_{\bar\t^*A} \,,\, \rd_{\bar\t^*A}\b \ra 
\;=\; \int_M \la F_A \,,\, \rd_A \bar\t_*\b \ra 
\;=\; - \int_{\pd M} \la \bar\t_*\b \wedge F_A \ra 
\;=\; 0 .
\end{align*}
Secondly, $X:=h\pd_s$ also vanishes in a neighbourhood of 
$(\pd\Om_k\setminus\pd\H)\times\Si$ and is tangential to the boundary 
$\pd\Om_k\times\Si$. So we have to verify (\ref{eqn2}) for all
$\b=\p h\cdot\ds$ with $\p\in\cT=\cC^\infty_\n(\Om_k\times\Si,\cg)$.
Again, $\bar\t_*\b$ extends trivially to $M$. Then the partial integration yields
\begin{align*}
\int_{\Om_k\times\Si} \la F_{\bar\t^*A} \,,\, \rd_{\bar\t^*A}\b \ra 
&= - \int_{\bar\t^{-1}(\pd M)} \la \b\wedge \bar\t^* F_A \ra \\
&= - \int_{(\Om_k\cap\pd\H)\times\Si} \la \p h \cdot \ds \wedge F_B \ra 
\;=\; 0 .
\end{align*}
The last step uses the fact that 
$B(s,0)=\t^*A|_{\{s\}\times\Si}\in\cL\subset\cA^{0,p}_{\rm flat}(\Si)$, and hence
$F_B$ vanishes on $\pd\H\times\Si$.
However, we have to approximate $A$ by smooth connections in order that 
Stokes' theorem holds and $F_B$ is welldefined.
So this calculation crucially uses the fact that a $W^{1,p}$-connection with boundary values
in the Lagrangian submanifold $\cL$ can be $W^{1,p}$-approximated by smooth
connections with boundary values in $\cL\cap\cA(\Si)$. This was proven in 
\cite[Corollary~4.5]{W Cauchy}.
So we have verified the assumptions of proposition~\ref{component reg} for 
both $\P=\bar\t^*A(\pd_s)$ and $\Psi=\bar\t^*A(\pd_t)$ and thus can deduce 
$\P,\Psi\in W^{k+1,q}(\Om\times\Si)$.
Moreover, under the additional assumptions of (ii) in the theorem we have the estimates
\begin{align}
\|\P-\P_0\|_{W^{k+1,q}(\Om\times\Si)}
&\leq C_s\left( 1 + C_k + C_k^3 \right) \;=:\, C_{k+1}^s ,\nonumber\\
\|\Psi-\Psi_0\|_{W^{k+1,q}(\Om\times\Si)}
&\leq C_t\left( 1 + C_k + C_k^3 \right) \;=:\, C_{k+1}^t.  \label{est2}
\end{align}
The constants $C_s$ and $C_t$ are uniform for all metrics in 
some small $\cC^{k+1}$-neighbourhood of $g_{0;s,t}$, so by a possibly smaller choice of 
$\d_{k+1}>0$ they become independent of $g_{s,t}$.
Note that in the above estimates we also have decomposed the reference connection in the 
tubular neighbourhood coordinates,
$\bar\t^*A_0 = \P_0\,\ds + \Psi_0\,\dt + B_0$.

It remains to consider the $\Si$-component $B$ in the tubular neighbourhood.
The boundary value problem (\ref{bvp}) becomes in the coordinates (\ref{coos}) 
\begin{equation}\label{bvp comp}
\left\{\begin{aligned} 
\rd_{B_0}^*(B-B_0) &= \nabla_s(\P-\P_0) + \nabla_t(\Psi-\Psi_0) , \\
*F_B &= \pd_t\P - \pd_s\Psi  + [\Psi,\P] ,\\
\pd_s B + * \pd_t B &= \rd_B\P + * \rd_B\Psi , \\
\Psi(s,0) - \Psi_0(s,0) &= 0\quad\forall (s,0)\in\pd\Om_k , \\
B(s,0) &\in\cL \quad\forall (s,0)\in\pd\Om_k . \\
\end{aligned}\right.
\end{equation}
Here $\rd_B$ is the exterior derivative on $\Si$ that is associated with the 
connection~$B$, $\rd_{B_0}^*$ is the coderivative associated with $B_0$, 
$*$ is the Hodge operator on $\Si$ with respect to the metric $g_{s,t}$, and 
$\nabla_s \P := \pd_s \P + [\P_0,\P]$,
$\nabla_t \P := \pd_t \P + [\Psi_0,\P]$.
We rewrite the first two equations in (\ref{bvp comp}) as a system
of differential equations for $\a:=B-B_0$ on $\Si$. For each $(s,t)\in \Om_k$
\begin{equation} \label{Sigma eqn}
\rd_\Si^*\a(s,t) = \x(s,t) , \qquad\qquad
\rd_\Si \a(s,t) = *\z(s,t) .
\end{equation}
Here we have abbreviated
\begin{align*}
\x &=     *[B_0\wedge *(B-B_0)] 
        + \nabla_s(\P-\P_0) + \nabla_t(\Psi-\Psi_0) , \\
\z &= -*\rd_\Si B_0 - * \half[B\wedge B] + \pd_t\P - \pd_s\Psi  + [\Psi,\P] .
\end{align*}
These are both functions in $W^{k,q}(\Om\times\Si,\cg)$ due to the smoothness 
of $A_0$ and the previously established regularity of $\P$ and $\Psi$.
(This uses the Sobolev embedding $W^{k,q}\cdot W^{k,q}\hookrightarrow W^{k,q}$
due to $W^{k,q}\hookrightarrow L^\infty$.)
So lemma~\ref{Sigma reg} asserts that $\nabla_\Si (B-B_0)$ is of class
$W^{k,q}$ on $\Om\times\Si$, and under the assumptions of (ii) in the theorem we obtain the 
estimate
\begin{align} 
&\|\nabla_\Si (B-B_0)\|_{W^{k,q}(\Om\times\Si)} \nonumber\\
&\leq  C \bigl( \|\x\|_{W^{k,q}} 
+ \|\z\|_{W^{k,q}} +  \|B-B_0\|_{W^{k,q}} \bigr) \nonumber\\
&\leq  C \bigl( 1 + \| B-B_0 \|_{W^{k,q}}
          + \|\P-\P_0\|_{W^{k+1,q}}
          + \|\Psi-\Psi_0\|_{W^{k+1,q}}  \nonumber\\
&\qquad\quad\,
          + \| B-B_0 \|_{W^{k,q}}^2
          + \|\P-\P_0\|_{W^{k,q}}\|\Psi-\Psi_0\|_{W^{k,q}} \bigr) \nonumber\\
&\leq  C \bigl( 1 + C_k + C^s_{k+1} + C^t_{k+1} + C_k^2 \bigr) 
 \;=:\, C^\Si_{k+1} .
\label{est3}
\end{align}
Here $C$ denotes any constant that is uniform for all metrics in a 
$\cC^{k+1}$-neighbourhood of the fixed $g_{0;s,t}$, so this might again require
a smaller choice of $\d_{k+1}>0$ in order that the constant $C^\Si_{k+1}$
becomes independent of the metric $g_{s,t}$.

Now we have established the regularity and estimate for all derivatives of $B$ 
of order $k+1$ containing at least one derivative in $\Si$-direction.
(Note that in the case $k=1$ we even have $L^q$-regularity with $q=2p$ where only
$L^p$-regularity was claimed. This additional regularity will be essential
for the following argument.)
It remains to consider the pure $s$- and $t$- derivatives of $B$ and 
establish the $L^p$-regularity and -estimate for $\nabla_\H^{k+1}B$ on
$\Om_{k+1}\times\Si$, where $\nabla_\H$ is the standard covariant derivative on $\H$ with
respect to the metric $\ds^2+\dt^2$.
The reason for this regularity, as we shall show, is the fact that 
$B\in W^{k,q}(\Om,\cA^{0,p}(\Si))$ satisfies a Cauchy-Riemann equation with 
Lagrangian boundary conditions,
\begin{equation}\label{Lagrangian bvp}
\left\{\begin{aligned} 
\pd_s B + * \pd_t B &= G , \\
B(s,0) \in\cL &\quad\forall (s,0)\in\pd\Om .
\end{aligned}\right.
\end{equation}
The inhomogeneous term is
$$
G := \rd_B\P + * \rd_B\Psi \;\in\; W^{k,q}(\Om,\cA^{0,p}(\Si)) .
$$
Here one uses the fact that
$W^{k,q}(\Om\times\Si,\rT^*\Si\otimes\cg)\subset W^{k,q}(\Om,\cA^{0,p}(\Si))$
since the smooth $1$-forms are dense in both spaces and the norm on the second space is 
weaker than the $W^{k,q}$-norm on $\Om\times\Si$, c.f.\ \cite[Lemma~2.2]{W Cauchy}.

Now one has to apply the regularity result \cite[Theorem~1.2]{W Cauchy}
for the Cauchy-Riemann equation on the complex Banach space $\cA^{0,p}(\S)$. As reference 
complex structure $J_0$ we use the Hodge $*$ operator on $\Si$ with respect to the 
fixed family of metrics $g_{0;s,t}$ on $\Si$ (that varies smoothly with $(s,t)\in\Om$).
The smooth family $J$ of complex structures in the equation is given by the Hodge operators 
with respect to the metrics $g_{s,t}$.
The Lagrangian submanifold $\cL\subset\cA^{0,p}(\Si)$ is totally real with respect to any
Hodge operator, and it is modelled on a closed subspace of $L^p(\Si,\R^n)$ for some $n\in\N$
(see \cite[Lemma~4.2, Corollary~4.4]{W Cauchy} ) .
In the case (ii) of the theorem moreover a family of connections 
$B_0\in\cC^\infty(\Om,\cA(\Si))$ is given such that $B_0(s,0)\in\cL$ for all $(s,0)\in\pd\Om$
and $B$ satisfies
\begin{align*}
\|B-B_0\|_{L^\infty(\Om,\cA^{0,p}(\Si))} 
&=\|\bar\t^*(A-A_0)|_\Si\|_{L^\infty(\Om,\cA^{0,p}(\Si))} \\
&\leq C \|\bar\t_0^*(A-A_0)|_\Si\|_{L^\infty(\cU,\cA^{0,p}(\Si))} \;\leq\; C\d .
\end{align*}
Here one uses the fact that $\bar\t(\Om\times\Si)\subset\bar\t_0(\cU\times\Si)$
lies in a component of the fixed neighbourhood $\cV$ of $K\cap\pd X$.
The assumption of closeness to $A_0$ in $\cA^{0,p}(\Si)$ was formulated for
$\bar\t_0^*(A-A_0)|_\Si$. However, for a metric $g$ in a sufficiently small 
$\cC^2$-neighbourhood of the fixed metric $g_0$ the extensions $\bar\t$ and
$\bar\t_0$ are $\cC^1$-close and one obtains the above estimate with a constant $C$
independent of the metric.
So $B\in W^{k,q}(\Om,\cA^{0,p}(\Si))$ satisfies the assumptions of \cite[Theorem~1.2]{W Cauchy}
if $\d>0$ is chosen sufficiently small. (Note that this choice is independent of $k\in\N$.)

Now \cite[Theorem~1.2]{W Cauchy} asserts $B\in W^{k+1,p}(\Om_{k+1},\cA^{0,p}(\Si))$.
By  \cite[Lemma~2.2]{W Cauchy} this also proves
$\nabla_\H^{k+1}B \in L^p(\Om_{k+1},\cA^{0,p}(\Si))
                     =L^p(\Om_{k+1}\times\Si,\rT^*\Si\otimes\cg)$,
and this finishes the induction step 
$\bar\t^*A|_{\Om_{k+1}\times\Si}\in\cA^{k+1,p}(\Om_{k+1}\times\Si)$
for the regularity near the boundary.
The induction step for the estimate in case (ii) of the theorem now follows from the estimate
from \cite[Theorem~1.2]{W Cauchy},
\begin{align}
&\|\nabla_\H^{k+1} (B-B_0) \|_{L^{p}(\Om_{k+1}\times\Si)} \nonumber\\
& \leq \| B-B_0 \|_{W^{k+1,p}(\Om_{k+1},\cA^{0,p}(\Si))} \nonumber\\
& \leq C \bigl( 1 + \|G\|_{W^{k,q}(\Om,\cA^{0,p}(\Si))} 
                  + \|B-B_0\|_{W^{k,q}(\Om,\cA^{0,p}(\Si))} \bigr) \nonumber\\
& \leq C \bigl( 1 + C_k + C_k^2 + C_{k+1}^s +  C_{k+1}^t \bigr) \label{CH}
  \;=:\, C^\H_{k+1} .
\end{align}
Here the constant from \cite[Theorem~1.2]{W Cauchy} 
is uniform for a sufficiently small $\cC^{k+1}$-neighbourhood of complex structures. 
In this case, these are the families of Hodge operators on $\Si$ that depend on the 
metric $g_{s,t}$. Thus for sufficiently small $\d_{k+1}>0$ that constant (and also the 
further Sobolev constants that come into the estimate) becomes independent of the metric.
The final constant $C_{k+1}$ then results from all the separate estimates,
see the decomposition (\ref{coos}) and the estimates in (\ref{est2}), (\ref{est3}), and
(\ref{CH}),
$$
\|\bar\t^*(A-A_0)\|_{W^{k+1,p}(\Om_{k+1}\times\Si)}
\;\leq\; C_k + C_{k+1}^s + C_{k+1}^t + C_{k+1}^\Si + C_{k+1}^\H.
$$

\medskip

\noindent
\underline{\bf Proof of II):}\\
Except for the higher differentiability of $B$ in direction of $\H$ this iteration
works by the same decomposition and equations as in I).
The start of the induction $k=1$ is given by assumption.
For the induction step assume that the claimed  $W^{1,p_k}$-regularity and -estimates hold 
for some $k\in\N$ with $p_k\leq 4$. Then proposition~\ref{component reg} gives
$\P,\Psi\in W^{2,q_k}(\Om\times\Si)$ with corresponding estimates and 
$$
q_k = \left\{ \begin{array}{cl}
          \frac{4p_k}{8-p_k}  &\text{if}\; p_k<4 ,\\
            \\
           3 &\text{if}\; p_k = 4 .
      \end{array}\right.
$$
(In the case $p_k=4$ one applies the proposition only assuming $W^{1,p_k'}$-regularity for 
$p_k'=\frac {24}7<4$, then one obtains $W^{2,q_k}$-regularity with $q_k=3$.)
Now the right hand sides in (\ref{Sigma eqn}) lie in $W^{1,q_k}(\Om\times\Si)$, 
so lemma~\ref{Sigma reg} gives $W^{1,q_k}$-regularity and -estimates for
$\nabla_\Si B$ on $\Om\times\Si$.
Next, $B\in W^{1,p_k}(\Om,\cA^{0,p}(\Si))$ satisfies the Cauchy-Riemann equation
(\ref{Lagrangian bvp}) with the inhomogeneous term 
$G\in W^{1,q_k}(\Om\times\Si,\rT^*(\Om\times\Si)\otimes\cg)$.
Now we shall use the Sobolev embedding 
$W^{1,q_k}(\Om\times\Si)\hookrightarrow L^{r_k}(\Om\times\Si)$ with
$$
r_k \;=\; \frac{4 q_k}{4 - q_k} 
\;=\; \left\{ \begin{array}{cl}
           \frac{2 p_k}{4 - p_k}  &\text{if}\; p_k<4 ,\\
            \\
           12 &\text{if}\; p_k = 4 .
      \end{array}\right.
$$
Note that $r_k>p_k\geq p$ due to $p_k>2$, so that we have $G\in L^{r_k}(\Om,\cA^{0,p}(\Si))$.
We cannot apply \cite[Theorem~1.2]{W Cauchy} directly because 
that would require the initial regularity $B\in W^{1,2p}(\Om,\cA^{0,p}(\Si))$ for some $p>2$. 
However, we still proceed as in its proof and introduce the coordinates from 
\cite[Lemma~4.3]{W Cauchy} that straighten out the Lagrangian submanifold,
$$
\Th_{s,t} : \cW_{s,t} \to \cA^{0,p}(\Si) .
$$
Here $\cW_{s,t}\subset Y\times Y$ is a neighbourhood of zero, $Y$ is a closed
subspace of $L^p(\Si,\R^m)$ for some $m\in\N$, $\Th$ is in $\cC^{k+1}$-dependence
on $(s,t)$ in a neighbourhood $U\subset\Om$ of some $(s_0,0)\in\Om\cap\pd\H$
and it maps diffeomorphically to a neighbourhood of $B(s,t)$ or $B_0(s,t)$ in case (ii).
Thus one can write
$$
B(s,t) = \Th_{s,t}(v(s,t)) \qquad \forall (s,t)\in U
$$
with $v=(v_1,v_2)\in W^{1,p_k}(U,Y\times Y)$.
Moreover, we have already seen that both $B$ and $\nabla_\Si B$ are $W^{1,q_k}$-regular
on $U\times\Si$, so we have the regularity
$B\in W^{1,q_k}(U,\cA^{1,q_k}(\Si))\subset W^{1,q_k}(U,\cA^{0,s_k}(\Si))$
with corresponding estimates.
Here we have used the Sobolev embedding $W^{1,q_k}(\Si)\hookrightarrow L^{s_k}(\Si)$,
see \cite[Theorem 5.4]{Adams}, for
$$
s_k = \left\{ \begin{array}{cl}
          \frac{2 q_k}{2 - q_k} \;=\; \frac{4 p_k}{8 - 3 p_k} &\text{if}\; p_k<\frac 83 ,\\
           \\
           \frac{44 p_k - 80}{8-p_k} &\text{if}\; p_k \geq \frac 83 , \\
           \\
           \frac {31}2 &\text{if}\; p_k = 4 .
      \end{array}\right.
$$
(Here we have chosen suitable values of $s_k$ for later calculations in case $p_k\geq\frac 83$
and thus $q_k\geq 2$.)
The special structure of the coordinates $\Th$ in \cite[Lemma~4.3]{W Cauchy} 
(it also is a local diffeomorphism between $\cA^{0,s_k}(\Si)$ and a closed subset of 
$L^{s_k}(\Si,\R^{2m})$ since $s_k > p_k>2$) implies that 
$v\in W^{1,q_k}(U,L^{s_k}(\Si,\R^{2m}))$, which will be important later on.

The Cauchy-Riemann equation (\ref{Lagrangian bvp}) now becomes
\[
\left\{\begin{array}{l}
\pd_s v + I \pd_t v = f , \\
v_2(s,0)=0 \quad\forall s\in\R .
\end{array}
\right.
\]
Here $I=(\rd_v\Th)^{-1}*(\rd_v\th)\in W^{1,p_k}(U,{\rm End}(Y\times Y))$ and
$$
f \;=\; (\rd_v\Th)^{-1}(G-\pd_s\Th(v) - * \pd_t\Th(v))
\;\in\; L^{r_k}(U,Y\times Y).
$$
We now approximate $f$ in $L^{r_k}(U,Y\times Y)$ by smooth functions that vanish on $\pd U$, 
then partial integration in \cite[(9)]{W Cauchy} 
yields for all $\p\in\cC^\infty(U,Y^*\times Y^*)$
and a smooth cutoff function as in the proof of \cite[Theorem~1.2]{W Cauchy} 
\begin{align}
\int_U \la h v \,,\, \laplace \p \ra 
&= \int_U \la f \,,\, \pd_s (h\p) - \pd_t (h\cdot I^*\p) \ra 
 \;+ \int_U \la \tilde F \,,\, \p \ra  \nonumber\\
&\quad + \int_{\pd U\cap\pd\H}  \la v_1 \,,\, \pd_t (h\p_1) + \pd_s (h\p_2) \ra .
\label{hv eqn}
\end{align}
Here $\tilde F= (\laplace h) v + 2(\pd_s h) \pd_s v + 2 (\pd_t h) \pd_t v
                + h (\pd_t I) \pd_s v - h (\pd_s I) \pd_t v$ 
contains the crucial terms $(\pd_t I)(\pd_s v)$ and $(\pd_s I)(\pd_t v)$
and thus lies in $L^{\frac 12 p_k}(U,Y\times Y)$.
This is a weak Laplace equation with Dirichlet boundary conditions for $hv_2$, 
Neumann boundary conditions for $hv_1$, and with the inhomogeneous term in
$W^{-1,r_k}(U,Y\times Y)$. The latter is the dual space of 
$W^{1,r_k'}(U,Y^*\times Y^*)$ with $\frac 1{r_k}+\frac 1{r_k'}=1$.
(The inclusion $L^{\frac 12 p_k}(U)\hookrightarrow W^{-1,r_k}(U)$ is 
continuous as can be seen via the dual embedding that is due to 
\hbox{$\frac 12 - \frac 1{r_k'} \geq - 1 + \frac 1{p_k/2}$}.)
Recall that $Y\subset L^p(\Si,\R^m)$ is a closed subspace.
Since $r_k > p$ the special regularity result \cite[Lemma~2.1]{W Cauchy} for 
the Laplace equation with values in a Banach space cannot be applied to deduce
$hv\in W^{1,r_k}(U,Y\times Y)$.
However, the general regularity theory for the Laplace equation extends to functions
with values in a Hilbert space (c.f.\ \cite{W}).
So we use the embedding $L^p(\S)\hookrightarrow L^2(\S)$. Then (\ref{hv eqn}) is a weak Laplace
equation with the inhomogeneous term in $W^{-1,r_k}(U,L^2(\S,\R^{2m}))$ and enables us to deduce
$hv\in W^{1,r_k}(U,L^2(\Si,\R^{2m}))$ and thus $v\in W^{1,r_k}(\tU,L^2(\Si,\R^{2m}))$ with the 
corresponding estimates for some smaller domain $\tU$ (a finite union of these still covers a 
neighbourhood of $\Om\cap\pd\H$).
Furthermore, recall that $v\in W^{1,q_k}(U,L^{s_k}(\Si,\R^{2m}))$.
Now we claim that the following inclusion with the corresponding estimates holds for 
some suitable $p_{k+1}$
\begin{equation}\label{v pk+1}
W^{1,r_k}(\tU,L^2(\Si)) \;\cap\; W^{1,q_k}(\tU,L^{s_k}(\Si))
\;\subset\; W^{1,p_{k+1}}(\tU,L^{p_{k+1}}(\Si)) .
\end{equation}
To show (\ref{v pk+1}) it suffices to estimate the $L^{p_{k+1}}(\tU\times\Si)$-norm
of a smooth function by its $L^{r_k}(\tU,L^2(\Si))$- and $L^{q_k}(\tU,L^{s_k}(\Si))$-norms.
Let $\a>2$ and $t\in[1,2)$, then the H\"older inequality gives for all
$f\in\cC^\infty(\tU\times\Si,\R^{2m})$
\begin{align*}
\| f \|_{L^\a(\tU\times\Si)}^\a
&= \int_\tU\int_\Si |f|^t |f|^{\a-t} \\
&\leq \int_\tU \|f\|_{L^2(\Si)}^t \|f\|_{L^{2\frac{\a-t}{2-t}}(\Si)}^{\a-t} \\
&\leq \|f\|_{L^r(\tU,L^2(\Si))}^t 
      \|f\|_{L^{r\frac{\a-t}{r-t}}(\tU,L^{2\frac{\a-t}{2-t}}(\Si))}^{\a-t} \\ 
&\leq \|f\|_{L^r(\tU,L^2(\Si))}^\a
    + \|f\|_{L^{r\frac{\a-t}{r-t}}(\tU,L^{2\frac{\a-t}{2-t}}(\Si))}^\a .
\end{align*}
Here we abbreviated $r:=r_k>p_k>2$. Now we want
\begin{equation} \label{q,s}
q_k=\frac{r_k(\a-t)}{r_k-t} \qquad\text{and}\qquad
s_k=\frac{2(\a-t)}{2-t} .
\end{equation}
Indeed, in the case $p_k=4$ our choices $q_k=3$, $r_k=12$, and $s_k=\frac{31}2$ together with 
$t:=\frac 53$ and $\a:=\frac{17}{4}$ solve these equations. 
So we obtain $p_{k+1}=\a=\frac{17}{16} p_k$.
In case $p_k<4$ the first equation gives
\begin{equation}\label{alpha}
\a = \frac{4+t}{8-p_k} p_k .
\end{equation}
If moreover $p_k\geq\frac 83$, then we choose $t:=\frac 53$ to obtain
$\a \;=\; \frac{17}{24-3p_k} p_k \;\geq\; \frac{17}{16} p_k $.
This also solves (\ref{q,s}) with our choice $s_k=\frac{44 p_k - 80}{8-p_k}$,
so we obtain $p_{k+1}=\frac{17}{16} p_k$.
Finally, in case $\frac 83 > p_k >2$ one obtains from (\ref{q,s})
$$
t\;=\;\frac{p_k^2}{-p_k^2 + 7 p_k -8} \;\in\; [1,2).
$$
Inserting this in (\ref{alpha}) yields $\a=\th(p_k)\cdot p_k$ with
$$
\th(p_k) = \frac{3p_k-4}{-p_k^2+7p_k-8} .
$$
One then checks that $\th(2)=1$ and $\th'(p)>0$ for $p>2$, thus $\th(p)>1$ for $p>2$.
Moreover, $\th(\frac 83)=\frac 98$, so $\th(p')=\frac{17}{16}$ for some 
$p'\in (2,\frac 83)$. Now for $p\geq p'$ we extend the function constantly
to obtain a monotonely increasing function $\th:(2,4]\to(1,\frac{17}{16}]$.
With this modified function we finally choose $p_{k+1}=\th(p_k)\cdot p_k$ for all
$2<p_k\leq 4$.
This finishes the proof of (\ref{v pk+1}) and thus shows that
$v\in W^{1,p_{k+1}}(\tU,L^{p_{k+1}}(\Si))$.

In addition, note that our choice of $p_{k+1}\leq\a$ will always satisfy 
$p_{k+1}\leq r_k$. In case $p_k=4$ see the actual numbers, in case $p_k<4$ this
is due to (\ref{alpha}), $t\leq 2$, and $p_k>2$, 
$$
\a\;\leq\; \frac 6{8-p_k} p_k \;\leq\; \frac 2{4-p_k} p_k \;=\; r_k .
$$
Now we again use the special structure of the coordinates $\Th$ in \cite[Lemma~4.3]{W Cauchy} 
to deduce that $B=\Th\comp v\in W^{1,p_{k+1}}(\tU,\cA^{0,p_{k+1}}(\Si))$ with the corresponding
estimates.
Above, we already established the $W^{1,r_k}$- and thus $W^{1,p_{k+1}}$-regularity and 
-estimates for $\P$ and $\Psi$ as well as $B\in L^{p_{k+1}}(\tU,\cA^{1,p_{k+1}}(\Si))$.
(Recall the Sobolev embedding $W^{1,q_k}\hookrightarrow L^{r_k}$, and that $p_k\geq q_k$
and $r_k\geq p_{k+1}$, so we have $L^{r_k}(\tU,L^{r_k}(\Si))$-regularity of $B$ as well as 
$\nabla_\Si B$.)
Putting all this together we have established the $W^{1,p_{k+1}}$-regularity and
-estimates for $\bar\t^*A$ over $\tU_i\times\Si$, where the $\tU_i$ cover a neighbourhood
of $\Om_{k+1}\cap\pd\H$. The interior regularity again follows directly from
proposition~\ref{component reg}. 

This iteration gives a sequence $(p_k)$ with $p_{k+1}=\th(p_k)\cdot p_k \geq \th(p)\cdot p_k$.
So this sequence grows at a rate greater or equal to $\th(p)>\th(2)=1$ and hence reaches 
$p_N>4$ after finitely many steps. This finishes the proof of II) and the theorem.
\QED \\

\noindent
{\bf Proof of theorem~A : }\\
Fix a solution $A\in\cA^{1,p}_{\rm loc}(X)$ of (\ref{ASD bvp}) with $p>2$.
We have to find a gauge transformation $u\in\cG^{2,p}_{\rm loc}(X)$ such that
$u^*A\in\cA(X)$ is smooth. Recall that the manifold $X=\bigcup_{k\in\N} X_k$ is exhausted
by compact submanifolds $X_k$ meeting the assumptions of proposition~\ref{DKargument}.
So it suffices to prove for every $k\in\N$ that there exists a gauge transformation 
$u\in\cG^{2,p}(X_k)$ such that $u^*A|_{X_k}$ is smooth.

For that purpose fix $k\in\N$ and choose a compact submanifold $M \subset X$ that is large 
enough such that theorem~\ref{reg} applies to the compact subset $K:=X_k\subset M$.
Next, choose $A_0\in\cA(M)$ such that $\|A-A_0\|_{W^{1,p}(M)}$ and \hbox{$\|A-A_0\|_{L^q(M)}$} are 
sufficiently small for the local slice theorem, proposition~\ref{local slice thm}, to apply 
to $A_0$ with the reference connection $\hat A=A$. 
Here due to $p>2$ one can choose $q>4$ in the local slice theorem such that the
Sobolev embedding $W^{1,p}(M)\hookrightarrow L^q(M)$ holds.
Then by proposition~\ref{local slice thm} and remark~\ref{local slice rmk}~(i) one 
obtains a gauge transformation $u\in\cG^{2,p}(M)$ such that $u^*A$ is in relative Coulomb 
gauge with respect to $A_0$. Moreover, $u^*A$ also solves (\ref{ASD bvp}) since both the 
anti-self-duality equation and the Lagrangian submanifolds $\cL_i$ are gauge 
invariant.
The latter is due to the fact that $u$ restricts to a gauge transformation in $\cG^{1,p}(\Si_i)$ 
on each boundary slice $\t_i(\{s\}\times\Si_i)$ due to the Sobolev embedding 
$\cG^{2,p}(\cU_i\times\S)\subset W^{1,p}(\cU_i,\cG^{1,p}(\S_i))\hookrightarrow
\cC^0(\cU_i,\cG^{1,p}(\S_i))$.
So $u^*A\in\cA^{1,p}(M)$ is a solution of (\ref{bvp}) and theorem~\ref{reg}~(i) 
asserts that $u^*A|_{X_k}\in\cA(X_k)$ is indeed smooth.

Such a gauge transformation $u\in\cG^{2,p}(X_k)$ can be found for every $k\in\N$, hence
proposition~\ref{DKargument}~(i) asserts that there exists a gauge transformation
$u\in\cG^{2,p}_{\rm loc}(X)$ on the full noncompact manifold such that
$u^*A\in\cA(X)$ is smooth as claimed.
\QED \\

\noindent
{\bf Proof of theorem~B : } \\
Fix a smoothly convergent sequence of metrics $g^\n\to g$ that are compatible to
$\t$ and let $A^\n\in\cA^{1,p}_{\rm loc}(X)$ be a sequence of solutions of 
(\ref{ASD bvp}) with respect to the metrics $g^\n$.
Recall that the manifold $X=\bigcup_{k\in\N} X_k$ is exhausted
by compact submanifolds $X_k$ meeting the assumptions of proposition~\ref{DKargument}.
We will find a subsequence (again denoted $A^\n$) and a sequence of gauge
transformations $u^\n\in\cG^{2,p}_{\rm loc}(X)$ such that the sequence $u^{\n\;*}A^\n$ 
is bounded in the $W^{\ell,p}$-norm on $X_k$ for all $\ell\in\N$ and $k\in\N$. 
Then due to the compact Sobolev embeddings 
$W^{\ell,p}(X_k)\hookrightarrow\cC^{\ell-2}(X_k)$ one finds a further (diagonal) 
subsequence that converges uniformly with all derivatives on every compact 
subset of $X$.

By proposition~\ref{DKargument}~(ii) it suffices to construct the gauge
transformations and establish the claimed uniform bounds over $X_k$ for all 
$k\in\N$ and for any subsequence of the connections (again denoted $A^\n$).
So fix $k\in\N$ and choose a compact submanifold $M \subset X$ such that 
theorem~\ref{reg} holds with $K=X_k\subset M$.
Choose a further compact submanifold $M' \subset X$ such that 
theorem~\ref{reg} holds with $K=M\subset M'$.
Then by assumption of the theorem $\|F_{A^\n}\|_{L^p(M')}$ is uniformly bounded.
So the weak Uhlenbeck compactness theorem, proposition~\ref{weak comp}, provides a 
subsequence (still denoted $A^\n$), a limit connection $A_0\in\cA^{1,p}(M')$, 
and gauge transformations $u^\n\in\cG^{2,p}(M')$ such that 
$u^{\n\;*}A^\n \to A_0$ in the weak $W^{1,p}$-topology.
The limit $A_0$ then satisfies the boundary value problem (\ref{ASD bvp}) with respect to the 
limit metric $g$. So as in the proof of theorem~A one finds a gauge 
transformation $u_0\in\cG^{2,p}(M)$ such that $u_0^*A_0\in\cA(M)$ is smooth.
Now replace $A_0$ by $u_0^*A_0$ and $u^\n$ by $u^\n u_0\in\cG^{2,p}(M)$, then one still has a 
$W^{1,p}$-bound, $\|u^{\n\;*}A^\n - A_0\|_{W^{1,p}(M)}\leq c_0$ for some constant $c_0$,
see \cite[Lemma~A.5]{W Cauchy}.

Due to $p>2$ one can now choose $q>4$ in the local slice theorem such that the
Sobolev embedding $W^{1,p}(M)\hookrightarrow L^q(M)$ is compact.
Hence for a further subsequence of the connections $u^{\n\;*}A^\n \to A_0$ in the $L^q$-norm.
Let $\ep>0$ be the constant from proposition~\ref{local slice thm} for the reference 
connection $\hat A=A_0$, then one finds a further subsequence such that
$\|u^{\n\;*}A^\n - A_0\|_{L^q(M)}\leq\ep$ for all $\n\in\N$. 
So the local slice theorem provides further gauge transformations 
$\tilde u^\n\in\cG^{2,p}(M)$ such that the $\tilde u^{\n\;*}A^\n$ are in 
relative Coulomb gauge with respect to $A_0$. 
The gauge transformed connections still solve (\ref{ASD bvp}), hence 
the $\tilde u^{\n\;*}A^\n$ are solutions of (\ref{bvp}).
Moreover, we have 
$\|\tilde u^{\n\;*}A^\n - A_0\|_q \leq C_{\scriptscriptstyle CG} \| u^{\n\;*}A^\n - A_0\|_q$,
hence $\tilde u^{\n\;*}A^\n\to A_0$ in the $L^q$-norm, and 
$$
\bigl\| \tilde u^{\n\;*}A^\n - A_0 \bigr\|_{W^{1,p}(M)} \leq 
C_{\scriptscriptstyle CG}c_0 .
$$
The higher $W^{k,p}$-bounds will now follow from theorem~\ref{reg}, so we first have to 
verify its assumptions.
Fix the metric $g_0:=g$ and a compact neighbourhood 
$\cV=\bigcup_{i=1}^n \bar\t_{0,i}(\cU_i\times\Si_i)$ of $K\cap\pd X$. 
Then the $\bar\t_{0,i}^*(\tilde u^{\n\;*}A^\n-A_0)|_{\Si_i}$ are uniformly $W^{1,p}$-bounded 
and converge to zero in the $L^q$-norm on $\cU_i\times\cS_i$ as seen above. Now the embedding
$$
W^{1,p}(\cU_i\times\Si_i,\rT^*\Si_i\otimes\cg)
\hookrightarrow L^\infty(\cU_i,\cA^{0,p}(\Si_i))
$$
is compact by lemma~\ref{Sob emb}.
Thus one finds a subsequence such that the $\bar\t_{0,i}^*(\tilde u^{\n\;*}A^\n)|_{\Si_i}$
converge in $L^\infty(\cU_i,\cA^{0,p}(\Si_i))$. The limit can only be 
$\bar\t_{0,i}^*A_0|_{\Si_i}$ since this already is the $L^q$-limit.
Now in theorem~\ref{reg}~(ii) choose the constant $C_1=C_{\scriptscriptstyle CG}c_0$ 
and let $\d>0$ be the constant determined from $C_1$. Then one can take a subsequence such that
$$
\|\bar\t_{0,i}^*(\tilde u^{\n\;*}A^\n - A_0)_{\Si_i}\|_{L^\infty(\cU_i,\cA^{0,p}(\Si_i))}\leq\d 
\qquad\forall i=1,\ldots,n, \forall \n .
$$
Now theorem~\ref{reg}~(ii) provides the claimed uniform bounds as follows.
Fix $\ell\in\N$, then $\|g^\n-g\|_{\cC^{\ell+2}(M)}\leq\d_\ell$
for all $\n\geq\n_\ell$ with some $\n_\ell\in\N$, and thus
$$
\bigl\| \tilde u^{\n\;*}A^\n - A_0 \bigr\|_{W^{\ell,p}(X_k)} \leq C_\ell 
\qquad\forall\n\geq\n_\ell.
$$
This finally implies the uniform bound for this subsequence,
$$
\sup_{\n\in\N} \bigl\| \tilde u^{\n\;*}A^\n \bigr\|_{W^{\ell,p}(X_k)} < \infty.
$$
Here the gauge transformations $\tilde u^\n\in\cG^{2,p}(X_k)$ still depend on $k\in\N$ and
are only defined on $X_k$. But now proposition~\ref{DKargument}~(ii) provides a subsequence
of $(A^\n)$ and gauge transformations $u^\n\in\cG^{2,p}_{\rm loc}(X)$ defined on the full
noncompact manifold such that $u^{\n\;*}A^\n$ is uniformly bounded in every $W^{\ell,p}$-norm 
on every compact submanifold $X_k$.
Now one can iteratively use the compact Sobolev embeddings
$W^{\ell+2,p}(X_\ell)\hookrightarrow\cC^{\ell}(X_\ell)$ for each $\ell\in\N$
to find a further subsequence of the connections that converges in $\cC^\ell(X_\ell)$.
If in each step one fixes one further element of the sequence, then this iteration finally
yields a sequence of connections that converges uniformly with all derivatives on every compact 
subset of $X$ to a smooth connection $A\in\cA(X)$.
\QED

\section{Fredholm theory}
\label{Fredholm}
This section concerns the linearization of the boundary value problem
(\ref{ASD bvp intro}) in the special case of a compact $4$-manifold of the form $X=S^1\times Y$, 
where $Y$ is a compact orientable $3$-manifold whose boundary $\pd Y=\S$ is a disjoint union 
of connected Riemann surfaces.
The aim of this section is to prove theorem~C. 

So we equip $S^1\times Y$ with a product metric $\tg=\ds^2+g_s$ (where $g_s$ is an 
$S^1$-family of metrics on $Y$) and assume that this is compatible with the
natural space-time splitting of the boundary $\pd X = S^1\times\Si$.
This means that for some $\D>0$ there exists an embedding 
$$
\t: S^1\times[0,\D)\times\Si \hookrightarrow  S^1\times Y 
$$ 
preserving the boundary, $\t(s,0,z)=(s,z)$ for all $s\in S^1$ and $z\in\Si$, 
such that
$$
\t^*\tg = \ds^2 + \dt^2 + g_{s,t} .
$$ 
Here $g_{s,t}$ is a smooth family of metrics on $\Si$.
This assumption on the metric implies that the normal geodesics are independent 
of $s\in S^1$ in a neighbourhood of the boundary. So in fact, the embedding is 
given by $\t(s,t,z)=(s,\g_z(t))$, where $\g$ is the normal geodesic starting at 
$z\in\Si$.
This seems like a very restrictive assumption, but it suffices for our application
to Riemannian $4$-manifolds with a boundary space-time splitting.
Indeed, the neighbourhoods of the compact boundary components are isometric to
$S^1\times Y$ with $Y=[0,\D]\times\Si$ and a metric $\ds^2 + \dt^2 + g_{s,t}$.

Now fix $p>2$ and let $\cL\subset\cA^{0,p}_{\rm flat}(\Si)$ be a gauge invariant 
Lagrangian submanifold of $\cA^{0,p}(\Si)$ as in the introduction. Then for
$\tA\in\cA^{1,p}(S^1\times Y)$ we consider the nonlinear boundary value problem
\begin{equation} \label{ASD bvp F}
\left\{\begin{array}{l}
*F_\tA + F_\tA =0 ,\\
\tA|_{\{s\}\times\pd Y} \in\cL \quad\forall s\in S^1 .
\end{array}\right.
\end{equation}
Fix a smooth connection $\tA\in\cA(S^1\times Y)$ with Lagrangian boundary values
(but not necessarily a solution of this boundary value problem).
It can be decomposed as $\tA=A+\P\ds$ with $\P\in\cC^\infty(S^1\times Y,\cg)$ and with
\hbox{$A\in\cC^\infty(S^1\times Y,\rT^*Y\otimes\cg)$} satisfying
$A_s:=A(s)|_{\pd Y}\in\cL$ for all $s\in S^1$.
Similarly, a tangent vector $\ta$ to $\cA^{1,p}(S^1\times Y)$ decomposes as
$\ta=\a+\ph\ds$ with $\ph\in W^{1,p}(S^1\times Y,\cg)$ and
$\a\in W^{1,p}(S^1\times Y,\rT^*Y\otimes\cg)$.
Now let $E_A^{1,p}\subset W^{1,p}(S^1\times Y,\rT^*Y\otimes\cg)$ be the subspace
of $S^1$-families of $1$-forms $\a$ that satisfy the boundary conditions from
the linearization of (\ref{ASD bvp F}) and the Coulomb gauge,
$$
*\a(s)|_{\pd Y} = 0 \qquad\text{and}\qquad
 \a(s)|_{\pd Y}\in\rT_{A_s}\cL \qquad\text{for all}\: s\in S^1 .
$$
Then the linearized operator for the study of the moduli space of gauge equivalence classes
of solutions of (\ref{ASD bvp F}) is as in the introduction
$$
D_{(A,\P)} :
E_A^{1,p}\times W^{1,p}(S^1\times Y,\cg)  \longrightarrow 
L^p(S^1\times Y,\rT^*Y\otimes\cg) \times L^p(S^1\times Y,\cg) 
$$
given by
$$
D_{(A,\P)}(\a,\ph) \;=\; \bigl( \nabla_s\a  - \rd_A\ph + *\rd_A\a \,,\, 
                                \nabla_s\ph - \rd_A^*\a \bigr) .
$$
Here $\rd_A$ denotes the exterior derivative corresponding to the connection $A(s)$ 
on $Y$ for all $s\in S^1$, $*$ denotes the Hodge operator on $Y$ with respect to the
$s$-dependent metric $g_s$ on $Y$, and we use the notation 
$\nabla_s\a := \pd_s\a + [\P,\a]$.
Our main result, theorem~C~(i), is the Fredholm property of $D_{(A,\P)}$.
We now give an outline of its proof.

The first crucial point is the estimate, theorem~C~(ii), which ensures that $D_{(A,\P)}$ 
has a closed image and a finite dimensional kernel.
It can be rephrased as follows due to the identities
\begin{align*}
\rd_{\tA}^+\tilde\a 
&=   \half *\bigl( \nabla_s\a - \rd_A\ph + *\rd_A\a \bigr) 
   - \half \bigl( \nabla_s\a - \rd_A\ph + *\rd_A\a \bigr)\wedge\ds ,\\
\rd_\tA^*\ta &= - \nabla_s\ph + \rd_A^*\a .
\end{align*}

\begin{lem} \label{estimate}
There is a constant $C$ such that for all $\ta\in W^{1,p}(X,\rT^*X\otimes\cg)$ 
satisfying
$$
*\ta|_{\pd X} = 0 
\qquad\text{and}\qquad
\ta|_{\{s\}\times\pd Y}\in\rT_{A_s}\cL \quad\forall s\in S^1
$$
one has the estimate
$$
\|\ta\|_{W^{1,p}} \leq C \bigl( \| \rd_\tA^+ \ta \|_p 
                              + \| \rd_\tA^* \ta \|_p + \| \ta \|_p  \bigr).
$$
\end{lem}

The second part of the Fredholm theory for $D_{(A,\P)}$ is the identification of the 
cokernel with the kernel of a slightly modified linearized operator,
which will be used to prove that the cokernel is finite dimensional. 
To be more precise let $\s:S^1\times Y \to S^1\times Y$ denote the reflection given by
$\s(s,y):=(-s,y)$, where $S^1\cong\R/\Z$. Then we will establish the following duality:
$$
(\b,\z) \in  (\im D_{(A,\P)})^\perp
\;\Longleftrightarrow\;
(\b\comp\s,\z\comp\s) \in \ker D_{\s^*(A,\P)},
$$
where $D_{\s^*(A,\P)}$ is the linearized operator at the connection
$\s^*\tA=A\comp\s-\P\comp\s$ with respect to the metric $\s^*\tg$ on 
$S^1\times Y$.
Once we know that $\im D_{(A,\P)}$ is closed, this gives an isomorphism between
$({\rm coker}\,D_{(A,\P)})^* \cong (\im D_{(A,\P)})^\perp$ and 
$\ker D_{\s^*(A,\P)}$. 
Here $Z^*$ denotes the dual space of a Banach space $Z$, and for a subspace $Y\subset Z$
we denote by $Y^\perp\subset Z^*$ the space of linear functionals that vanish on $Y$.
Now the estimate in theorem~C~(ii) will also apply to $D_{\s^*(A,\P)}$,
and this implies that its kernel -- and hence the cokernel of $D_{(A,\P)}$ --
is of finite dimension.
The main difficulty in establishing the above duality is the regularity result theorem~C~(iii).

This regularity as well as the estimate in theorem~C~(ii) or lemma~\ref{estimate} will be proven 
analogously to the nonlinear regularity and estimates in section~\ref{reg,comp}. 
Again, the interior regularity and estimate is standard elliptic theory, and one has to use a 
splitting near the boundary. We shall show that the $S^1$- and the normal component both 
satisfy a Laplace equation with Neumann and Dirichlet boundary conditions 
respectively. 
The $\Si$-component will again gives rise to a (weak) Cauchy-Riemann equation in a 
Banach space, only this time the boundary values will lie in the tangent space of the 
Lagrangian.
In contrast to the required $L^p$-estimates we shall first show that the $L^2$-estimate for 
$L^p$-regular $1$-forms can be obtained by more elementary methods.
These were already outlined in \cite{Sa1} as a first indication for the Fredholm property
of the boundary value problem (\ref{ASD bvp F}).

Let $\ta\in W^{1,p}(X,\rT^*X\otimes\cg)$ be as supposed for some $p>2$.
From the first boundary condition $*\ta|_{\pd X}=0$ one obtains
$$
\|\nabla \ta\|_2^2
= \|\rd \ta\|_2^2 + \|\rd^* \ta\|_2^2 
 - \int_{\pd X} \tilde g(Y_\ta,\nabla_{Y_\ta}\n) .
$$
Here one has 
$\int_{\pd X} \tilde g(Y_\ta,\nabla_{Y_\ta}\n) \geq - C \|\ta\|_{L^2(\pd X)}^2$
since the vector field $Y_\ta$ is given by $\i_{Y_\ta}\tilde g = \ta$.
In this last term one uses the following estimate for general $1<p<\infty$.

Let $\t:[0,\D)\times\pd X \to X$ be a diffeomorphism to a tubular neighbourhood
of $\pd X$ in $X$. Then for all $\d>0$ one finds a constant $C_\d$ such that for 
all $f\in W^{1,p}(X)$ 
\begin{align}
&\|f\|_{L^p(\pd X)}^p \nonumber\\
&= \int_{\pd X}\int_0^1 \frac\rd\ds \Bigl( (s-1)|f(\t(s,z))|^p \Bigl)
 \,\ds\,\rd^3z \nonumber\\
&\leq \int_{\pd X}\int_0^1 |f(\t(s,z))|^p \,\ds\,\rd^3z 
    + \int_{\pd X}\int_0^1 p|f(\t(s,z))|^{p-1}|\pd_s f(\t(s,z))|
  \,\ds\,\rd^3z \nonumber\\
&\leq C \bigl( \| f \|_{L^p(X)}^p 
             + \| f \|_{L^p(X)}^{p-1} \|\nabla f \|_{L^p(X)} \bigr) \nonumber\\
&\leq \bigl( \d \| f \|_{W^{1,p}(X)} + C_\d \| f \|_{L^p(X)} \bigr)^p .
\label{bdy est}
\end{align}
This uses the fact that for all $x,y\geq 0$ and $\d>0$
$$
x^{p-1} y \;\leq\;
\left\{\begin{array}{ll}
\d^p y^p &; x\leq \d^{\frac p{p-1}} y \\
\d^{-\frac p{p-1}} x^p &; x\geq \d^{\frac p{p-1}} y 
\end{array}\right\}
\;\leq\; \bigl( \d y + \d^{-\frac 1{p-1}} x \bigr) ^p .
$$
So we obtain
\begin{equation} \label{L2 est0}
\|\ta\|_{W^{1,2}}
\leq C \bigl( \|\rd_\tA\ta\|_2 + \|\rd_\tA^*\ta\|_2  + \|\ta\|_2 \bigr).
\end{equation}
In fact, the analogous $W^{1,p}$-estimates hold true for general $p$, as is proven e.g.\ 
in \cite[Theorem 6.1]{W}. However, in the case $p=2$ one can calculate further for all $\d>0$
\begin{align}
\|\rd_\tA \ta\|_2^2 
&= \int_X \la \rd_\tA\ta \,,\, 2\rd_\tA^+\ta \ra  
 - \int_X \la \rd_\tA\ta \wedge \rd_\tA\ta \ra  \nonumber\\
&= 2 \|\rd_\tA^+\ta\|_2^2 - \int_X \la \ta \wedge [F_\tA\wedge \ta]\ra
  - \int_{\pd X} \la \ta \wedge \rd_\tA \ta\ra  \nonumber\\
&\leq 2 \|\rd_\tA^+\ta\|_2^2 + C_\d \|\ta\|_2^2 
      + \d \|\ta\|_{W^{1,2}}^2 .  \label{L2 est}
\end{align}
Here the boundary term above is estimated as follows. 
We use the universal covering of $S^1=\R/\Z$ to integrate over $[0,1]\times \pd Y$ 
instead of \hbox{$\pd X=S^1\times\pd Y$}.
Introduce $A:=(A_s)_{s\in S^1}$, which is a smooth path in $\cL$. Then
using the splitting $\ta|_{\pd X}=\a+\ph\ds$ with 
$\a:S^1\times\Si \to \rT^*\Si\otimes\cg$ and $\ph:S^1\times\Si \to \cg$
one obtains
\begin{align*}
&- \int_{\pd X} \la \ta \wedge \rd_\tA\ta\ra \\
&= - \int_0^1 \int_\Si  \la \ph \,,\,  \rd_A\a \ra \dvol_\Si\wedge\ds \;
   - \int_0^1 \int_\Si  \la \a \wedge (\rd_A\ph - \nabla_s\a )\ra \wedge\ds  \\
&= \int_0^1 \int_\Si  \la \a \wedge \nabla_s\a\ra \wedge\ds  \\
&\leq \d \|\ta\|_{W^{1,2}(X)}^2  + C'_\d \|\ta\|_{L^2(X)}^2 .
\end{align*}
Firstly, we have used the fact that $\rd_A\a|_\Si=0$ since 
$\a(s)\in\rT_{A_s}\cL\subset\ker\rd_{A_s}$ for all $s\in S^1$.
Secondly, we have also used that both $\a$ and $\rd_A\ph$ lie in $\rT_A\cL$, hence the 
symplectic form $\int_\Si \la \a\wedge \rd_A\ph \ra$ vanishes for all $s\in S^1$.
This is not strictly true since $\ta$ only restricts to an 
$L^p$-regular $1$-form on $\pd X$. However, as $1$-form on $[0,1]\times Y$
it can be approximated as follows by smooth $1$-forms that meet the Lagrangian 
boundary condition on $[0,1]\times\Si$.

We use the linearization of the coordinates in \cite[Lemma~4.3]{W Cauchy}
at $A_s$ for all $s\in[0,1]$. Since the path $s\mapsto A_s\in\cL\cap\cA(\Si)$ is smooth, 
this gives a smooth path of diffeomorphisms $\Th_s$ for any $q>2$,
$$
\Th_s:
\begin{array}{ccc}
Z\times Z &\longrightarrow& L^q(\Si,\rT^*\Si\otimes\cg) \\
(\x,v,\z,w) &\longmapsto& \rd_{A_s}\x + \tsum_{i=1}^m v^i  \g_i(s) 
                   + *\rd_{A_s}\z + \tsum_{i=1}^m w^i *\g_i(s) ,
\end{array}
$$
where $Z:=W^{1,q}_z(\Si,\cg) \times \R^m$ and the
$\g_i\in\cC^\infty([0,1]\times\Si,\rT^*\Si\otimes\cg)$ satisfy
$\g_i(s)\in\rT_{A_s}\cL$ for all $s\in [0,1]$.
We perform above estimate on $[0,1]\times Y$ since we can not necessarily achieve 
$\g_i(0)=\g_i(1)$.
In these coordinates, we mollify to obtain the required smooth approximations of $\ta$ near 
the boundary.
Furthermore, we use these coordinates for $q=3$ to write the smooth approximations
on the boundary as $\a(s) = \rd_{A_s}\x(s) + \tsum_{i=1}^m v^i(s) \g_i(s) $ with 
$\|\x(s)\|_{W^{1,3}(\Si)} + |v(s)| \leq C \|\a(s)\|_{L^3(\Si)}$.
Then for all $s\in[0,1]$
\begin{align*}
\int_\Si \la \a(s) \wedge \nabla_s\a(s) \ra 
&= \int_\Si \la \a\wedge\bigl( \rd_{A_s}\pd_s \x + \tsum_{i=1}^m \pd_s v^i \cdot \g_i \bigr)\ra\\
&\quad + \int_\Si \la \a \wedge \bigl( [\P,\a] + [\pd_s A,\x] 
                              + \tsum_{i=1}^m v^i \cdot \pd_s\g_i  \bigr) \ra \\
&\leq C \|\a(s)\|_{L^2(\Si)} \|\a(s)\|_{L^3(\Si)} .
\end{align*}
Here the crucial point is that $\rd_A\pd_s\x$ and $\pd_s v^i \cdot\g_i$ are tangent to the 
Lagrangian, hence the first term vanishes.
Now one uses (\ref{bdy est}) for $p=2$ and the Sobolev trace theorem
(the restriction $W^{1,2}(X)\to L^3(\pd X)$ is continuous by e.g.\ \cite[Theorem~6.2]{Adams} ) 
to obtain the estimate,
\begin{align*}
\int_0^1 \int_\Si  \la \a \wedge \nabla_s\a\ra \wedge\ds  
&\leq C \|\ta\|_{L^2(\pd X)} \|\ta\|_{L^3(\pd X)}  \\
&\leq \tfrac\d 2 \|\ta\|_{W^{1,2}(X)}^2 
     + C_\d \|\ta\|_{L^2(X)} \|\ta\|_{W^{1,2}(X)} \\
&\leq \d \|\ta\|_{W^{1,2}(X)}^2 
     + C'_\d \|\ta\|_{L^2(X)}^2 .
\end{align*}
This proves (\ref{L2 est}). Now $\d>0$ can be chosen arbitrarily small, so
the term $\|\ta\|_{W^{1,2}}$ can be absorbed into the left hand side of (\ref{L2 est0}),
and thus one obtains the claimed estimate
$$
\|\ta\|_{W^{1,2}}
\;\leq\; C \bigl( \|\rd_\tA^+\ta\|_2 + \|\rd_\tA^*\ta\|_2 + \|\ta\|_2 \bigr).
$$

\noindent
{\bf Proof of theorem~C~(ii) or lemma \ref{estimate} : } \\
We will use lemma~\ref{Hodge prop} for the manifold $M:=S^1\times Y$ in several different
cases to obtain the estimate for different components of $\ta$. 
The first weak equation in lemma~\ref{Hodge prop} is the same in all cases.
For all $\e\in\cC^\infty(M;\cg)$
\begin{align*}
\int_M \la \ta \,,\, \rd\e \ra
&= \int_M \la \rd^*\ta \,,\, \e \ra
 + \int_{\pd M} \la \e \,,\, *\ta \ra     \\
&= \int_M \la \rd_\tA^*\ta + *[\tA\wedge *\ta] \,,\, \e \ra  
\quad= \int_M \la f \,,\, \e \ra .
\end{align*}
Here one uses the fact that $*\ta|_{\pd M}=0$ .
Then $f\in L^p(M,\cg)$ and
\begin{equation} \label{est f}
\|f\|_p \leq \|\rd_\tA^*\ta\|_p + 2 \|\tA\|_\infty \|\ta\|_p .
\end{equation}
For the second weak equation lemma~\ref{Hodge prop} we obtain for all $\l\in\Om^1(M;\cg)$
\begin{align}
\int_M \la \ta \,,\, \rd^*\rd \l \ra
&= \int_M \la \ta \,,\, \rd^*\rd \l + \rd^* * \rd\l \ra \label{2nd eqn}\\
&= \int_M \la \g \,,\, \rd \l \ra  
 - \int_{S^1\times\pd Y} \la \ta \wedge * \rd \l \ra 
 - \int_{S^1\times\pd Y} \la \ta \wedge \rd \l \ra ,   \nonumber
\end{align}
where $\g\,=\,\rd\ta + * \rd\ta \,=\, 2\rd_\tA^+\ta -2[\tA\wedge\ta]^+ 
       \,\in L^p(M,\L^2\rT^*M\otimes\cg)$ 
and
\begin{equation} \label{est g}
\|\g\|_p \leq 2 \|\rd_\tA^+\ta\|_p 
                   + 4 \|\tA\|_\infty \|\ta\|_p .
\end{equation}
Now recall that there is an embedding 
$\t:S^1\times[0,\D)\times\Si \hookrightarrow S^1\times Y$
to a tubular neighbourhood of $S^1\times\pd Y$ such that 
$\t^*\tg = \ds^2 +\dt^2 + g_{s,t}$ for a family $g_{s,t}$ of metrics on $\Si$.
One can then cover $M=S^1\times Y$ with $\t(S^1\times[0,\frac\D 2]\times\Si)$ 
and a compact subset $V\subset M\setminus\pd M$.

For the claimed estimate of $\ta$ over $V$ it suffices to use lemma~\ref{Hodge prop}
for vector fields $X\in\G(\rT M)$ that equal to coordinate vector fields on $V$ and
vanish on $\pd M$.
So one has to consider (\ref{2nd eqn}) for $\l=\p\cdot \i_X\tilde g$ with
$\p\in\cC^\infty_\d(M,\cg)$. Then both boundary terms vanish and hence lemma~\ref{Hodge prop} 
directly asserts, with some constants $C$ and $C_V$, that
\begin{align*}
\|\ta\|_{W^{1,p}(V)} 
&\leq C \bigl( \|f\|_{L^p(M)} + \|\g\|_{L^p(M)} + \|\ta\|_{L^p(M)} \bigr)  \\
&\leq C_V  \bigl( \| \rd_\tA^+\ta \|_{L^p(M)} 
                + \| \rd_\tA^*\ta \|_{L^p(M)} 
                + \| \ta \|_{L^p(M)}     \bigr). \\
\end{align*}
So it remains to prove the estimate for $\ta$ near the boundary 
$\pd M=S^1\times\Si$.
For that purpose we can use the decomposition
$\t^*\ta = \ph\ds + \psi\dt + \a$, where
$\ph,\psi\in W^{1,p}(S^1\times[0,\D)\times\Si,\cg)$ and
$\a\in W^{1,p}(S^1\times[0,\D)\times\Si,\rT^*\Si\otimes\cg)$.
Let $\Om:=S^1\times[0,\frac 34 \D]$ and let $K:=S^1\times[0,\frac\D 2]$.
Then we will prove the estimate for $\ph$ and $\psi$ on $\Om\times\Si$ and for
$\a$ on $K\times\Si$.

Firstly, note that $\psi=\ta(\t_*\pd_t)\comp\t$, where $-\t_*\pd_t|_{\pd M}=\n$ 
is the outer unit normal to $\pd M$. 
So one can cut off $\t_*\pd_t$ outside of $\t(\Om\times\Si)$ to obtain a vector 
field $X\in\G(\rT M)$ that satisfies the assumption of 
lemma~\ref{Hodge prop}, that is $X|_{\pd M}=-\n$ is perpendicular to 
the boundary. 
Then one has to test (\ref{2nd eqn}) with $\l=\p\cdot \i_X\tg$
for all $\p\in\cC^\infty_\d(M,\cg)$. Again both boundary terms 
vanish. Indeed, on $S^1\times\pd Y$ we have $\p\equiv 0$ and $\i_X\tg=\t_*\dt$,
hence $\rd\l|_{\R\times \pd Y}=0$ and
$*\rd\l|_{\R\times\pd Y} = -\frac{\pd\p}{\pd\n} * \t_*(\dt\wedge\dt)=0$. 
Thus lemma~\ref{Hodge prop} yields the following estimate.
\begin{align*}
\|\psi\|_{W^{1,p}(\Om\times\Si)} 
&\leq C \|\ta(X)\|_{W^{1,p}(M)} \\
&\leq C \bigl( \|f\|_{L^p(M)} + \|\g\|_{L^p(M)} + \|\ta\|_{L^p(M)} \bigr)  \\
&\leq C_t  \bigl( \| \rd_\tA^+\ta \|_{L^p(M)} 
                  + \| \rd_\tA^*\ta \|_{L^p(M)} 
                  + \| \ta \|_{L^p(M)}     \bigr) . 
\end{align*}
Here $C$ denotes any finite constant and the bounds on the derivatives of $\t$ enter into 
the constant $C_t$.

Next, for the regularity of $\ph=\ta(\pd_s)\comp\t$ one can apply 
lemma~\ref{Hodge prop} with the tangential vector field $X=\pd_s$.
Recall that $\t$ preserves the $S^1$-coordinate.
One has to verify the second weak equation for all $\p\in\cC^\infty_\n(M,\cg)$, 
i.e.\ consider (\ref{2nd eqn}) for $\l=\p\cdot \i_X \tg = \p\cdot \ds$.
The first boundary term vanishes since one has
$*\rd\l|_{S^1\times\pd Y}=-\frac{\pd\p}{\pd\n}\dvol_{\pd Y}=0$.
For the second term one can choose any $\d>0$ and then finds a constant $C_\d$
such that for all $\p\in\cC^\infty_\n(M,\cg)$
\begin{align*}
\left| \int_{S^1\times\pd Y} \la \ta \wedge \rd \l \ra \right|
&= \left| \int_{S^1} \int_\Si \la \a(s,0) \wedge \rd_\Si(\p\comp\t)(s,0) \ra
                                                          \wedge \ds \right|  \\
&= \left| \int_{S^1\times\pd Y} \la \ta \wedge [\tA ,\p] \ra \wedge \ds \right|\\
&\leq \|\ta\|_{L^p(\pd M)}\|\tA\|_\infty\|\p\|_{L^{p^*}(\pd M)}\\
&\leq \bigl( \d \|\ta\|_{W^{1,p}(M)} + C_\d \|\ta\|_{L^p(M)} \bigr) 
      \|\p\|_{W^{1,p^*}(M)} .
\end{align*}
This uses the fact that $\a(s,0)$ and $\rd_{A_s} (\p\comp\t) |_{(s,0)\times\Si}$ 
both lie in the tangent space $\rT_{A_s}\cL$ to the Lagrangian, on which the 
symplectic form vanishes, that is
$\int_\Si \la \a \wedge \rd_A (\p\comp\t) \ra =0$.
Moreover, we have used the trace theorem for Sobolev spaces, in particular
the estimate (\ref{bdy est}).
Now lemma~\ref{Hodge prop} and remark~\ref{Hodge rmk} 
yield with $c_1=\|f\|_p$, 
$c_2 = \|\g\|_{L^p(M)} + \d \|\ta\|_{W^{1,p}(M)} + C_\d \|\ta\|_{L^p(M)}$,
and using (\ref{est f}), (\ref{est g})
\begin{align*}
&\|\ph\|_{W^{1,p}(\Om\times\Si)} \\
&\leq C \bigl( \|f\|_{L^p(M)} + c_2 + \|\ta\|_{L^p(M)} \bigr)  \\
&\leq \d \|\ta\|_{W^{1,p}(M)}
    + C_s(\d) \bigl( \| \rd_\tA^+\ta \|_{L^p(M)}  
                   + \| \rd_\tA^*\ta \|_{L^p(M)} 
                   + \| \ta \|_{L^p(M)}          \bigr) . 
\end{align*}
Here again $\d>0$ can be chosen arbitrarily small and the constant $C_s(\d)$ 
depends on this choice.

It remains to establish the estimate for the $\Si$-component $\a$ near the 
boundary.
In the coordinates $\t$ on $\Om\times\Si$, the forms
$\rd_\tA^*\ta$ and $\rd_\tA^+\ta$ become
\begin{align*}
\t^*\rd_\tA^*\ta 
&= -\pd_s\ph - \pd_t\psi + \rd_\Si^*\a - \t^*(*[\tA\wedge*\ta]), \\
\t^*\rd_\tA^+\ta 
&= \half\bigl( - (\pd_s\a + *_{\scriptscriptstyle\Si}\pd_t\a)\wedge\ds
+ *_{\scriptscriptstyle\Si}(\pd_s\a + *_{\scriptscriptstyle\Si}\pd_t\a)\wedge\dt
\bigr) \\
&\quad 
+ \half \bigl(\rd_\Si\a + (*_{\scriptscriptstyle\Si}\rd_\Si\a)\ds\wedge\dt\bigr)
+ \t^*([\tA\wedge\ta]^+) .
\end{align*}
So one obtains the following bounds: The components in the mixed direction of $\Om$ and $\S$ of
the second equation yields for some constant $C_1$ 
\begin{align*}
\|\pd_s\a+*_{\scriptscriptstyle\Si}\pd_t\a\|_{L^p(\Om\times\Si)}
&\leq \bigl\|\t^*\rd_\tA^+\ta\bigr\|_{L^p(\Om\times\Si)} 
    + \bigl\|\t^*([\tA\wedge\ta]^+)\bigr\|_{L^p(\Om\times\Si)}  \\
&\leq C_1 \bigl( \|\rd_\tA^+\ta\|_{L^p(M)} + \|\ta\|_{L^p(M)} \bigr) .
\end{align*}
Similarly, a combination of the first equation and the $\S$-component of the second equation
can be used for every $\d>0$ to find a constant $C_2(\d)$ such that
\begin{align*}
&\|\rd_\Si\a\|_{L^p(\Om\times\Si)} + \|\rd_\Si^*\a\|_{L^p(\Om\times\Si)} \\
&\leq C \bigl( \|\rd_\tA^+\ta\|_{L^p(M)} + \|\rd_\tA^*\ta\|_{L^p(M)}
            + \|\ta\|_{L^p(M)} \\
&\qquad\qquad\qquad\quad\;\;\; + \|\ph\|_{W^{1,p}(\Om\times\Si)} 
            + \|\psi\|_{W^{1,p}(\Om\times\Si)} \bigr) \\
&\leq \d \|\ta\|_{W^{1,p}(M)}
      + C_2(\d) \bigl( \| \rd_\tA^+\ta \|_{L^p(M)} 
                     + \| \rd_\tA^*\ta \|_{L^p(M)}  
                     + \| \ta \|_{M)}   \bigr) .
\end{align*}
Now firstly, lemma~\ref{Sigma reg} provides an $L^p$-estimate for the 
derivatives of $\a$ in $\Si$-direction,
\begin{align*}
&\|\nabla_\Si\a\|_{L^p(\Om\times\Si)} \\
&\leq C \bigl( \|\rd_\Si\a\|_{L^p(\Om\times\Si)} 
             + \|\rd_\Si^*\a\|_{L^p(\Om\times\Si)} 
             + \|\a\|_{L^p(\Om\times\Si)} \bigr) \\
&\leq \d \|\ta\|_{W^{1,p}(M)}
      + C_\Si(\d) \bigl( \| \rd_\tA^+\ta \|_{L^p(M)}
                      + \| \rd_\tA^*\ta \|_{L^p(M)}  
                      + \| \ta \|_{L^p(M)}   \bigr) ,
\end{align*}
where again $C_\Si(\d)$ depends on the choice of $\d>0$.
For the derivatives in $s$- and $t$-direction, we will now apply 
\cite[Theorem~1.3]{W Cauchy} on the Banach space
$X=L^p(\Si,\rT^*\Si\otimes\cg)$ with the complex structure
$*_{\scriptscriptstyle\Si}$ determined by the metric $g_{s,t}$ on $\Si$ and hence 
depending smoothly on $(s,t)\in\Om$.
The Lagrangian submanifold $\cL\subset X$ is totally real with respect to all Hodge operators
and it is modelled on a closed subspace of $L^p(\S,\R^n)$ as seen in 
\cite[Lemma~4.2, Corollary~4.4]{W Cauchy}
Now $\a\in W^{1,p}(\Om,X)$ satisfies the Lagrangian boundary condition
$\a(s,0)\in\rT_{A_s}\cL$ for all $s\in S^1$, where 
$s\mapsto A_s$ is a smooth loop in $\cL$.
Thus \cite[Corollary~1.4]{W Cauchy} yields a constant $C$
such that the following estimate holds,
\begin{align*}
\|\nabla_\Om\a\|_{L^p(K\times\Si)}
&\leq \| \a \|_{W^{1,p}(K,X)} \\
&\leq C \bigl( \| \pd_s\a + *_{\scriptscriptstyle\Si}\pd_t\a \|_{L^p(\Om,X)} 
             + \| \a \|_{L^p(\Om,X)} \bigr) \\
&\leq C_K \bigl( \| \rd_\tA^+\ta \|_{L^p(M)}
               + \| \ta \|_{L^p(M)}   \bigr) .
\end{align*}
Here $C_K$ also includes the above constant $C_1$.
Now adding up all the estimates for the different components of $\ta$ gives
for all $\d>0$
\begin{align*}
\|\ta\|_{W^{1,p}}  
&\leq (C_V+C_t+C_s(\d)+C_\Si(\d)+C_K) 
     \bigl( \| \rd_\tA^+\ta \|_p + \| \rd_\tA^*\ta \|_p + \| \ta \|_p   \bigr) \\
&\quad +2\d \|\ta\|_{W^{1,p}} .
\end{align*}
Finally, choose $\d=\frac 14$, then the term $\|\ta\|_{W^{1,p}}$ 
can be absorbed into the left hand side and this finishes the proof of the lemma.
\QED \\

\noindent
{\bf Proof of theorem~C~(iii) : } \\
Let $\b\in L^q(S^1\times Y,\rT^*Y\otimes\cg)$ and $\z\in L^q(S^1\times Y,\cg)$
be as supposed in theorem~C. Then there exists a constant $C$ such that for all
$\a\in E_A^{1,p}$ and $\ph\in W^{1,p}(S^1\times Y,\cg)$
\begin{align}
&\left| \int_{S^1}\int_Y \la \nabla_s\a -\rd_A\ph + *\rd_A\a \,,\, \b \ra 
    \;+ \int_{S^1}\int_Y \la \nabla_s\ph -\rd_A^*\a \,,\, \z \ra  \right| 
\nonumber\\
&= \left| \int_{S^1\times Y} \la D_{(A,\P)}(\a,\ph) \,,\, (\b,\z) \ra  \right|
\nonumber\\
&\leq C \|(\a,\ph)\|_{q^*} . \label{test}
\end{align}
The higher regularity of $\z$ is most easily seen if we go back to the notation
$\ta=\a+\ph\ds$ with $D_{(A,\P)}(\a,\ph)=( 2\g \,,\, -\rd_\tA^*\ta )$ for
$\rd_\tA^+\ta=*\g - \g\wedge\ds$ .
Abbreviate $M:=S^1\times Y$, then we have for all 
$\ta\in\cC^\infty(M,\rT^*M\otimes\cg)$ with $*\ta|_{\pd M}=0$
and $\ta|_{\{s\}\times\pd Y}\in\rT_{A_s}\cL$ for all $s\in S^1$
$$
\left| \int_M \la 2 \rd_\tA^+\ta \,,\, \b\wedge\ds \ra 
     + \int_M \la \rd_\tA^*\ta \,,\, \z \ra \right|
\leq C \|\ta\|_{q^*} .
$$
Now use the embedding $\t:S^1\times[0,\D)\times\Si \hookrightarrow M$ to
construct a connection $\hat A\in\cA(M)$ such that $\t^*\hat A (s,t,z) = A_s(z)$ 
near the boundary (this can be cut off and then extends trivially to all of $M$).
Then $\ta:=\rd_{\hat A}\p$ satisfies the above boundary conditions for all
$\p\in\cC^\infty_\n(M,\cg)$ since
$\rd_{\hat A}\p (\n) = \frac{\pd\p}{\pd\n} + [\hat A(\n),\p] = 0$
and $\rd_{\hat A}\p|_{\{s\}\times\pd Y} = \rd_{A_s}\p \in \rT_{A_s}\cL$
for all $s\in S^1$.
Thus we obtain for all $\p\in\cC^\infty_\n(M,\cg)$ in view of 
$\laplace\p = \rd^*(\ta-[\hat A,\p])$, denoting all constants by $C$,
\begin{align*}
&\left| \int_M \la \laplace \p \,,\, \z \ra \right| \\
&= \left| \int_M \la \rd_\tA^*\ta + *[\tA\wedge *\ta] - \rd^*[\hat A,\p]
               \,,\, \z \ra \right| \\
&\leq C \|\ta\|_{q^*} 
+ \left| \int_M \la -2 \rd_\tA^+\rd_{\hat A}\p \,,\, \b\wedge\ds \ra \right|
+ \left| \int_M \la  *[\tA\wedge *\rd_{\hat A}\p] - \rd^*[\hat A,\p]\,,\, \z \ra 
  \right| \\
&\leq C \|\p\|_{W^{1,q^*}} .
\end{align*}
The regularity theory for the Neumann problem, e.g.\ proposition~\ref{Laplace reg},
then asserts that $\z\in W^{1,q}(M)$.

To deduce the higher regularity of $\b$ we will mainly use 
lemma~\ref{Hodge prop}. The first weak equation in the lemma is given by choosing
$\a=0$ in (\ref{test}). For all $\e\in\cC^\infty(M,\cg)$
\begin{align*}
\left| \int_M \la \b \,,\, \rd\e \ra \right|
&= \left| \int_{S^1}\int_Y \la \b \,,\, \rd_A\e - [A,\e] \ra \right| \\
&\leq C \|\e\|_{q^*} + \left| \int_{S^1}\int_Y \la \nabla_s\z \,,\, \e \ra \right|
      + \left| \int_{S^1}\int_Y \la [\b\wedge *A] \,,\, \e \ra \right| \\
&\leq C \|\e\|_{q^*} .
\end{align*}
For the second weak equation let $\ph=0$ and $\a=*\rd\l - \pd_s\l$
for $\l=\p\cdot \i_X \tg$ with $\p\in\cT$ in the function space 
$\cC^\infty_\d(M,\cg)$ or $\cC^\infty_\n(M,\cg)$
corresponding to the vector field $X\in\cC^\infty(M,\rT Y)$. 
If the boundary conditions for $\a\in E_A^{1,p}$ are 
satisfied, then we obtain with $\rd=\rd_Y$
\begin{align*}
&\left| \int_M \la \b \,,\, \rd_M^*\rd_M\l \ra  \right| \\
&= \left| \int_{S^1}\int_Y \la \b \,,\, 
          *\rd*\rd\l - \pd_s^2\l - *(\pd_s*)\pd_s\l \ra \right| \\
&= \left| \int_{S^1}\int_Y \la \b \,,\, 
          *\rd_A\a - *[A\wedge*\rd\l] + *\rd_A\pd_s\l  \right.\\ 
& \qquad\qquad\qquad 
  +\nabla_s\a - [\P,\pd_s\l] -\nabla_s * \rd\l - *(\pd_s*)\pd_s\l \ra  \biggr| \\
&\leq C \|\l\|_{W^{1,q^*}} 
+ \left| \int_{S^1}\int_Y \la \z \,,\, \rd_A^*\a \ra \right|
+ \left| \int_{S^1}\int_Y \la \b \,,\, *\rd_A\pd_s\l - \nabla_s * \rd\l \ra \right|
\\
&\leq C \|\p\|_{W^{1,q^*}} .
\end{align*}
Here we have used the identity
$$
*\rd_A\pd_s\l - \nabla_s * \rd\l 
= *[A\wedge\pd_s\l] - [\P,*\rd\l] -(\pd_s *)\rd\l .
$$
Moreover, we have used partial integration with vanishing boundary term \hbox{$*\a|_{\pd Y}=0$} 
to obtain
$$
\int_{S^1}\int_Y \la \z \,,\, \rd_A^*\a \ra 
= \int_{S^1}\int_Y \la \rd_A\z \,,\, *\rd\l - \pd_s\l \ra .
$$
Now let $X\in\cC^\infty(M,\rT Y)$ be perpendicular to the boundary 
$\pd M=S^1\times\pd Y$, then for all $\p\in\cC_\d^\infty(M)$ the boundary 
conditions for $\a=*\rd\l - \pd_s\l\in E_A^{1,p}$ are satisfied.
Indeed, on the boundary $\pd M=S^1\times\pd Y$ the $1$-form $\l=\p\cdot \i_X\tg$ vanishes,
we have $\i_X\tg = h \cdot \t_*\dt$ for some smooth function $h$, and moreover
$\rd\p=-\frac{\pd\p}{\pd\n} \cdot \t_*\dt$. Hence 
\begin{align*}
*\a|_{\pd Y} & =  \rd\l|_{\pd Y} - * \pd_s\l|_{\pd Y}   = 0 , \\
 \a|_{\pd Y} & = *\rd\l|_{\pd Y} -   \pd_s\l|_{\pd Y}   
               = - \tfrac{\pd\p}{\pd\n} h * (\t_*\dt \wedge \t_*\dt) = 0 .
\end{align*}
Thus for all vector fields $X\in\cC^\infty(M,\rT Y)$ that are perpendicular to 
the boundary, lemma~\ref{Hodge prop} asserts that $\b(X)\in W^{1,q}(M,\cg)$.
In particular, this implies $W^{1,q}$-regularity of $\b$ on all compact subsets
$K\subset{\rm int}\,M$. So it remains to prove the regularity on the neighbourhood
$\t(S^1\times[0,\frac \D 2)\times\Si)$ of the boundary $\pd M$.
In these coordinates we decompose 
$$
\t^*\b = \x\dt + \hat\b .
$$
Now firstly, lemma~\ref{Hodge prop} applies as described above to assert the regularity
$\x=\b(\t_*\pd_t)\comp\t\in W^{1,q}(\Om\times\Si,\cg)$ on $\Om:=S^1\times[0,\tfrac 34\D)$.
Here a vector field $X$ that is perpendicular to the boundary is constructed by cutting off 
$\t_*\pd_t$ outside of $\t(\Om\times\Si)$.

So it remains to consider $\hat\b\in L^q(\Om\times\Si,\rT^*\Si\otimes\cg)$
and establish its $W^{1,q}$-regularity.
In order to derive a weak inequality for $\hat\b$ on $\Om\times\Si$ from (\ref{test}) we 
use $\ta=\t_*(\ph\ds + \psi\dt + \hat\a )$ with $\ph\in\cC^\infty_\d(\Om\times\Si,\cg)$, 
$\psi\in\cC^\infty_\d(\Om\times\Si,\cg)$, and 
$\hat\a\in W^{1,p}(\Om\times\Si,\rT^*\Si\otimes\cg)$ such that 
$\hat\a (s,\tfrac 34\D,\cdot)= 0$ and $\hat\a (s,0,\cdot)\in\rT_{A_s}\cL$ for all $s\in S^1$. 
This $\ta$ satisfies the boundary conditions for (\ref{test}) and it can be
extended trivially to a $W^{1,p}$-regular $1$-form on all of $M$.
Thus we obtain
\begin{align*}
&\left| \int_{\Om\times\Si} 
    \la \nabla_s\hat\a + *\nabla_t\hat\a - \rd_A\ph -*\rd_A\psi \,,\, \hat\b \ra  
 \right| \\
&\leq \left| \int_{\Om\times\Si} \la - \nabla_s\psi + \nabla_t\ph - *\rd_A\hat\a
                               \,,\, \x \ra \;
    + \int_{\Om\times\Si} \la \nabla_s\ph + \nabla_t\psi - \rd_A^*\hat\a
                  \,,\, \z \ra \right| \\
&\quad + C \| \ph\,\ds + \psi\,\dt + \hat\a \|_{q^*} .
\end{align*}
Here we have introduced the decomposition $\t^*\tA=\P\ds+\Psi\dt+A$, where
$A\in\cC^\infty(\Om\times\Si,\rT^*\Si\otimes\cg)$ with $A(s,0)=A_s\in\cL$ for
all $s\in S^1$. We have also used the notation $\nabla_t\ph = \pd_t\ph + [\Psi,\ph]$, and 
moreover $\rd_A$ and $*$ denote the differential and Hodge operator on $\Si$.
Now firstly put $\hat\a=0$, then we obtain for all 
$\ph,\psi\in\cC^\infty_\d(\Om\times\Si,\cg)$ by partial integration
\begin{align*}
\left| \int_{\Om\times\Si} \la \rd_A\ph \,,\, \hat\b \ra \right|
&\leq \bigl( C + \| \nabla_t \x - \nabla_s \z \|_{L^q(\Om\times\Si)} \bigr)
      \| \ph \|_{L^{q^*}(\Om\times\Si)} , \\
\left| \int_{\Om\times\Si} \la *\rd_A\psi \,,\, \hat\b \ra  \right|
&\leq \bigl( C + \| \nabla_s \x + \nabla_t \z \|_{L^q(\Om\times\Si)} \bigr)
      \| \psi \|_{L^{q^*}(\Om\times\Si)} .
\end{align*}
This shows that the weak derivatives $\rd_\Si^*\hat\b$ and $\rd_\Si\hat\b$ are
of class $L^q$ on $\Om\times\Si$. So we have verified the assumptions of lemma~\ref{Sigma reg} 
for $\hat\b$ and conclude that $\nabla_\Si\hat\b$ also is of class $L^q$. 
So it remains to deduce the $L^q$-regularity of $\pd_s\hat\b$ and $\pd_t\hat\b$ on 
$S^1\times[0,\frac \D 2]\times\Si$ from the above inequality for $\ph=\psi=0$, namely from
\begin{equation}\label{hat a eqn}
\left| \int_{\Om\times\Si} \la \nabla_s\hat\a + *\nabla_t\hat\a \,,\, \hat\b \ra 
\right|
\leq \bigl( C + \| \rd_A\z + *\rd_A\x \|_{L^q(\Om\times\Si)} \bigr)
      \| \hat\a \|_{L^{q^*}(\Om\times\Si)} .
\end{equation}
This holds for all $\hat\a\in W^{1,p}(\Om\times\Si,\rT^*\Si\otimes\cg)$ 
such that $\hat\a (s,\tfrac 34\D,\cdot)= 0$ 
and $\hat\a (s,0,\cdot)\in\rT_{A_s}\cL$ for all $s\in S^1$. 
We now have to employ different arguments according to whether $q>2$ or $q<2$.\\

\noindent
\underline{\bf Case $\mathbf{q>2}$ :}\\
In this case the regularity of $\pd_s\hat\b$ and $\pd_t\hat\b$
will follow from \cite[Theorem~1.3]{W Cauchy} on the Banach 
space $X=L^q(\Si,\rT^*\Si\otimes\cg)$ with the complex structure given by the Hodge
operator on $\Si$ with respect to the metric $g_{s,t}$.
From (\ref{hat a eqn}) one obtains the following estimate for some constant $C$ 
and all $\hat\a$ as above:
\begin{equation} \label{test b}
\left| \int_\Om \int_\Si \la \hat\b \,,\, \pd_s\hat\a + \pd_t(*\hat\a) \ra \right|
\leq C \|\hat\a\|_{L^{q^*}(\Om,X^*)} .
\end{equation}
Note that this extends to the $W^{1,q^*}(\Om,L^{q^*}(\Si))$-closure of the 
admissible $\hat\a$ from above. In particular the estimate above holds for all 
$\hat\a\in W^{1,q}(\Om,X)$ that are supported in $\Om$ and satisfy 
$\hat\a (s,0,\cdot)\in\rT_{A_s}\cL$ for all $s\in S^1$. 
To see that these can be approximated by smooth $\hat\a$ with Lagrangian 
boundary conditions one uses the Banach submanifold coordinates for $\cL$ 
given by \cite[Lemma~4.3]{W Cauchy} as before.
Here the Lagrangian submanifold $\cL\subset X$ is totally real with respect to all 
Hodge operators as before, and it is the $L^q$-restriction or -completion of the original 
submanifold in $\cA^{0,p}(\S)$, hence it is modelled on $W^{1,q}_z(\S,\cg)\times\R^m$, a 
closed subspace of $L^q(\S,\R^n)$ (see \cite[Lemma~4.2,~4.3]{W Cauchy}).
However, in order to be able to apply \cite[Theorem~1.3]{W Cauchy}, we need to extend 
this estimate to all $\hat\a\in W^{1,\infty}(\overline{\Om},X^*)$ with $\supp\a\subset\Om$ and 
$\a(s,0)\in (*\rT_{A_s}\cL)^\perp$ for all $s\in S^1$.
This is possible since any such $\hat\a$ can be approximated in 
$W^{1,q^*}(\Om,X^*)$ by $\hat\a_i\in\cC^\infty(\Om,X)$ that are compactly 
supported in $\Om$ and satisfy the above stronger boundary condition 
$\hat\a_i (s,0)\in\rT_{A_s}\cL$ for all $s\in S^1$. 

Indeed, \cite[Lemma~2.2]{W Cauchy} provides such an approximating 
sequence $\a_i$ without the Lagrangian boundary conditions. From the proof via mollifiers 
one sees that the approximating sequence can be chosen with compact support in $\Om$.
Now for all $s\in S^1$ one has the topological splitting
$X=\rT_{A_s}\cL \oplus *\rT_{A_s}\cL$ and thus
$X^*=(\rT_{A_s}\cL)^\perp \oplus (*\rT_{A_s}\cL)^\perp$.
Since $q>2$ the embedding $X\hookrightarrow X^*$ is continuous.
This identification uses the $L^2$-inner product on $X$ which equals the metric 
$\o(\cdot,*\cdot)$ given by the symplectic form $\o$ and the complex structure $*$.
So due to the Lagrangian condition this embedding maps
$\rT_{A_s}\cL \hookrightarrow (*\rT_{A_s}\cL)^\perp$ and
$*\rT_{A_s}\cL \hookrightarrow (\rT_{A_s}\cL)^\perp$.
We write $\hat\a=\g+\d$ and $\a_i=\g_i+\d_i$ according to these splittings
to obtain $\g,\d\in\cC^\infty(\Om,X^*)$ and $\g_i,\d_i\in\cC^\infty(\Om,X)$
such that $*\rT_A\cL \ni \g_i\to\g \in (\rT_A\cL)^\perp$ 
and $\rT_A\cL\ni\d_i\to\d\in (*\rT_A\cL)^\perp$ with convergence in $W^{1,q^*}(\Om,X^*)$.
The boundary condition on $\hat\a$ gives $\g|_{t=0}\equiv 0$.
Moreover, $\pd_t\g$ is uniformly bounded in $X^*$, so one can find a constant
$C$ such that $\|\g(s,t)\|_{X^*}\leq Ct$ for all $t\in [0,\frac 34\D]$ and
hence for sufficiently small $\ep>0$
$$
\|\g\|_{L^{q^*}(S^1\times[0,\ep],X^*)}\leq \tfrac C{1+q^*}\ep^{1+\frac 1{q^*}} .
$$ 
Now let $\d>0$ be given and choose $1>\ep>0$ such that 
$\|\g\|_{L^{q^*}(S^1\times[0,\ep],X^*)} \leq \ep\d$
and $\|\g\|_{W^{1,q^*}(S^1\times[0,\ep],X^*)}\leq \d$.
Next, choose $i\in\N$ sufficiently large such that
$\|\g_i-\g\|_{W^{1,q^*}(\Om,X^*)}\leq \ep\d$, and let 
$h\in\cC^\infty([0,\frac 34 \D],[0,1])$ be a cutoff function with
$h(0)=0$, $h|_{t\geq\ep}\equiv 0$, and $|h'|\leq \frac 2 \ep$.
Then \hbox{$\hat\a_i:=h\g_i+\d_i\in\cC^\infty(\Om,X)$} \hbox{satisfies}
 the Lagrangian
boundary condition $\hat\a_i(s,0)\in\rT_{A_s}\cL$ and approximates $\hat\a$ 
in view of the following estimate,
\begin{align*}
\|\hat\a_i-\hat\a\|_{W^{1,q^*}(\Om,X^*)}
&\leq \| h (\g_i-\g) \|_{W^{1,q^*}(\Om,X^*)}
     +\| (1-h) \g \|_{W^{1,q^*}(\Om,X^*)}  \\
&\leq \| \g_i-\g \|_{W^{1,q^*}(\Om,X^*)}
     +\tfrac 2\ep \| \g_i-\g \|_{L^{q^*}(\Om,X^*)}  \\
&\qquad     +\| \g \|_{W^{1,q^*}(S^1\times[0,\ep],X^*)} 
     +\tfrac 2\ep \| \g \|_{L^{q^*}(S^1\times[0,\ep],X^*)}  \\
&\leq 6\d .
\end{align*}
This approximation shows that (\ref{test b}) holds indeed true for all
$\hat\a\in\cC^1(\overline\Om,X^*)$ with $\supp\a\subset\Om$ and $\a(s,0)\in (*\rT_{A_s}\cL)^\perp$ 
for all $s\in S^1$.
Thus \cite[Theorem~1.3]{W Cauchy} asserts that 
$\hat\b\in W^{1,q}(K,X)$ for $K:=S^1\times[0,\frac \D 2]$, and hence 
$\pd_s\hat\b$ and $\pd_t\hat\b$ are of class $L^q$ on $S^1\times[0,\frac \D 2]\times\Si$ 
as claimed.\\

\noindent
\underline{\bf Case $\mathbf{q<2}$ :}\\
In this case we cover $S^1$ by two intervals, $S^1=I_1\cup I_2$ such that 
there are isometric embeddings $(0,1)\hookrightarrow S^1$ identifying 
$[\tfrac 14,\tfrac 34]$ with $I_1$ and $I_2$ respectively.
Abbreviate $K:=[\tfrac 14,\tfrac 34]\times[0,\frac \D 2]$ and let 
$\Om'\subset (0,1)\times[0,\frac 34\D]$ be a compact submanifold of the half
space $\H$ such that $K\subset{\rm int}\,\Om'$.
Then for each of the above identifications $S^1\setminus\{pt\}\cong (0,1)$ one 
has $L^q$-regularity of $\hat\b$ on $\Om'\times\Si$ by assumption and of $*\rd_A\x + \rd_A\z$
from above.
Now the task is to establish in both cases the $L^q$-regularity of 
$\pd_s\hat\b$ and $\pd_t\hat\b$ on $K\times\Si$ using (\ref{hat a eqn}). 
For that purpose choose a cutoff function $h\in\cC^\infty(\H,[0,1])$ supported
in $\Om'$ such that $h|_K\equiv 1$. Then it suffices to find a constant $C$ 
such that for all $\g\in\cC^\infty_0(\Om'\times\Si,\rT^*\Si\otimes\cg)$ 
(these are compactly supported in ${\rm int}(\Om')\times\Si$)
$$
\left| \int_{\Om'\times\Si} \la \pd_s\g \,,\, h\hat\b \ra \right|
\leq C \|\g\|_{q^*} .
$$
This gives $L^q$-regularity of the weak derivative $\pd_s(h\hat\b)$ and hence
of $\pd_s\hat\b$ on $K\times\Si$. For the regularity of $\pd_t\hat\b$ one has to 
replace $\pd_s\g$ by $\pd_t\g$, then the argument is the same as the following argument
for $\pd_s\hat\b$.

We linearize the submanifold chart maps along 
$(A_s)_{s\in(0,1)}\in\cL\cap\cA(\Si)$ given by \cite[Lemma~4.3]{W Cauchy} 
for the Lagrangian $\cL\subset\cA^{0,q^*}(\Si)$. Note that this uses the 
$L^{q^*}$-completion of the actual Lagrangian in $\cA^{0,p}(\Si)$.
Abbreviate $Z:= W^{1,q^*}_z(\Si,\cg)\times\R^m$ and let $*_{s,t}$ denote the Hodge 
operator on $\Si$ with respect to the metric $g_{s,t}$.
Then one obtains a smooth family of bounded isomorphisms 
$$
\Th_{s,t}: \;
Z \times Z \overset{\sim}{\longrightarrow} L^{q^*}(\Si,\rT^*\Si\otimes\cg) \;=:\, X
$$
defined for all $(s,t)\in\Om'$ by
$$
\Th_{s,t}(\x,v,\z,w) \;=\; 
\rd_{A_s}\x + \tsum_{i=1}^m v^i  \g_i(s) 
+ *_{s,t}\rd_{A_s}\z + \tsum_{i=1}^m w^i *_{s,t}\g_i(s) .
$$
Here $\g_i\in\cC^\infty((0,1)\times\Si,\rT^*\Si\otimes\cg)$ with 
$\g_i(s)\in\rT_{A_s}\cL$ for all $s\in(0,1)$.
Abbreviate $Z^\infty:=\cC^\infty_z(\Si,\cg)\times\R^m\subset Z$, then $\Th_{s,t}$ maps
$Z^\infty\times Z^\infty$ into the set of smooth $1$-forms $\Om^1(\Si,\cg)$.
So given any $\g\in\cC^\infty_0(\Om'\times\Si,\rT^*\Si\otimes\cg)$ 
we have $f:=\Th^{-1}\comp\g\in\cC^\infty_0(\Om',Z^\infty\times Z^\infty)$
and for some constant $C$
$$
\|f\|_{L^{q^*}(\Om',Z\times Z)}
\;\leq\; C \|\g\|_{L^{q^*}(\Om',X)}
\;=\; C \|\g\|_{L^{q^*}(\Om'\times\Si)}.
$$
Write $f=(f_1,f_2)$ with $f_i\in\cC^\infty_0(\Om',Z^\infty)$ and note that 
$\int_{\Om'} \pd_s f_1 = 0$ due to the compact support.
So one can solve $\laplace_{\Om'}\p_1 = \pd_s f_1$ by 
$\p_1\in\cC^\infty_\n(\Om',Z^\infty)$ with $\int_{\Om'}\p_1=0$ 
and $\laplace_{\Om'}\p_2 = \pd_s f_2$ by $\p_2\in\cC^\infty_\d(\Om',Z^\infty)$.
(For the $\cC^\infty_z(\Si,\cg)$-component of $Z^\infty$ one has solutions of the 
Laplace equation on every $\Om'\times\{x\}$ that depend smoothly on $x\in\Si$.)
Now let $\P:=(\p_1,\p_2)\in\cC^\infty(\Om',Z\times Z)$ and consider the $1$-form
$$
\hat\a_\g \,:=\; h\cdot \Th(-\pd_s\P + J_0 \pd_t \P)
\;\in\;\cC^\infty(\Om',X) .
$$
This extends to a $1$-form on $\Om\times\Si$ that is admissible in 
(\ref{hat a eqn}). Indeed, $\hat\a_\g$ vanishes for $s$ close to $0$ or $1$ and 
thus trivially extends to $s\in S^1$. The Lagrangian boundary condition is met 
since for all $s\in S^1$
$$
\hat\a_\g(s,0)
=h(s,0) \cdot \Th_{s,0}(-\pd_s\p_1 - \pd_t\p_2 , -\pd_s\p_2 + \pd_t\p_1)
\in \Th_{s,0}(Z,0) = \rT_{A_s}\cL .
$$
So (\ref{hat a eqn}) provides a constant $C$ such that for all $\hat\a_\g$ of the above form 
$$
\left| \int_{\Om'\times\Si} \la \hat\b \,,\, \pd_s\hat\a_\g + \pd_t(*\hat\a_\g) \ra 
\right|
\;\leq\; C \| \hat\a_\g \|_{L^{q^*}(\Om,X)}
$$
Moreover, one has for all $\g\in\cC^\infty_0(\Om'\times\Si,\rT^*\Si\otimes\cg)$ 
and the associated $f$, $\P$ and $\hat\a_\g$ and denoting all constants by $C$
$$
\| \hat\a_\g \|_{L^{q^*}(\Om,X)}
\;\leq\; C \| \P \|_{W^{1,q^*}(\Om',Z\times Z)}
\;\leq\; C \| f \|_{L^{q^*}(\Om',Z\times Z)} 
\;\leq\; C \| \g \|_{L^{q^*}(\Om'\times\Si)} .
$$
Here the second inequality follows from $\laplace_{\Om'}\P=\pd_s f$ and 
\cite[Lemma~2.1]{W Cauchy} as follows.
In the $\R^m$-component of $Z$, this is the usual elliptic estimate for the 
Dirichlet or Neumann problem; for the components in the infinite dimensional part 
$Y:=W^{1,q^*}_z(\Si,\cg)$ of $Z$ (still denoted by $\p_i$ and $f_i$) this uses the following 
estimate.
For all $\psi\in\cC^\infty_\n(\Om',Y^*)$ in the case $i=1$ and 
for all $\psi\in\cC^\infty_\d(\Om',Y^*)$ in the case $i=2$ 
\begin{align*}
\left| \int_{\Om'\times\Si} \la \p_i \,,\, \laplace_{\Om'}\psi \ra \right|\;
&=\; \left| \int_{\Om'\times\Si} \la \laplace_{\Om'}\p_i \,,\, \psi \ra \right|  
\;=\; \left| \int_{\Om'\times\Si} \la \pd_s f_i \,,\, \psi \ra \right|  \\
&=\; \left| \int_{\Om'\times\Si} \la f_i \,,\, \pd_s\psi \ra \right|  
\;\leq\; \| f_i \|_{L^{q^*}(\Om',Y)} \| \psi \|_{W^{1,q}(\Om',Y^*)}  .
\end{align*}
Now a calculation shows that
$$
\pd_s\hat\a_\g + \pd_t(*\hat\a_\g)
= h \cdot \Th(\laplace\P) + \pd_s(h\cdot\Th)(-\pd_s\P + J_0 \pd_t \P)
  + \pd_t(h\cdot\Th)(\pd_t\P - J_0 \pd_s \P) .
$$
We then use $\laplace\P=\pd_s f$ to obtain, denoting all constants by $C$,
\begin{align*}
& \left| \int_{\Om'\times\Si} \la h\cdot\hat\b \,,\, \pd_s\g \ra \right| \\
&= \left| \int_{\Om'\times\Si} \la \hat\b \,,\, 
                    h \cdot \Th(\laplace\P) + h\cdot \pd_s\Th (f) \ra \right|\\
&\leq \left| \int_{\Om'\times\Si} \la \hat\b \,,\, 
                 \pd_s\hat\a_\g + \pd_t(*\hat\a_\g) \ra \right|  \\
&\quad + C \| \hat\b \|_{L^q(\Om',X^*)} \bigl( 
    \| -\pd_s\P + J_0 \pd_t \P) \|_{L^{q^*}(\Om',Z\times Z)}
   +\| f \|_{L^{q^*}(\Om',Z\times Z)} \bigr) \\
&\leq C \| \g \|_{L^{q^*}(\Om'\times\Si)} .
\end{align*}
This proves the $L^q$-regularity of $\pd_s\hat\b$ (and analogously of 
$\pd_t\hat\b$) on $S^1\times[0,\frac \D 2]\times\Si$ in the case $q<2$ and thus 
finishes the proof of theorem~C~(iii).
\QED \\

\noindent
{\bf Proof of theorem~C~(i) : } \\
Lemma~\ref{estimate} and the subsequent remark imply that for some constant $C$ and for all 
$(\a,\ph)$ in the domain of $D_{(A,\P)}$
$$
\|(\a,\ph)\|_{W^{1,p}} \leq C \bigl( \| D_{(A,\P)} (\a,\ph) \|_p 
                                   + \| (\a,\ph) \|_p  \bigr).
$$
Note that the embedding $W^{1,p}(X)\hookrightarrow L^p(X)$ is compact, 
so this estimate asserts that $\ker D_{(A,\P)}$ is finite dimensional
and $\im D_{(A,\P)}$ is closed (see e.g.\ \cite[3.12]{Zeidler}).
So it remains to consider the cokernel of $D_{(A,\P)}$. 
We abbreviate $Z:=L^p(S^1\times Y,\rT^*Y\otimes\cg) \times L^p(S^1\times Y,\cg)$, then
${\rm coker}\,D_{(A,\P)} = Z / \im D_{(A,\P)}$ is a Banach space since 
$\im D_{(A,\P)}$ is closed. So it has the same dimension as its dual space
$(Z/ \im D_{(A,\P)})^* \cong (\im D_{(A,\P)})^\perp$.
Now let $\s:S^1\times Y \to S^1\times Y$ denote the reflection $\s(s,y):=(-s,y)$ 
on $S^1\cong\R/\Z$, then we claim that there is an isomorphism
\begin{equation} \label{map}
\begin{array}{ccc}
(\im D_{(A,\P)})^\perp &\overset{\sim}{\longrightarrow}& \ker D_{\s^*(A,\P)} \\
(\b,\z) & \longmapsto & (\b\comp\s,\z\comp\s) .
\end{array}
\end{equation}
Here $D_{\s^*(A,\P)}=D_{(A',\P')}$ is the linearized operator at the reflected connection
$\s^*\tA=A'+\P'\ds$ with respect to the metric $\s^*\tg$ on $X$.
Note that $\ker D_{\s^*(A,\P)}$ is finite dimensional since the estimate in 
lemma~\ref{estimate} also holds for the operator $D_{\s^*(A,\P)}$.
So this would indeed prove that ${\rm coker}\,D_\tA$ is of finite dimension 
and hence $D_\tA$ is a Fredholm operator.

To establish the above isomorphism consider any 
$(\b,\z)\in (\im D_{(A,\P)})^\perp$, that is
$\b\in L^{p^*}(S^1\times Y,\rT^*Y \otimes\cg)$ and
$\z\in L^{p^*}(S^1\times Y,\cg)$ such that for all
$\a\in E_A^{1,p}$ and $\ph\in W^{1,p}(S^1\times Y,\cg)$
$$
\int_{S^1\times Y} \la D_{(A,\P)}(\a,\ph) \,,\, (\b,\z) \ra  = 0 .
$$
Iteration of theorem~C~(iii) implies that $\b$ and $\z$ are in fact 
$W^{1,p}$-regular:
We start with $q=p^*<2$, then the lemma asserts $W^{1,p^*}$-regularity.
Next, the Sobolev embedding theorem gives $L^{q_1}$-regularity for some 
$q_1\in(\frac 43,2)$ with $q_1>p^*$. Indeed, the Sobolev embedding holds for any 
$q_1\leq\frac{4p^*}{4-p^*}$, and $\frac 43 <\frac{4p^*}{4-p^*}$ as well as 
$p^* <\frac{4p^*}{4-p^*}$ holds due to $p^*>1$.
So the lemma together with the Sobolev embeddings can be iterated to give
$L^{q_{i+1}}$-regularity for $q_{i+1}=\frac{4q_i}{4-q_i}$ as long as 
$4>q_i>2$ or $2>q_i\geq p^*$.
This iteration yields $q_2\in(2,4)$ and $q_3>4$. Thus another iteration of the
lemma gives $W^{1,q_3}$- and thus also $L^p$-regularity of $\b$ and $\z$.
Finally, since $p>2$ the lemma applies again and asserts the claimed
$W^{1,p}$-regularity of $\b$ and $\z$.
Now by partial integration
\begin{align}
0 
&=\int_{S^1\times Y} \la D_{(A,\P)}(\a,\ph) \,,\, (\b,\z) \ra  \nonumber\\
&=\int_{S^1}\int_Y \la \nabla_s\a - \rd_A\ph + *\rd_A\a \,,\, \b \ra 
 \;+ \int_{S^1}\int_Y \la \nabla_s\ph - \rd_A^*\a \,,\, \z \ra \nonumber\\
&=\int_{S^1}\int_Y \la \a \,,\, -\nabla_s\b - \rd_A\z + *\rd_A\b \ra 
 \;+ \int_{S^1}\int_Y \la \ph \,,\, -\nabla_s\z - \rd_A^*\b \ra \nonumber\\
&\quad 
 \;- \int_{S^1}\int_\Si \la \a \wedge \b \ra 
 \;- \int_{S^1}\int_\Si \la \ph \,,\, *\b \ra  .   \label{pi}
\end{align}
Testing this with all 
$\a\in\cC_0^\infty(S^1\times Y,\rT^*Y\otimes\cg)\subset E_A^{1,p}$
and $\ph\in\cC_0^\infty(S^1\times Y,\cg)$ implies 
$-\nabla_s\b - \rd_A\z + *\rd_A\b = 0$ and $-\nabla_s\z - \rd_A^*\b = 0$.
Then furthermore we deduce $*\b(s)|_{\pd Y}=0$ for all $s\in S^1$ 
from testing with $\ph$ that run through all of $\cC^\infty(S^1\times\Si,\cg)$ on 
the boundary. Finally, $\int_{S^1}\int_\Si \la \a \wedge \b \ra = 0 $ remains from
(\ref{pi}).
Since both $\a$ and $\b$ restricted to $S^1\times\Si$ are continuous paths in 
$\cA^{0,p}(\Si)$, this implies that for all $s\in S^1$ and every 
$\a\in\rT_{A_s}\cL$
$$
0 \;=\;\int_\Si \la \a \wedge \b(s) \ra \;=\; \o (\a,\b(s)) ,
$$
where $\o$ is the symplectic structure on $\cA^{0,p}(\S)$.
Since $\rT_{A_s}\cL$ is a Lagrangian subspace, this proves
$\b(s)|_{\pd Y}\in\rT_{A_s}\cL$ for all $s\in S^1$ and thus $\b\in E_A^{1,p}$, or
equivalently $\b\comp\s\in E_{A\comp\s}$.
So $(\b\comp\s,\z\comp\s)$ lies in the domain of $D_{\s^*(A,\P)}$.
Now note that $\s^*\tA = A\comp\s - (\P\comp\s) \ds$, thus one obtains
$(\b\comp\s,\z\comp\s) \in \ker D_{\s^*(A,\P)}$ since
$$
D_{\s^*(A,\P)}(\b\comp\s,\z\comp\s)
\;=\; \bigl( ( -\nabla_s\b -\rd_A\z + *\rd_A\b ) \comp\s  \,,\, 
             ( -\nabla_s\z -\rd_A^*\b ) \comp \s  \bigr)
\;=\;0.
$$
This proves that the map in (\ref{map}) indeed maps into $\ker D_{\s^*(A,\P)}$. 
To see the surjectivity of this map consider any 
$(\b,\z) \in \ker D_{\s^*(A,\P)}$.
Then the same partial integration as in (\ref{pi}) shows that 
$(\b\comp\s,\z\comp\s)\in(\im D_{(A,\P)})^\perp$, and thus $(\b,\z)$ is the 
image of this element under the map (\ref{map}).
So this establishes the isomorphism (\ref{map}) and thus shows that $D_{(A,\P)}$ 
is Fredholm.
\QED

\appendix
\section{Dirichlet and Neumann problem}

Throughout this paper we use various regularity results for the Laplace operator.
For convenience these are summarized in this appendix.

We deal with (homogeneous) Dirichlet boundary conditions and with possibly inhomogeneous
Neumann boundary conditions. Often, the equations are formulated weakly with the help
of the following test function spaces:
\begin{align*}
\cC^\infty_\d(M) &= \bigl\{ \p\in\cC^\infty(M) \st \p|_{\pd M} =0 \bigr\},  \\
\cC^\infty_\n(M) &= \bigl\{ \p\in\cC^\infty(M) \st 
                                \tfrac{\pd\p}{\pd\n}\bigr|_{\pd M} =0 \bigr\}.
\end{align*}
Here and throughout this appendix let $M$ be a manifold with boundary.
We abbreviate $\laplace:=\rd^*\rd$, and denote by $\tfrac{\pd\p}{\pd\n}$ the Lie 
derivative in the direction of the outer unit normal.
Moreover, we use the notation $\N=\{1,2,\ldots\}$ and $\N_0=\{0,1,\ldots\}$.
The regularity theory for the Dirichlet and Neumann problem that is used in this paper
can be summarized as follows.
References are for example \cite{GT} and \cite[Theorems 2.3',3.2,D.2]{W}.

\begin{prp}\label{Laplace reg}
Let $k\in\N$, then there exists a constant $C$ such that the following holds.
Let $f\in W^{k-1,p}(M)$ and $G\in W^{k,p}(M)$ and suppose that $u\in W^{k,p}(M)$ is a weak 
solution of the Dirichlet problem (or the Neumann problem with inhomogeneous boundary conditions),
that is for all $\psi\in\cC^\infty_\d(M)$ (or for all $\psi\in\cC^\infty_\n(M)$)
$$
\int_M  u \cdot \laplace\psi\; = \int_M  f \cdot \psi  +  \int_{\pd M}  G \cdot \psi  .
$$
Then $u\in W^{k+1,p}(M)$ and
$$
\|u\|_{W^{k+1,p}} \leq C \bigl( \|f\|_{W^{k-1,p}} + \|G\|_{W^{k,p}} + \|u\|_{W^{k,p}} \bigr).
$$
In the special case $k=0$ there exists a constant $C$ such that the following holds:
Suppose that $u\in L^p(M)$ and that there exists a constant $c$ such that for all 
$\psi\in\cC^\infty_\d(M)$ (or for all $\psi\in\cC^\infty_\n(M)$)
$$
\int_M  u \cdot \laplace\psi\; \leq c \| \psi \|_{W^{1,p^*}} .
$$
Then $u\in W^{1,p}(M)$ and
$$
\|u\|_{W^{1,p}} \leq C \bigl( c + \|u\|_{L^p} \bigr).
$$
\end{prp}

We also frequently encounter Laplace equations for $1$-forms, where the components
satisfy different boundary conditions. In these cases the following lemma allows to 
obtain regularity results for the components separately.
The proof relies on the above standard regularity theory for the Laplace operator.

\begin{lem}  \label{Hodge prop}
Let $(M,g)$ be a compact Riemannian manifold (possibly with boundary), let 
$k\in\N_0$ and $1<p<\infty$. 
Let $X\in\G(\rT M)$ be a smooth vector field that is either perpendicular to the boundary, 
i.e.\ $X|_{\pd M}=h\cdot\n$ for some $h\in\cC^\infty(\pd M)$, or tangential, i.e.\ 
$X|_{\pd M}\in\G(\rT\pd M)$.
In the first case let $\cT=\cC^\infty_\d(M)$, in the latter case let
$\cT=\cC^\infty_\n(M)$.
Then there exists a constant $C$ such that the following holds:

Let $f\in W^{k,p}(M)$, $\g\in W^{k,p}(M,\L^2\rT^*M)$, and suppose that
the $1$-form $\a\in W^{k,p}(M,\rT^*M)$ satisfies
\begin{align*}
\int_M \la \a \,,\, \rd\e \ra &= \int_M f \cdot \e 
\quad\qquad\forall\e\in\cC^\infty(M), \\
\int_M \la \a \,,\, \rd^*\o \ra &= \int_M \la \g \,,\, \o \ra  
\qquad\forall \o=\rd(\p\cdot \i_X g)\,,\;\p\in\cT .
\end{align*}
Then $\a(X)\in W^{k+1,p}(M)$ and
$$
\|\a(X)\|_{W^{k+1,p}} \leq C \bigl(  \|f\|_{W^{k,p}} + \|\g\|_{W^{k,p}}
                                  + \|\a\|_{W^{k,p}} \bigr).
$$
\end{lem}

\begin{rmk} \label{Hodge rmk} \rm
In the case $k=0$ let $\frac 1p + \frac 1{p^*} = 1$, then the weak equations 
for $\a$ can be replaced by the following:
There exist constants $c_1$ and $c_2$ such that
\begin{align*}
\left| \int_M \la \a \,,\, \rd\e \ra \right| &\leq  c_1 \;\| \e \|_{p^*}
\;\qquad\qquad\forall\e\in\cC^\infty(M), \\
\left| \int_M \la \a \,,\, \rd^*\rd (\p\cdot \i_X g) \ra \right|
&\leq c_2 \;\|\p\|_{(W^{1,p^*})^*}  \qquad\forall\p\in\cT  .
\end{align*}
The estimate then becomes \;
$\|\a(X)\|_{W^{1,p}} \leq C \bigl(  c_1 + c_2 + \|\a\|_p \bigr)$ .
\end{rmk}

\noindent
{\bf Proof of lemma \ref{Hodge prop} and remark \ref{Hodge rmk} : } \\
Let $\a^\n\in\cC^\infty(M,\rT^*M)$ be an $L^p$-approximating sequence for $\a$
such that $\a^\n\equiv 0$ near $\pd M$.
Then one obtains for all $\p\in\cT$
\begin{align*}
\int_M \a(X) \cdot \laplace\p
&=\lim_{\n\to\infty}\left(
    \int_M \la \cL_X\a^\n \,,\, \rd\p \ra 
  - \int_M \la  \i_X\rd\a^\n \,,\, \rd\p \ra  \right) \\
&=\lim_{\n\to\infty}\left(
  - \int_M \la \a^\n \,,\, \cL_X\rd\p \ra 
  - \int_M \la \a^\n \,,\, {\rm div}X \cdot \rd\p \ra  \right. \\ 
&\qquad\qquad \left.
  - \int_M \la \a^\n \,,\, \i_{Y_{\rd\p}}\cL_X g \ra  
  - \int_M \la \rd\a^\n \,,\, \i_X g \wedge \rd\p \ra  \right) \\
&=  \int_M \la \a \,,\, \rd( - \cL_X\p - {\rm div}X \cdot \p ) \ra 
  - \int_M \la \a \,,\, \rd^*( \i_X g \wedge \rd\p ) \ra   \\
&\quad
  + \int_M \la \a \,,\, \p \cdot \rd({\rm div}X)  - \i_{Y_{\rd\p}}\cL_X g \ra \\
&=  \int_M \la f  \,,\,  - \cL_X\p - {\rm div}X \cdot \p \ra
  + \int_M \la \g \,,\, \rd(\p\cdot \i_X g ) \ra \\
&\quad
  - \int_M \la \a \,,\, \rd^*(\p\cdot \rd \i_X g ) \ra 
  + \int_M \la \a \,,\, \p \cdot \rd({\rm div}X) - \i_{Y_{\rd\p}}\cL_X g \ra .
\end{align*}
Here the vector field $Y_{\rd\p}$ is given by $\i_{Y_{\rd\p}} g = \rd\p$.
In the case $k\geq 1$ further partial integration yields for all $\p\in\cT$
$$
\int_M \a(X) \cdot \laplace\p
=  \int_M  F \cdot \p + \int_{\pd M} G \cdot \p ,
$$
where $F\in W^{k-1,p}(M)$, $G\in W^{k,p}(M)$, and for some constant
$C$
\begin{align*}
\|F\|_{W^{k-1,p}} + \|G\|_{W^{k,p}}
\leq C \bigl( \|f\|_{W^{k,p}} + \|\g\|_{W^{k,p}} + \|\a\|_{W^{k,p}} \bigr) .
\end{align*}
So the regularity proposition~\ref{Laplace reg} for the weak Laplace equation with either 
Neumann ($\cT=\cC^\infty_\n(M)$) or Dirichlet ($\cT=\cC^\infty_\d(M)$) boundary conditions  
proves that $\a(X)\in W^{k+1,p}(M)$ with the according estimate.

In the case $k=0$ one works with the following inequality:
Let $\frac 1{p^*}+\frac 1p =1$, then there is a constant $C$ such that for all 
$\p\in\cT$
$$
\left| \int_M \a(X) \cdot \laplace\p \right|
\leq C \bigl( \|f\|_p + \|\g\|_p + \|\a\|_p \bigr) \| \p \|_{W^{1,p^*}} .
$$
(Under the assumptions of remark~\ref{Hodge rmk}, one simply replaces 
$\|f\|_p$ and $\|\g\|_p$ by $c_1$ and $c_2$ respectively.)
The regularity and estimate for $\a(X)$ then follow from proposition~\ref{Laplace reg}.
\QED \\

 \pagebreak
 
 \bibliographystyle{alpha}

\end{document}